\documentclass{article}

    \PassOptionsToPackage{numbers, compress}{natbib}

    \usepackage[preprint]{neurips_2025}

\usepackage[utf8]{inputenc} 
\usepackage[T1]{fontenc}    
\usepackage{hyperref}       
\usepackage{url}            
\usepackage{booktabs}       
\usepackage{amsfonts}       
\usepackage{nicefrac}       
\usepackage{microtype}      
\usepackage{xcolor}         

\usepackage{algorithm}
\usepackage{algorithmic}

\usepackage{amsmath}
\usepackage{amssymb}
\usepackage{mathtools}
\usepackage{amsthm}

\theoremstyle{plain}
\newtheorem{theorem}{Theorem}[section]

\newtheorem{lemma}[theorem]{Lemma}
\newtheorem{remark}[theorem]{Remark}

\theoremstyle{definition}

\newtheorem{assumption}[theorem]{Assumption}

\theoremstyle{remark}

\usepackage{amsfonts}
\usepackage{dsfont}
\usepackage{wrapfig}
\usepackage[symbol]{footmisc}

\allowdisplaybreaks

\usepackage{mdframed}

\usepackage{xcolor}
\usepackage{colortbl}
\usepackage{color}
\usepackage{graphicx}
\usepackage{multirow}
\usepackage{pifont}

\usepackage[font=small]{caption}
\usepackage[font=small]{subcaption}

\usepackage[algo2e]{algorithm2e} 
\usepackage{verbatim}
\usepackage{xspace} %

\usepackage{enumerate}
\usepackage{enumitem}

\newcommand{\dotprod}[2]{\left\langle #1,#2 \right\rangle}
\newcommand{\Obound}[1]{\mathcal{O}\left( #1 \right)}
\newcommand{\OboundTilde}[1]{\tilde{\mathcal{O}}\left( #1 \right)}

\newcommand{\norms}[1]{\left\| #1 \right\| }
\newcommand{\expect}[1]{\mathbb{E}\left[ #1 \right]}

\newcommand{\clip}[1]{ \text{clip}_c\left( #1 \right)}
\newcommand{\cliplam}[1]{ \text{clip}_\lambda\left( #1 \right)}

\newcommand{\prob}[1]{\mathcal{P}\left[ #1 \right]}

\definecolor{PineGreen}{HTML}{008B72}
\newcommand{\greencheck}{\color{PineGreen}\ding{51}}
\newcommand{\redx}{\color{red} \ding{55}}

  \providecommand{\prob}[1]{{\rm Pr}\left[#1\right] }

  \DeclareMathOperator*{\arginf}{arg\,inf}

  \providecommand{\bb}{\mathbf{b}}

  \providecommand{\ee}{\mathbf{e}}

  \renewcommand{\gg}{\mathbf{g}}

  \usepackage{bm}
  \newcommand{\bxi}{\boldsymbol{\xi}}

\usepackage[textsize=tiny]{todonotes}

\providecommand{\mycomment}[3]{\todo[caption={},size=footnotesize,color=#1!20, inline]{\textbf{#2: }#3}}%
\providecommand{\inlinecomment}[3]{%
  {\color{#1}#2: #3}}%
\newcommand\commenter[2]%
{%
  \expandafter\newcommand\csname i#1\endcsname[1]{\inlinecomment{#2}{#1}{##1}}
  \expandafter\newcommand\csname #1\endcsname[1]{\mycomment{#2}{#1}{##1}}
}

\newcommand{\circledOne}{\text{\ding{172}}}
\newcommand{\circledTwo}{\text{\ding{173}}}
\newcommand{\circledThree}{\text{\ding{174}}}
\newcommand{\circledFour}{\text{\ding{175}}}
\newcommand{\circledFive}{\text{\ding{176}}}
\newcommand{\circledSix}{\text{\ding{177}}}
\newcommand{\circledSeven}{\text{\ding{178}}}
\newcommand{\circledEight}{\text{\ding{179}}}
\newcommand{\circledNine}{\text{\ding{180}}}

\commenter{AL}{blue}

\usepackage{url}

\usepackage[many]{tcolorbox}    	

\definecolor{main}{HTML}{5989cf}    
\definecolor{sub}{HTML}{cde4ff}     

\tcbset{
    sharp corners,
    colback = white,
    before skip = 0.2cm,    
    after skip = 0.3cm      
}                           

\newtcolorbox{boxA}{
    fontupper = \bf,
    boxrule = 1.5pt,
    colframe = black 
}

\newtcolorbox{boxB}{
    fontupper = \bf\color{main}, 
    boxrule = 1.5pt,
    colframe = main,
    rounded corners,
    arc = 5pt   
}

\newtcolorbox{boxC}{
    colback = sub, 
    boxrule = 0pt  
}

\newtcolorbox{boxD}{
    colback = white, 
    colframe = black, 
    boxrule = 0pt, 
    toprule = 3pt, 
    bottomrule = 3pt 
}

\newtcolorbox{boxE}{
    enhanced, 
    boxrule = 0pt, 
    borderline = {0.75pt}{0pt}{main}, 
    borderline = {0.75pt}{2pt}{sub} 
}

\newtcolorbox{boxF}{
    colback = white,
    enhanced,
    boxrule = 1.5pt, 
    colframe = white, 
    borderline = {1.5pt}{0pt}{black, dashed} 
}

\newtcolorbox{boxG}{
    enhanced,
    boxrule = 0pt,
    colback = sub,
    borderline west = {1pt}{0pt}{main}, 
    borderline west = {0.75pt}{2pt}{main}, 
    borderline east = {1pt}{0pt}{main}, 
    borderline east = {0.75pt}{2pt}{main}
}

\newtcolorbox{boxH}{
    colback = white, 
    colframe = black, 
    boxrule = 0pt, 
    leftrule = 6pt 
}

\newtcolorbox{boxI}{
    colback = sub, 
    colframe = main, 
    boxrule = 0pt, 
    toprule = 6pt 
}

\newtcolorbox{boxJ}{
    sharpish corners, 
    colback = sub, 
    colframe = main, 
    boxrule = 0pt, 
    toprule = 4.5pt, 
    enhanced,
    fuzzy shadow = {0pt}{-2pt}{-0.5pt}{0.5pt}{black!35} 
}

\newtcolorbox{boxK}{
    sharpish corners, 
    boxrule = 0pt,
    toprule = 4.5pt, 
    enhanced,
    fuzzy shadow = {0pt}{-2pt}{-0.5pt}{0.5pt}{black!35} 
}

\newtcolorbox{boxL}{
    fontupper = \color{main},
    rounded corners,
    arc = 6pt,
    colback = sub, 
    colframe = main!50, 
    boxrule = 0pt, 
    bottomrule = 4.5pt 
}

\newtcolorbox{boxM}{
    rounded corners,
    arc = 6pt,
    colback = blue!5, 
    colframe = black, 
    boxrule = 0pt, 
    bottomrule = 1.5pt,
    toprule = 1.5pt,
    enhanced,
    fuzzy shadow = {0pt}{-3pt}{-0.5pt}{0.5pt}{black!35}
}

\title{Power of Generalized Smoothness in Stochastic Convex Optimization: First- and Zero-Order Algorithms}

\author{%
  Aleksandr Lobanov \\
  MIPT, Skoltech, HSE\\
  \texttt{lobbsasha@mail.ru} \\
  \And
  Alexander Gasnikov \\
  Innopolis, MIPT, ISP RAS \\
  \texttt{gasnikov@yandex.ru} \\
}

\begin{document}

\maketitle

\begin{abstract}
  This paper is devoted to the study of stochastic optimization problems under the generalized smoothness assumption. By considering the unbiased gradient oracle in \textit{Stochastic Gradient Descent}, we provide strategies to achieve in bounds the summands describing linear rate. In particular, in the case $L_0 = 0$, we obtain in the \textbf{convex setup} the iteration complexity: $N = \Obound{L_1R \log\frac{1}{\varepsilon} + \frac{L_1 c R^2}{\varepsilon}}$ for \textit{Clipped Stochastic Gradient Descent} and ${N = \Obound{L_1R \log\frac{1}{\varepsilon}}}$ for \textit{Normalized Stochastic Gradient Descent}. Furthermore, we generalize the convergence results to the case with a biased gradient oracle, and show that the power of $(L_0,L_1)$-smoothness extends to \textit{zero-order algorithms}. Finally, we demonstrate the possibility of linear convergence in the convex setup through numerical experimentation, which has aroused some interest in the machine learning community.
\end{abstract}

\section{Introduction}\label{sec:Introduction}
In many real-world scenarios, systems are often noisy and complex, making deterministic optimization infeasible. Therefore, this work focuses on a stochastic optimization problem:
\begin{equation}
    \label{eq:init_problem}
    \min_{x \in \mathbb{R}^d} \left\{f(x) := \mathbb{E}_{\xi \sim \mathcal{D}} \left[f(x,\xi) \right] \right\},
\end{equation}
where $f: \mathbb{R}^d \rightarrow \mathbb{R}$ is a convex function and where we assume that optimization algorithms only have access to the gradient oracle $\gg: \mathbb{R}^d \times \mathcal{D} \rightarrow \mathbb{R}^d$ with stochastic gradient $\expect{\nabla f(x,\xi)} = \nabla f(x)$ and bias $\bb(x)$ terms:
\begin{equation}
    \label{eq:gradient_oracle}
    \gg(x, \xi) = \nabla f(x,\xi) + \bb(x).
\end{equation}

Frequently, to solve problem~\eqref{eq:init_problem} one uses what is likely already a classic optimization algorithm, namely Stochastic Gradient Descent (SGD) \cite{Bottou_1998} or its variations, which have demonstrated their effectiveness in different settings, for instance, federated learning~\cite{Yuan_2020, Kairouz_2021,Woodworth_2021}, deep learning \cite{Dean_2012,Zhang_2015,Dimlioglu_2024}, reinforcement learning \cite{Bello_2017,Lee_2024} and others. Among the variants of SGD, it is worth noting the Normalized Stochastic Gradient Descent (NSGD) \cite{Hazan_2015,Zhao_2024} which has received widely attention from the community because it addresses challenges in optimization for machine learning \cite{Bengio_1994}. And it's also worth noting the Clipped Stochastic Gradient Descent (ClipSGD) \cite{Goodfellow_2016}, which is commonly used to stabilize the training of deep learning models \cite{Pascanu_2013,Gorbunov_2020}.

Many standard literatures analyze stochastic optimization algorithms with unbiased gradient oracle~\eqref{eq:gradient_oracle}. In particular, SGD \cite{Lacoste_2012,Bottou_2018}, NSGD \cite{Zhao_2021,Hubler_2024}, ClipSGD \cite{Gorbunov_2020,Koloskova_2023}. However, there are a number of applications where gradient oracle~\eqref{eq:gradient_oracle} is biased. For example, sparsified SGD \cite{Alistarh_2018}, delayed SGD \cite{Stich_2019}, etc. Zero-order algorithms \cite{Nesterov_2017,Demidovich_2023} occupy a special place in the class of stochastic methods with biased gradient oracle~\eqref{eq:gradient_oracle}. They are motivated by various applications, including multi-armed bandit \cite{Shamir_2017, Lattimore_2021}, online optimization \cite{Agarwal_2010, Bach_2016, Tsybakov_2022}, hyperparameter tuning \cite{Hernandez_2014, Nguyen_2022}.

In our work, we investigate the convergence of first-order algorithms: ClipSGD, NSGD, and zero-order algorithms: ZO-ClipSGD, ZO-NSGD, assuming convexity and $(L_0,L_1)$-smoothness.

We emphasize the following points:

\textbf{Algorithm step size.} 
Zero-order algorithms do not have access to the exact (stochastic) gradient in particular, as well as every algorithm with a biased gradient oracle, so we focus on creating first-order methods whose step size does not depend on knowledge of the gradient at a given point. We use the developed first-order algorithms as a basis for creating zero-order methods \cite{Gasnikov_2023}. 

\begin{wrapfigure}{r}{0.45\textwidth}
  \vskip -0.20in
  \centering
  \includegraphics[width=0.45\textwidth]{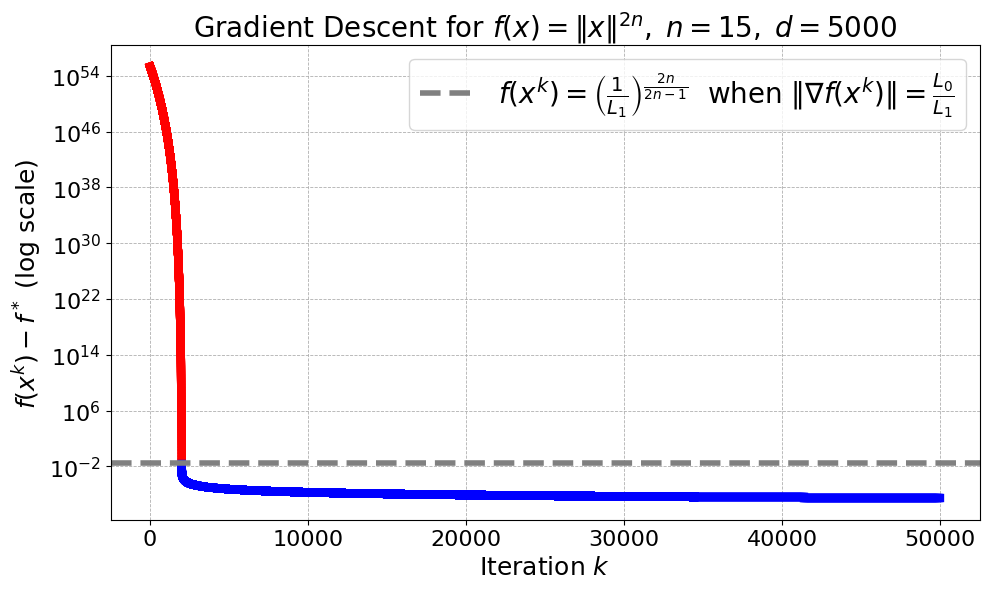}
  \caption{Changing regimes demonstration}
  \label{fig:Change_regimes}
  \vskip -0.25in
\end{wrapfigure}

\textbf{Linear convergence.} Historically \cite{Nesterov_2018}, stochastic optimization first- and zero-order algorithms have achieved the desired accuracy with a linear rate of convergence only in strongly convex case and under assumption of standard smoothness. However, the work of \cite{Lobanov_2024} showed that if the generalized smoothness assumption is satisfied in a \underline{\textit{deterministic convex}} optimization problem, then gradient descent has two regimes: linear convergence rate as long as $\norms{\nabla f(x^k)} \geq \frac{L_0}{L_1}$, and a sublinear convergence rate in the other case (see the example of the power of norm function in Figure~\ref{fig:Change_regimes}). Considering these  points, our work answers the following question:

\begin{center}
    \textit{Can linear convergence rate in stochastic convex optimization be achieved for first- and zero-order algorithms with constant step size?}
\end{center}

\subsection{Main Contributions}
    More specifically, our contributions are the following:
    \begin{itemize}
        \item We provide strategies to obtain summands that describe the linear convergence~rate.~In~particular, we show that using clipping or normalization techniques can achieve~the~desired~results.

        \item We improve convergence results for ClipSGD and NSGD with unbiased gradient oracle~\eqref{eq:gradient_oracle} in the convex setting assuming $(L_0,L_1)$-smoothness (see Table~\ref{tab:table_compare}). Moreover, we show that in the case $L_0 = 0$, NSGD can converge in the convex setup with a linear convergence rate to the desired accuracy, requiring $N = \OboundTilde{L_1 R}$ iterations and~$B = \Obound{\frac{\sigma^2 M R^3}{\varepsilon^{3}}}$~batch~size.

        \item We generalize ClipSGD, and NSGD to the case of a biased gradient oracle, showing how the bias accumulates over iterations (this result may be of independent interest).

        \item We provide the first convergence results for the zero-order algorithms ZO-ClipSGD (Algorithm~\ref{algo:ZO_ClipSGD}), and ZO-NSGD (Algorithm~\ref{algo:ZO_NSGD}) in the convex and $(L_0,L_1)$-smooth setting. We show that the power of generalized smoothness extends to zero-order methods as well, achieving summands characterizing the linear convergence rate (see Table~\ref{tab:table_compare}).

        \item We demonstrate on a numerical example of logistic regression (which is of particular interest to the machine learning community) that indeed, zero- and first-order stochastic algorithms can converge with linear rates in a convex setup.
    \end{itemize}

\subsection{Formal Setting and Assumptions}
    In this subsection, we introduce and discuss main assumptions and notations used~throughout~paper.

    \vspace{-1em}
    \paragraph{Notations.} We use $\dotprod{x}{y}:= \sum_{i=1}^{d} x_i y_i$ to denote standard inner product of $x,y \in \mathbb{R}^d$. We denote Euclidean norm in $\mathbb{R}^d$ as $\| x\| := \sqrt{ \sum_{i=1}^d x_i^2}$. In particular, this norm ${\| x \| := \sqrt{\dotprod{x}{x}}}$ is related to the inner product. We use $\mathcal{P}[\cdot]$ to define probability measure which is always known from the context, $\mathbb{E}[\cdot]$ denotes mathematical expectation. We use the following notation $B^d(r):=\left\{ x \in \mathbb{R}^d : \| x \| \leq r \right\}$ to denote Euclidean ball ($l_2$-ball) and  $S^d(r):=\left\{ x \in \mathbb{R}^d : \| x \| = r \right\}$ to denote Euclidean sphere. We denote $0\leq M < \infty$ as the upper bound of the gradient norm $\norms{\nabla f(x^k)}$. For simplicity, we denote $f^* := f(x^*)$ and $R = \norms{x^0 - x^*}$. We use $\tilde{O} (\cdot)$ to hide the logarithmic coefficients.

    \begin{table*}
    \begin{minipage}{\textwidth}
        \caption{Comparison of convergence results of SGD variants to the most related work \cite{Gaash_2025} in the convex and $(L_0, L_1)$-smooth setup. Notation: $\eta \leq (L_0 + L_1 c)^{-1}$ -- step size; $c>0$ -- clipping radius;  $\mathcal{R} = \left(\eta + \frac{MR}{c^2} + \frac{R}{c}  \right)$; $\varepsilon =$ desired accuracy; $d =$ dimension; SLCR $=$ summand with linear convergence rate.}
    \label{tab:table_compare} 
    \centering
    \resizebox{\linewidth}{!}{ 
    \begin{tabular}{lcccccc}\toprule
    \multirow{2}{*}{Algorithm} & Number of Iterations  & Batch Size & Maximum Noise Level & \multirow{2}{*}{SLCR?} & \multirow{2}{*}{Reference} \\
    & $\#N$ & $\#B$ & $\#\Delta$ & & & \\ \midrule
    \multirow{3}{*}{ClipSGD} & $\Obound{\frac{L_0 R^2}{\varepsilon} + \frac{\sigma^2 R^2}{\varepsilon^2}  + L_1^2 R^2}$ & \redx & \redx & \redx & \citet{Gaash_2025}\\
     & $\Obound{\frac{R}{\eta c} \log \frac{1}{\varepsilon} + \frac{R^2}{\eta \varepsilon}}$ & $\Obound{\frac{\sigma^2 \mathcal{R}}{\varepsilon}}$ & \redx & \greencheck & Theorem~\ref{th:clipSGD} \textbf{(Ours)}\\
     NSGD & $\Obound{(L_1 R + \frac{L_0 R^2}{\varepsilon}) \log \frac{1}{\varepsilon}}$ & $\Obound{\frac{\sigma^2 M R^3}{\varepsilon^3}}$ & \redx & \greencheck & Theorem~\ref{th:NSGD} \textbf{(Ours)}\\ \midrule
     ZO-ClipSGD & $\Obound{\frac{R}{\eta c} \log \frac{1}{\varepsilon} + \frac{R^2}{\eta \varepsilon}}$ & $\Obound{\frac{d M R \tilde{\sigma}^2}{\varepsilon c^2}}$ & $\Obound{\frac{\varepsilon}{\sqrt{d} R (L_0 + L_1 M)}\min \left\{ \tilde{\sigma}, \frac{\varepsilon}{\sqrt{d} R}  \right\}}$ & \greencheck & Theorem~\ref{th:ZO_ClipSGD} \textbf{(Ours)}\\
     ZO-NSGD & $\Obound{(L_1 R + \frac{L_0 R^2}{\varepsilon}) \log \frac{1}{\varepsilon}}$ & $\Obound{\frac{d M R^3 \tilde{\sigma}^2}{\varepsilon^3}}$ & $\Obound{\frac{\varepsilon^{3/2}}{\sqrt{d} R^{3/2} (L_0 + L_1 M)}\min \left\{ \tilde{\sigma}, \frac{\varepsilon^{3/2}}{\sqrt{d} R^{3/2}}  \right\}}$ & \greencheck & Theorem~\ref{th:ZO_NSGD} \textbf{(Ours)}\\
    \bottomrule
    \end{tabular}}
    \end{minipage}
\end{table*}

    \vspace{-1em}
    \paragraph{Assumptions on objective function.} Throughout this paper, we refer to the standard $L$-smoothness assumption, which is widely used in the literature \citep[e.g.][]{Polyak_1987} and has the following form:
    \begin{assumption}[$L$-smoothness]\label{ass:L_smooth}
        Function $f$ is $L$-smooth if for any $x,y \in \mathbb{R}^d$ is satisfied:
        \begin{equation*}
            \norms{\nabla f(y) -\nabla f(x)} \leq L \norms{y - x}.
        \end{equation*}
    \end{assumption}

    Despite the widespread use of Assumption~\ref{ass:L_smooth}, our work focuses on the more general smoothness assumption, which has recently attracted increased interest. In particular, in \cite{Zhang_2019} it was shown that norm of Hesse matrix correlates with norm of gradient function when training neural networks, and in \cite{Lobanov_2024} it was shown that using generalized smoothness it is possible to significantly improve the convergence of algorithms. $(L_0,L_1)$-smoothness \cite{Zhang_2019,Zhang_2020} has been proposed as a natural relaxation of standard~smoothness~assumption.

    \begin{assumption}[$(L_0, L_1)$-smoothness]\label{ass:L0_L1_smooth}
        A function ${f: \mathbb{R}^d \rightarrow \mathbb{R}}$ is $(L_0, L_1)$-smooth if the following inequality is satisfied for any $x,y \in \mathbb{R}^d$ with $\norms{y - x} \leq \frac{1}{L_1}$:
        \begin{equation*}
            \norms{\nabla f(y) -\nabla f(x)} \leq \left(L_0 + L_1 \norms{\nabla f(x)} \right) \norms{y - x}.
        \end{equation*}
    \end{assumption}
    Assumption~\ref{ass:L0_L1_smooth} in the case $L_1 =0$ covers the standard Assumption~\ref{ass:L_smooth}. Moreover, $(L_0,L_1)$-smoothness is strictly more general than $L$-smoothness, see the examples in \cite{Zhang_2019, Chen_2023, Koloskova_2023, Gorbunov_2024}.

    \begin{remark}[Clarification regarding $L_0 = 0$]\label{rem:Smoothness_only_L0}
        In this paper we often emphasize the case $L_0=0$ in Assumption~\ref{ass:L0_L1_smooth}. It is worth noting that the class of functions that do not reach their infimum $x^*$ (converge to an asymptote) satisfies this case. Explicit examples of functions with $L_0 = 0$ are the exponent of the inner product and the logistic function (see \cite{Gorbunov_2024} for details). 
    \end{remark}

    \vspace{-1em}
    \paragraph{Assumptions on gradient oracle.}
    In our analysis, we consider cases with both unbiased and biased gradient oracle~\eqref{eq:gradient_oracle}. Therefore, we assume that the bias and variance of gradient oracle~\eqref{eq:gradient_oracle}~are~bounded:

    \begin{assumption}[Bounded bias] \label{ass:bounded_bias}
        There exists constant ${\zeta \geq 0}$ such that the bias is bounded if $\forall x \in \mathbb{R}^d$:
        \begin{equation*}
            \norms{\bb(x)} \leq \zeta.
        \end{equation*}
    \end{assumption}

    \begin{assumption}[Bounded variance] \label{ass:bounded_var}
        There exists constant $\sigma^2 \geq 0$ such that the variance is bounded if $\forall x \in \mathbb{R}^d$:
        \begin{equation*}
            \expect{\norms{\gg(x,\xi) - \expect{\gg(x,\xi)}}^2} \leq \sigma^2.
        \end{equation*}
    \end{assumption}
    Assumption~\ref{ass:bounded_bias} is organic \citep[see, e.g.][]{Lobanov_JOTA}, and the case $\zeta = 0$ corresponds to the unbiased gradient oracle~\eqref{eq:gradient_oracle}. Assumption~\ref{ass:bounded_var} is often used by the community \citep[e.g.][]{Juditsky_2010,Lan_2012}, and is sometimes called heavy-tailed noise \cite{Gorbunov_2020}.

\subsection{Paper Organization}
    Next, our paper has the following structure. In Section~\ref{sec:Related Works}, we discuss related work. In Section~\ref{sec:Clipped Stochastic Gradient Descent}, We start to present the main results of our work, in particular, we provide the first strategy for obtaining a summand characterizing the linear rate in the convergence estimate. In Section~\ref{sec:Normalized Stochastic Gradient Descent}, we analyze NSGD in the convex setting, showing in which regime linear convergence can be observed. We provide a first analysis of zero-order algorithms under $(L_0,L_1)$-smoothness in Section~\ref{sec:Zero-Order Algorithms}. In Section~\ref{sec:Discussion and Future Works}, we discuss the results obtained. While, in Section~\ref{sec:Numerical Experiments}, we show experimentally about the possibility of linear convergence in the convex setting. Finally, Section~\ref{sec:Conclusion} concludes our paper. All missing proofs of Lemmas and Theorems are provided in the supplementary materials (Appendix).

\section{Related Works}\label{sec:Related Works}

In this section, we will discuss the most related works.
\vspace{-1em}
\paragraph{Algorithms under $(L_0,L_1)$-smoothness.} Generalized smoothness was first introduced in \cite{Zhang_2019}, which analyzed ClipSGD in the non-convex setting. A number of works \cite{Wang_2023,Li_2023,Faw_2023,Li_2023_Convex,Hong_2024,Xie_2024,Wang_2024} followed that also focused on the non-convex setup, including ClipSGD \cite{Zhang_2020,Zhang_2020_Improved,Koloskova_2023}, NSGD \cite{Zhao_2021,Hubler_2024_NSGD}. After that, there was interest in research on algorithms in the convex deterministic setting: Clipped Gradient Descent \cite{Koloskova_2023}, Normalized Gradient Descent \cite{Chen_2023_NSGD}, Gradient Descent with Polyak step size $\eta_k = \frac{f(x^k) - f^*}{\| \nabla f(x^k) \|^2}$ \cite{Takezawa_2024}, and $\eta_k = \frac{1}{L_0 + L_1 \| \nabla f(x^k) \|}$ \cite{Gorbunov_2024,Vankov_2024}. Moreover, in \cite{Lobanov_2024}, it was theoretically shown that it is possible to significantly improve the convergence of algorithms in the (strongly) convex setting by achieving linear convergence rate. However, much less attention has been paid to the stochastic convex setting. Perhaps the only results are \cite{Gorbunov_2024,Gaash_2025}, which considers SGD and ClipSGD achieving only a sublinear convergence rate. Moreover, in the case $L_0 = 0$ in the Assumption, the algorithms from \cite{Gorbunov_2024,Gaash_2025} cannot converge to the desired accuracy. \textit{In our work, we focus on the stochastic convex setup, showing that existing convergence results can be significantly improved.}   
\vspace{-1em}
\paragraph{Zero-order algorithms.} The work of \cite{Gasnikov_2022} showed that to achieve optimal estimates of iteration $N$ and oracle $T$ complexity in zero-order algorithms, one should base it on a first-order algorithm using a gradient approximation as the biased gradient oracle~\eqref{eq:gradient_oracle}, which uses only information about the objective function $f$. Using this technique a number of works have achieved the best convergence results in various settings including distributed optmization \cite{Akhavan_2021}, federated optimization \cite{Patel_2022}, overparameterization \cite{Lobanov_2023}, Polyak-Lojasiewicz condition \cite{Gasnikov_2024}, etc. However, all these works assumed standard smoothness (Assumption~\ref{ass:L_smooth}) and achieved only sublinear convergence rates. \textit{In our work, we present convergence results for zero-order algorithms under $(L_0,L_1)$-smoothness.}   
\vspace{-0.5em}
\section{Clipped Stochastic Gradient Descent}\label{sec:Clipped Stochastic Gradient Descent}
\vspace{-0.5em}
In this section we begin to present the main results of our work. In particular, we analyze the convergence of SGD variants under convexity and $(L_0,L_1)$-smoothness with step size independent of the gradient norm. We assume that the gradient oracle~\eqref{eq:gradient_oracle} is unbiased $\zeta = 0$, i.e.,~Assumption~\ref{ass:bounded_var}~takes:
\begin{equation*}
    \expect{\norms{\nabla f(x,\xi) - \nabla f(x)}^2} \leq \sigma^2.
\end{equation*}

As a first strategy to obtain the summands that characterize the linear rate, we consider the clipping technique. Applying this technique we produce the ClipSGD, which has the following form:

\begin{algorithm}[H]
    \caption{Clipped Stochastic Gradient Descent Method (ClipSGD)}
    \label{algo:clipSGD}
    \begin{algorithmic}
        \STATE {\bfseries Input:} initial point $x_0 \in \mathbb{R}^d$, iterations $N$, batch size $B$, step size $\eta_k >0$~and~clipping~radius~$c>0$
        \vspace{-1em}
        \FOR{$k=0$ {\bfseries to} $N-1$}
            \STATE  {~~~1.} Draw fresh i.i.d. samples $\xi^k_1,...,\xi^k_B$
            \STATE  {~~~2.} $\nabla f(x^k, \bxi^k) = \frac{1}{B} \sum_{i=1}^B \nabla f(x^k, \xi^k_i)$ 
            \STATE  {~~~3.} ${\text{clip}_c (\nabla f(x^k, \bxi^k)) = \min \left\{ 1, \frac{c}{\|\nabla f(x^k, \bxi^k)\|} \right\} \nabla f(x^k, \bxi^k)}$
            \STATE  {~~~4.} $x^{k+1} \gets x^{k}  - \eta_{k} \cdot \text{clip}_c (\nabla f(x^k, \bxi^k))$
        \ENDFOR
        \STATE {\bfseries Return:} $x^{N}$
    \end{algorithmic}
\end{algorithm}  
\vskip -0.2in
Algorithm~\ref{algo:clipSGD} uses the clipped stochastic gradient $\text{clip}_c (\nabla f(x, \bxi))$, which normalizes the gradient only if $\norms{\nabla f(x,\bxi)} > c$. Next theorem provides the convergence result for ClipSGD.

\begin{boxM}
    \begin{theorem}\label{th:clipSGD}
        Let function $f$ satisfy Assumption~\ref{ass:L0_L1_smooth} ($(L_0,L_1)$-smoothness) and unbiased gradient oracle~\eqref{eq:gradient_oracle} satisfy Assumption~\ref{ass:bounded_var} (bounded variance), then Algorithm~\ref{algo:clipSGD} with constant step size ${\eta_k = \eta \leq \left[4 (L_0 + L_1 c)\right]^{-1}}$ and  arbitrary clipping radius $c$ guarantees~error:
        \begin{align*}
            \expect{f(x^N)} - f^* \lesssim  \left( 1 - \frac{\eta c}{R} \right)^{K} &(f(x^0) - f^*) +  \frac{R^2}{\eta (N - K)} + \frac{\sigma^2}{B} \left(\eta + \frac{MR}{c^2} + \frac{R}{c}  \right),
        \end{align*}
        
        where $0 \leq K < N$ is number of iterations for which $\|\nabla f(x^k)\| \leq \frac{c}{2}$ is satisfied.
    \end{theorem}
\end{boxM}

It should be noted that the results of Theorem~\ref{th:clipSGD} are given with a choice of step size independent of the gradient at the current point. This choice of step allows us to separate the constants $L_0$ and $L_1$ in the final estimates. The summand $\frac{R^2}{\eta (N - K)}$ is a typical ClipSGD characterizing the sublinear rate (see e.g. \cite{Gorbunov_2020}). However, it is worth noting that by substituting $\eta = (L_0 + L_1 c)^{-1}$, then the summand with $L_1: \frac{L_1 c R^2}{N-K}$ already improves existing results both assuming standard smoothness \cite{Gorbunov_2020} and generalized smoothness \cite{Gaash_2025}. The first summand $\left( 1 - \frac{\eta c}{R} \right)^{K} (f(x^0) - f^*)$, which characterizes the linear rate, deserves special attention. \textit{To the best of our knowledge, this is the first result for ClipSGD showing such a summand in a convex setting.} Moreover, by substituting $\eta = (L_0 + L_1 c)^{-1}$, it is not hard to see that in the regime $L_0 = 0$ (see Remark~\ref{rem:Smoothness_only_L0}), Algorithm~\ref{algo:clipSGD} at a batch size $B = \Obound{\frac{\sigma^2}{\varepsilon} \left(\eta + \frac{MR}{c^2} + \frac{R}{c} \right)}$ requires only $N = \Obound{L_1 R \log \frac{1}{\varepsilon} + \frac{L_1 c R^2}{\varepsilon}}$ iterations. This iteration complexity significantly outperforms standard results in the $L$-smoothness setting (Assumption~\ref{ass:L_smooth}), since \cite{Lan_2012} shows a lower bound consisting only of a sublinear convergence rate. Moreover, comparing to the closest work to the problem setting, then even when $\sigma = 0$ \cite{Gaash_2025} does not guarantee convergence to the desired accuracy, offering an estimate of $N= \Obound{L_1^2 R^2}$ that is independent of accuracy . 

\section{Normalized Stochastic Gradient Descent}\label{sec:Normalized Stochastic Gradient Descent}
In the previous section, we showed that it is possible to obtain in the final convergence estimate a summand characterizing the linear rate. In addition, we highlighted the regime $L_0 = 0$, in which ClipSGD has the following iteration complexity: $N = \Obound{L_1 R \log \frac{1}{\varepsilon} + \frac{L_1 c R^2}{\varepsilon}}$. However, with this iteration complexity, it cannot be said that the algorithm can converge with linear rate to the desired accuracy. Such an estimate can characterize that the algorithm converges with linear rate as long as the gradient norm is large $\norms{\nabla f(x^k)} \geq c$, and then slows down to the sublinear rate. However, it is worth noting that the summand responsible for the sublinear rate depends on the clipping radius: $\frac{L_1 c R^2}{\varepsilon}$. That is, if we take $c$ smaller, the ClipSGD will take longer to converge to the linear rate. Thus, noticing that the regime $\norms{\nabla f(x^k)} \geq c$ is a gradient normalization, then considering NSGD, it seems that one can achieve a true linear convergence rate to the desired accuracy. 

In this section, we consider a normalization technique to obtain the summand characterizing the linear rate. Applying this technique we produce the NSGD, which has the following form (see~Algorithm~\ref{algo:NSGD}):

\begin{algorithm}[H]
        \caption{Normalized Stochastic Gradient Descent Method (NSGD)}
        \label{algo:NSGD}
        \begin{algorithmic}
            \STATE {\bfseries Input:} initial point $x_0 \in \mathbb{R}^d$, iterations number $N$, batch size $B$ and step size $\eta_k >0$ 
            \FOR{$k=0$ {\bfseries to} $N-1$}
                \STATE {~~~1.} Draw fresh i.i.d. samples $\xi^k_1,...,\xi^k_B$
                \STATE  {~~~2.} $\nabla f(x^k, \bxi^k) = \frac{1}{B} \sum_{i=1}^B \nabla f(x^k, \xi^k_i)$ 
                \STATE {~~~3.} $x^{k+1} \gets x^{k}  - \eta_{k} \cdot \frac{\nabla f(x^k, \bxi^k)}{\norms{\nabla f(x^k, \bxi^k)}}$
            \ENDFOR
            \STATE {\bfseries Return:} $x^{N}$
        \end{algorithmic}
    \end{algorithm}  
    \vskip -0.2in
    The following theorem provides a convergence result for Algorithm~\ref{algo:NSGD}.

    \begin{boxM}
        \begin{theorem}\label{th:NSGD}
            Let function $f$ satisfy Assumption~\ref{ass:L0_L1_smooth} ($(L_0,L_1)$-smoothness) and unbiased gradient oracle~\eqref{eq:gradient_oracle} satisfy Assumption~\ref{ass:bounded_var} (bounded variance), then Algorithm~\ref{algo:NSGD} with hyperparameter $\lambda > 0$ and constant step size ${\eta_k = \eta \leq \lambda / \left[2 (L_0 + L_1 \lambda)\right]}$ guarantees:
            \begin{align*}
                \expect{f(x^{N})} - f^* \lesssim  \left( 1 - \frac{\eta}{R} \right)^{N} &(f(x^0) - f^*) + \frac{\sigma^2 M R}{B \lambda^2} + \lambda R.
            \end{align*}
        \end{theorem}
    \end{boxM}

    From Theorem~\ref{th:NSGD} we can see that we have indeed got rid of the summand characterizing the sublinear rate from the deterministic part. \textit{Thus, we see that it is normalization that allows us to achieve the summand characterizing the linear rate} $\left( 1 - \frac{\eta}{R} \right)^{N} (f(x^0) - f^*)$. However, note that by substituting $\eta = \lambda / \left[2 (L_0 + L_1 \lambda)\right]$, we obtain a summand with $L_0: \left( 1 - \frac{\lambda}{R L_0} \right)^{N} (f(x^0) - f^*)$, which is in fact sublinear since it depends on the hyperparameter $\lambda$ (it follows from the third summand that $\lambda \sim \varepsilon/R$), and with $L_1: \left( 1 - \frac{1}{R L_1} \right)^{N} (f(x^0) - f^*)$, which is indeed linear since it does not depend on $\lambda$ in any way. That is, Algorithm~\ref{algo:NSGD}, which uses batch parallelization, requires $N = \Obound{(L_1 R + \frac{L_0 R^2}{\varepsilon}) \log \frac{1}{\varepsilon}}$ iterations. Similar to the reasoning in the previous section, it is worth highlighting the regime $L_0 = 0$ (see Remark~\ref{rem:Smoothness_only_L0}). Then we obtain a very surprising result on iteration complexity, namely, to achieve the desired accuracy NSGD converge with a linear rate of $N = \Obound{L_1 R \log \frac{1}{\varepsilon}}$ iterations. \textit{This estimate \underline{breaks all existing bounds on first-order algorithms \cite{Lan_2012}}, given the specificity of the problem formulation, namely convexity}. However, to achieve this rate over iterations, NSGD requires a batch size $B = \Obound{\frac{\sigma^2 M R^3}{\varepsilon^3}}$. We emphasize that the fact that NSGD requires a large batch size is not surprising (see, e.g., \cite{Cutkosky_2020}), in contrast to the true linear convergence rate in the convex setup.

\section{Zero-Order Algorithms}\label{sec:Zero-Order Algorithms}
In this section, we consider another class of algorithms: optimization algorithms that have access only to an objective function value $f$ possibly with some bounded adversarial noise $|\delta(x)| \leq \Delta$:
\begin{equation}
    \label{eq:zero_oracle}
    \tilde{f}(x,\xi) = f(x, \xi) + \delta(x).
\end{equation}
In \eqref{eq:zero_oracle}, $\Delta$ means the maximum possible allowable noise level at which the desired accuracy can still be achieved. In \cite{Anonymous_2025}, the importance of considering $\Delta$ as a third optimality criterion for zero-order algorithms was shown. In particular, in some applications \cite{Bogolubsky_2016}, the larger noise level $\Delta$ is, the cheaper the call to the inexact oracle $\tilde{f}$ in \eqref{eq:zero_oracle}.

Since this class of algorithms does not have access to the stochastic gradient $\nabla f(x,\xi)$, the gradient oracle~\eqref{eq:gradient_oracle} will be the gradient approximation with $L_2$ randomization \cite{Shamir_2017,Nesterov_2017}:
\begin{equation}
    \label{eq:approximation_gradient}
    \gg(x, \{e, \xi\}) = \frac{d}{2 \gamma} \left( \tilde{f}(x + \gamma e, \xi) - \tilde{f}(x - \gamma e, \xi) \right) e,
\end{equation}
where $\gamma > 0$ is a smoothing parameter, $e$ is a random vector uniformly distributed in $S^d(1)$.  

Due to the fact that the gradient approximation is the biased gradient oracle~\eqref{eq:gradient_oracle}, in order to create zero-order algorithms by basing on the results in Sections~\ref{sec:Clipped Stochastic Gradient Descent} and \ref{sec:Normalized Stochastic Gradient Descent}, it is necessary to first generalize the results of Theorems~\ref{th:clipSGD}, \ref{th:NSGD} (note that in these regimes it is not necessary to know the (stochastic) norm of the gradient with step size) to the case of gradient~oracle~with~bias. 

Next, we present convergence results for the following two zero-order algorithms: ZO-ClipSGD (Algorithm~\ref{algo:ZO_ClipSGD}) and ZO-NSGD (Algorithm~\ref{algo:ZO_NSGD}).

    \subsection{ZO-ClipSGD Method}
    The first algorithm we consider in this section is ZO-ClipSGD. This algorithm is a modification of ClipSGD (Algorithm~\ref{algo:clipSGD}), which uses instead of the original $\norms{\nabla f(x,\xi)}$, the stochastic gradient approximation~\eqref{eq:approximation_gradient}, which is the biased gradient oracle~\eqref{eq:gradient_oracle}. The ZO-ClipSGD has the following form:
    \begin{algorithm}[H]
        \caption{Zero-Order Clipped Stochastic Gradient Descent Method (ZO-ClipSGD)}
        \label{algo:ZO_ClipSGD}
        \begin{algorithmic}
            \STATE {\bfseries Input:} initial point $x_0 \in \mathbb{R}^d$, iterations $N$, batch size $B$, step size $\eta_k >0$ and clipping~radius~$c>0$  
            \vspace{-1em}

            \FOR{$k=0$ {\bfseries to} $N-1$}
                \STATE {~~~1.} Draw fresh i.i.d. samples $\xi^k_1,...,\xi^k_B$ and $e^k_1,...,e^k_B$
                \vspace{0.2em}
                \STATE {~~~2.} $\gg(x^k, \{\ee^k, \bxi^k\}) = \frac{1}{B} \sum_{i=1}^B \gg(x^k, \{e^k_i, \xi^k_i\})$ via~\eqref{eq:approximation_gradient}
                \vspace{0.4em}
                \STATE {~~~3.} ${\text{clip}_c (\gg(x^k, \{\ee^k, \bxi^k\})) = \min \left\{ 1, \frac{c}{\|\gg(x^k, \{\ee^k, \bxi^k\})\|} \right\} \gg(x^k, \{\ee^k, \bxi^k\})}$
                \vspace{0.1em}
                \STATE {~~~4.} $x^{k+1} \gets x^{k}  - \eta_{k} \cdot \text{clip}_c (\gg(x^k, \{\ee^k, \bxi^k\}))$
            \ENDFOR
            \STATE {\bfseries Return:} $x^{N}$
        \end{algorithmic}
    \end{algorithm}  
    \vskip -0.2in
    Before proceeding to present the convergence results of Algorithm~\ref{algo:ZO_ClipSGD}, we note that the gradient approximation~\eqref{eq:approximation_gradient} is a biased gradient oracle, so we cannot directly use the results obtained in Theorem~\ref{th:clipSGD}. Thus, in order to obtain estimates for the iteration complexity $N$, oracle complexity $T$ and maximum noise $\Delta$, we first need to generalize the results of Theorem~\ref{th:clipSGD} to the case with a biased gradient oracle~\eqref{eq:gradient_oracle}. 
    \begin{lemma}
        \label{lem:biased_ZO_ClipSGD}
       Let function $f$ satisfy Assumption~\ref{ass:L0_L1_smooth} and \underline{biased gradient oracle}~\eqref{eq:gradient_oracle} (${\zeta > 0}$) satisfy Assumption~\ref{ass:bounded_var}, then Algorithm~\ref{algo:clipSGD} with step size ${\eta_k = \eta \leq \left[4 (L_0 + L_1 c)\right]^{-1}}$ guarantees~error:
       \begin{equation*}
             \expect{f(x^N)} - f^*\lesssim  \left( 1 - \frac{\eta c}{R} \right)^{K} (f(x^0) - f^*) +  \frac{R^2}{\eta (N - K)} +  \mathcal{R} \cdot\left( \frac{\sigma^2}{B} + \zeta^2 \right) + R \zeta,
       \end{equation*}
        where $c$ is arbitrary clipping radius, $\mathcal{R} = \left(\eta + \frac{MR}{c^2} + \frac{R}{c}  \right)$, $0 \leq K < N$ is the number of iterations for which the condition $\|\nabla f(x^k)\| \leq \frac{c}{3}$ is satisfied. 
    \end{lemma}
    Lemma~\ref{lem:biased_ZO_ClipSGD} shows how the bias accumulates over iterations, thus converging to the error floor. By estimating the second moment and the bias of the gradient approximation~\eqref{eq:approximation_gradient} and substituting them into Lemma~\ref{lem:biased_ZO_ClipSGD}, we find three optimality criteria for ZO-ClipSGD. 
    \begin{boxM}
        \begin{theorem}\label{th:ZO_ClipSGD}
            Let function $f$ satisfy Assumption~\ref{ass:L0_L1_smooth}, gradient approximation~\eqref{eq:approximation_gradient} satisfy Assumption~\ref{ass:bounded_var}, then Algorithm~\ref{algo:ZO_ClipSGD} with step size ${\eta_k = \eta \leq \left[4 (L_0 + L_1 c)\right]^{-1}}$ converges to desired $\varepsilon$ accuracy after:
            \begin{equation*}
                N = \Obound{\frac{R}{\eta c}\log \frac{1}{\varepsilon} + \frac{R^2}{\eta \varepsilon}}; \quad \quad \quad
                T = \Obound{\frac{d \tilde{\sigma}^2 M R^2}{\varepsilon c^2 \eta} \left( \frac{1}{c}\log \frac{1}{\varepsilon} + \frac{R}{\varepsilon} \right)}
            \end{equation*}
           
            number of iterations and zero-order oracle calls at
            \begin{equation*}
                \Delta \lesssim \frac{\varepsilon}{\sqrt{d} R (L_0 + L_1 M)}\min \left\{ \tilde{\sigma}, \frac{\varepsilon}{\sqrt{d} R}  \right\}
            \end{equation*}
            maximum noise level, where $c>0$ is clipping radius, $\expect{\norms{\nabla f(x,\xi)}^2} \leq \tilde{\sigma}^2$.
        \end{theorem}
    \end{boxM}
    It is not hard to see from Theorem~\ref{th:ZO_ClipSGD} that in the generalized smoothness condition, the iteration complexity at ZO-ClipSGD is the same as the first-order method of Algorithm~\ref{algo:clipSGD}. This effect is similar in the standard smoothness condition as well. The oracle complexity is $d$ times worse than its first-order counterpart due to the restriction to the oracle (Algorithm~\ref{algo:ZO_ClipSGD} uses only the zero-order oracle~\eqref{eq:zero_oracle}). It is worth noting that the maximum noise level $\Delta$ outperforms \cite{Kornilov_2023}, showing that generalized smoothness not only allows us to reach the summands characterizing the linear rate, but also improves the estimate on the maximum noise level (it is $\Delta$ that affects the error floor, in other words, the accuracy of the solution, to control the asymptote).  See the proof in the Appendix~\ref{app:ZO_ClipSGD}.

    \subsection{ZO-NSGD Method}
    Similar to the first-order algorithms, in this subsection we answer the question whether linear convergence can be achieved by the zero-order method in a convex setup. To answer this question, we consider the Zero-Order Normalized Stochastic Gradient Descent Method.

    \begin{algorithm}[H]
        \caption{Zero-Order Normalized Stochastic Gradient Descent Method (ZO-NSGD)}
        \label{algo:ZO_NSGD}
        \begin{algorithmic}
            \STATE {\bfseries Input:} initial point $x_0 \in \mathbb{R}^d$, iterations number $N$, batch size $B$, step size $\eta_k >0$ 
            \FOR{$k=0$ {\bfseries to} $N-1$}
                \STATE {~~~1.} Draw fresh i.i.d. samples $\xi^k_1,...,\xi^k_B$ and $e^k_1,...,e^k_B$
                \vspace{0.2em}
                \STATE {~~~2.} $\gg(x^k, \{\ee^k, \bxi^k\}) = \frac{1}{B} \sum_{i=1}^B \gg(x^k, \{e^k_i, \xi^k_i\})$ via~\eqref{eq:approximation_gradient}  
                \vspace{0.4em}
                \STATE {~~~3.} $x^{k+1} \gets x^{k}  - \eta_{k} \cdot \frac{\gg(x^k, \{\ee^k, \bxi^k\})}{\norms{\gg(x^k, \{\ee^k, \bxi^k\})}}$
            \ENDFOR
            \STATE {\bfseries Return:} $x^{N}$
        \end{algorithmic}
    \end{algorithm}  
    \vskip -0.2in
    In a similar way as in the previous subsection, we first generalize Theorem~\ref{th:NSGD} to the case with a biased gradient oracle to substitute estimates on the bias and the second moment of the gradient approximation to find optimality criteria for the gradient-free algorithm ZO-NSGD.
    \begin{lemma}
        \label{lem:biased_ZO_NSGD}
        Let function $f$ satisfy Assumption~\ref{ass:L0_L1_smooth} ($(L_0,L_1)$-smoothness) and \underline{biased gradient oracle}~\eqref{eq:gradient_oracle} $(\zeta > 0)$ satisfy Assumption~\ref{ass:bounded_var} (bounded variance), then Algorithm~\ref{algo:ZO_NSGD} with hyperparameter $\lambda > 0$ and step size ${\eta_k = \eta \leq \lambda / \left[2 (L_0 + L_1 \lambda)\right]}$ guarantees:
            \vspace{-0.8em}
            \begin{align*}
                \expect{f(x^N)} - f^* \lesssim  \left( 1 - \frac{\eta}{R} \right)^{N} (f(x^0) - f^*) + \frac{M R}{\lambda^2} \left( \frac{\sigma^2}{B} + \zeta^2 \right) + \lambda R.
            \end{align*}
    \end{lemma}

    From Lemma~\ref{lem:biased_ZO_NSGD} we can see how the inaccuracy accumulates over the iteration. The summand with $\zeta^2$ is unimprovable for first-order unaccelerated algorithms (see, e.g., \cite{Devolder_2013,Gasnikov_2024}). Now, having obtained the results for the biased NSGD we can use them to derive convergence results for Algorithm~\ref{algo:ZO_NSGD}.

    \begin{boxM}
        \begin{theorem}\label{th:ZO_NSGD}
            Let function $f$ satisfy Assumption~\ref{ass:L0_L1_smooth} ($(L_0,L_1)$-smoothness) and gradient approximation~\eqref{eq:approximation_gradient} satisfy Assumption~\ref{ass:bounded_var} (bounded variance), then ZO-NSGD with step size ${\eta_k = \eta \leq \lambda / \left[2 (L_0 + L_1 \lambda)\right]}$ converges to  desired $\varepsilon$ accuracy~$\left( \expect{f(x^N)} - f^* \leq \varepsilon \right)$~after:
            \begin{align*}
                N = \Obound{\frac{ R}{\eta} \log \frac{1}{\varepsilon}}; \quad \quad \quad  T = \OboundTilde{\frac{d \tilde{\sigma}^2 M R^4}{\varepsilon^3 \eta}}
            \end{align*}
            number of iterations and zero-order oracle calls~\eqref{eq:zero_oracle} at
            \begin{equation*}
                \Delta \lesssim \frac{\varepsilon^{3/2}}{\sqrt{d} R^{3/2} (L_0 + L_1 M)}\min \left\{ \tilde{\sigma}, \frac{\varepsilon^{3/2}}{\sqrt{d} R^{3/2}}  \right\}
            \end{equation*}
            maximum noise level, where $\lambda>0$ is hyperparameter, $\expect{\norms{\nabla f(x,\xi)}^2} \leq \tilde{\sigma}^2$.
        \end{theorem}
    \end{boxM}
    It is not hard to see that given a restricted oracle~\eqref{eq:zero_oracle}, Theorem~\ref{th:ZO_NSGD} shows that ZO-NSGD requires ${N = \Obound{\frac{R}{\eta} \log \frac{1}{\varepsilon}}}$ iterations to achieve the desired accuracy, which corresponds to a linear rate. However, it is worth noting that, as in Theorem~\ref{th:NSGD}, if we take the maximum step size $\eta = \lambda / \left[2 (L_0 + L_1 \lambda)\right]$, then the summand with $L_0$ still shows sublinear convergence $\OboundTilde{\frac{L_0 R^2}{\varepsilon}}$, despite the presence of the summand with $L_1$. But in the $L_0 = 0$ regime, \textit{we can say unambiguously that the zero-order algorithm ZO-NSGD can converge with a linear rate to the desired accuracy in the convex setup if we take the batch size $B = \Obound{\frac{d \tilde{\sigma}^2 M R^3}{\varepsilon^3}}$.} \textit{Theorem~\ref{th:ZO_NSGD} shows that the power of generalized smoothness (together with Batch parallelization) extends to zero-order algorithms.} Comparing the result of ZO-NSGD with Theorem~\ref{th:ZO_ClipSGD}, we can see that for a significant improvement in iteration complexity $N$ we “pay” for a deterioration in both oracle complexity $T$ and maximum noise level $\Delta$, which seems quite natural. See the proof of Lemma~\ref{lem:biased_ZO_NSGD} and Theorem~\ref{th:ZO_NSGD}~in~Appendix~\ref{app:ZO_NSGD}.

\section{Discussion and Future Works}\label{sec:Discussion and Future Works}
In Sections~\ref{sec:Clipped Stochastic Gradient Descent}-\ref{sec:Zero-Order Algorithms}, we gave two strategies for obtaining summands that characterize the linear rate in convergence estimates of algorithms for the convex setting: clipping (see Section~\ref{sec:Clipped Stochastic Gradient Descent}) and normalization (see Section~\ref{sec:Normalized Stochastic Gradient Descent}) techniques. Although these summands are quite unexpected for the convex setup and improve the estimates on the iteration complexity, we cannot claim linear convergence in general, since the convergence is dominated by the summands characterizing the sublinear rate. However, as we noted in Theorems~\ref{th:NSGD} and \ref{th:ZO_NSGD}, in the regime $L_0 = 0$, the NSGD and ZO-NSGD methods break all existing bounds on iteration complexity, demonstrating that it is possible to converge with linear rate to the desired accuracy with the condition of using Batch parallization. This result is pleasantly surprising and opens up a number of directions for future research.

In this paper we have focused on iteration complexity, so we see a careful analysis of the optimality criterion in the aggregate as future work. In particular, it seems interesting to show that $M$ is indeed bounded, e.g., using the technique from \cite{Li_2023_Convex} and evaluating with respect to the smoothness constants $L_0$ and $L_1$. The existence of the regime $L_0 = 0$ which allows one to achieve a linear convergence rate in the convex setting prompts the following question: can the iteration complexity be improved by assuming, for example, strong convexity, the PL condition, etc.? It is also interesting to see if similar effects are found in accelerated, adaptive algorithms, variational inequalities, distributed learning, nonsmooth (or increased smoothness) problems, overparameterization,~online~optimization,~etc.
\vspace{-0.5em}
\section{Numerical Experiments}\label{sec:Numerical Experiments}
\vspace{-0.5em}
In this section, we numerically analyze the algorithms presented in this paper and show that linear convergence in stochastic convex optimization is possible. For this illustration, we have chosen a problem that is of particular interest in the machine learning community: the logistic regression problem on w1a dataset \cite{Platt_1998}. We consider the following convex problem statement~\eqref{eq:init_problem}:
\begin{equation*}
    \min \limits_{x \in \mathbb{R}^d} f(x) = \frac{1}{M} \sum \limits_{i = 1}^M f_i(x), \quad \quad f_i(x) = \log \left(1 + \exp(-y_i \cdot (Ax)_i)\right),
\end{equation*}
where $f_i(x)$ is the loss on the $i$-th data point, $A \in \mathbb{R}^{M \times d}$ is an instances matrix, $y \in \{-1, 1\}^{M}$ is a label vector and $x \in \mathbb{R}^d$ is a vector of weights. It is easy to show, that logistic regression function is $L$-smooth (see Assumption~\ref{ass:L_smooth}) with $L = \frac{1}{4M} \sqrt{\lambda_{\max} (A^TA)}$, where $\lambda_{\max} (A^TA)$ denotes the largest eigenvalue of the matrix $A^TA$. Moreover, such a problem statement is a special case of \eqref{eq:init_problem} with $\xi$ being a random variable with the uniform distribution on $\{1, . . . , M\}$.

\begin{figure}[H]
\centering
\begin{minipage}{.48\textwidth}
  \centering
  \includegraphics[width=0.99\linewidth]{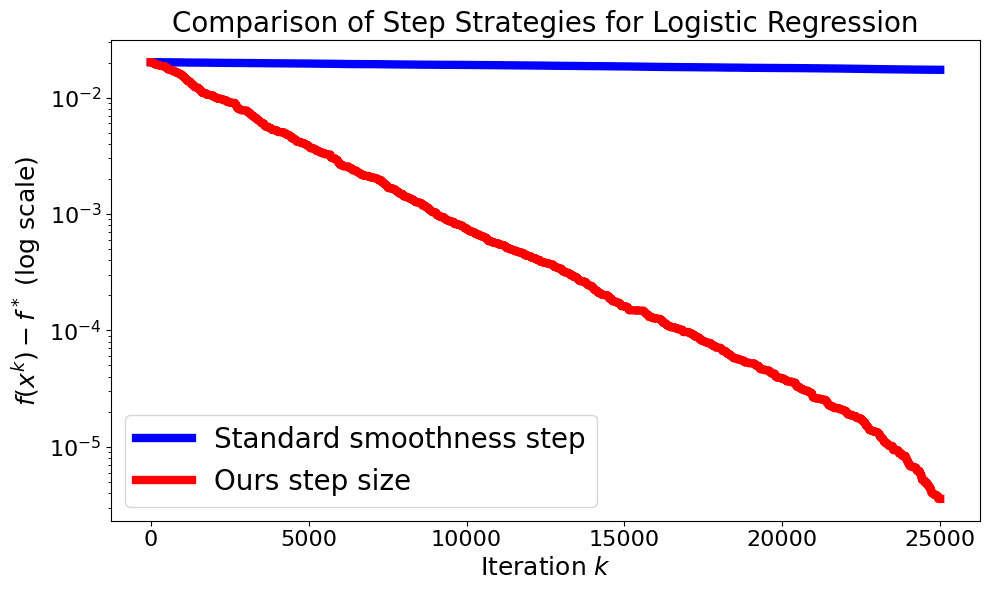}
    \caption{Linear convergence demonstration.}
    \label{fig:NSGD_compare1}
\end{minipage}%
\begin{minipage}{.04\textwidth}
\phantom{02}
\end{minipage}%
\begin{minipage}{.48\textwidth}
  \centering
  \includegraphics[width=0.99\linewidth]{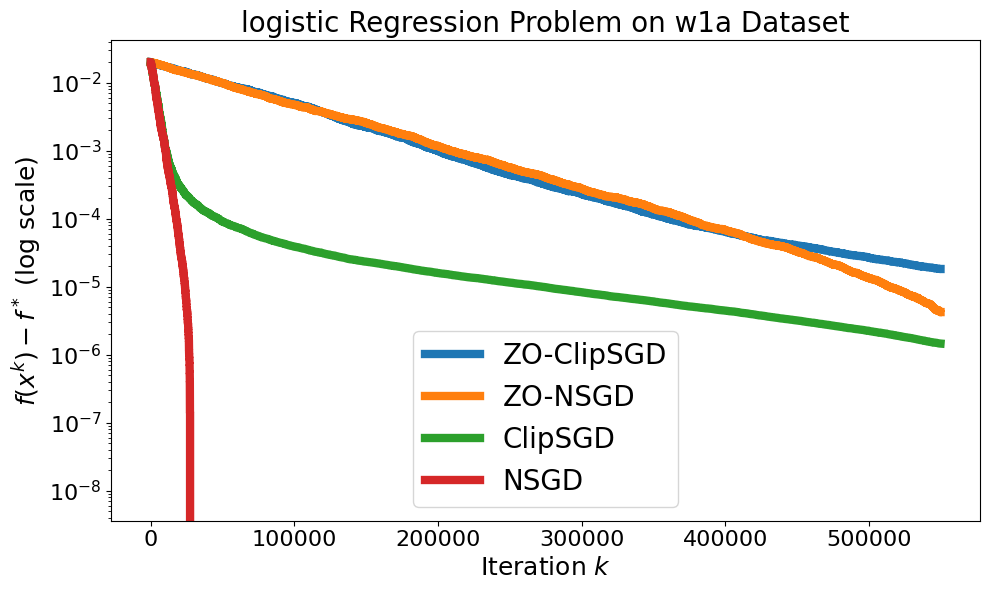}
    \caption{Comparison of the considered algorithms.}
    \label{fig:All_algorithms}
\end{minipage}
\vspace{-1em}
\end{figure}
In all tests we used the following parameters: $M = 2477$ - number of data, $d = 300$ - problem dimension, $B = 10$ - batch size, $h = 10^{-5}$ - smoothing parameter, $\Delta = 10^{-9}$ - noise level. Figure~\ref{fig:NSGD_compare1} shows a comparison of two step strategies of the NSGD algorithm. The \textcolor{blue}{blue line} (see "Standard smoothness step"), which corresponds to the theoretical step size in the $L$-smoothness setting demonstrates slow (sublinear) convergence. In particular, NSGD with this choice of step advanced the function from $0.02014151259$ to $0.01731953499$ in $25000$ iterations. The \textcolor{red}{red line} (see "Ours step size"), which corresponds to the next step size $\eta = \frac{1}{\norms{A}_{1}}$, demonstrates linear convergence that significantly outperforms the strategy with the theoretical step size. Moreover, \textit{Figure~\ref{fig:NSGD_compare1} demonstrates that indeed the first-order algorithm can converge with linear rate~to~the~desired~accuracy~in~a~convex~setup!} 

Figure~\ref{fig:All_algorithms} demonstrates the convergence dynamics of all the algorithms considered in this paper. In particular, as in Figure~\ref{fig:NSGD_compare1}, NSGD (see \textcolor{red}{red line}) shows a real linear convergence. ClipSGD (see \textcolor{green}{green line}), which used a step size $\eta = \frac{1}{c \cdot \norms{A}_{1}}$, where $c = 10^{-1}$ is the clipping radius, shows two modes of convergence: as long as $\norms{\nabla f(x^k)} \geq c$ the algorithm converges with a linear rate matching the NSGD, as soon as $\norms{\nabla f(x^k)} < c$ the algorithm slows down to the sublinear rate. The dynamics are similar for zero-order algorithms: ZO-NSGD (see \textcolor{orange}{orange line}), ZO-ClipSGD (see \textcolor{blue}{blue line}). Expectedly, these algorithms converge slower on the first iterations than their first-order counterparts due to restricted access to the oracle~\eqref{eq:zero_oracle}. However, it is worth noting that \textit{ZO-NSGD also exhibits linear convergence, thereby outperforming the first-order ClipSGD algorithm after $55000$ iterations.}
\vspace{-0.5em}

\section{Conclusion}\label{sec:Conclusion}
\vspace{-0.5em}
In this paper, we considered a stochastic convex optimization problem under the generalized smoothness condition of the objective function. We are the first who have provided strategies to achieve summands characterizing linear rate, thereby improving the iteration complexity (see Sections~\ref{sec:Clipped Stochastic Gradient Descent} and \ref{sec:Normalized Stochastic Gradient Descent}). In Section \ref{sec:Zero-Order Algorithms}, we showed that this effect of generalized smoothness extends to zero-order algorithms as well. Moreover, we highlight the regime $L_0 = 0$, under which we theoretically guarantee linear convergence for the NSGD and ZO-NSGD algorithms (subject to the use of batch parallization). This is the first result demonstrating linear convergence in such a problem setting, thus opening up a number of future works (see Section~\ref{sec:Discussion and Future Works}). Finally, in Section~\ref{sec:Numerical Experiments} we show using numerical experiments that linear convergence in such a problem formulation is also possible in practice.

\medskip

\bibliographystyle{plainnat-fixed}
{\small
\bibliography{3_references}
}

\newpage
\appendix

\vbox{%
    \hsize\textwidth
    \linewidth\hsize
    \vskip 0.1in
    \vskip 0.29in
    \vskip -\parskip
    \hrule height 1pt
    \vskip 0.19in
    \vskip 0.09in
    \begin{center}
        \LARGE \bf APPENDIX \\ Power of Generalized Smoothness in Stochastic Convex Optimization: First- and Zero-Order Algorithms
    \end{center}
    \vskip 0.29in
    \vskip -\parskip
    \hrule height 1pt
    \vskip 0.09in}


\section{Auxiliary Results}
In this section we provide auxiliary materials that are used in the proof of Theorems.
\subsection{Basic inequalities and assumptions}
\paragraph{Basic inequalities.} For all $a,b \in \mathbb{R}^d$ ($d \geq 1$) the following equality holds:
\begin{equation}
    \label{eq:qudrat_raznosti}
    2 \dotprod{a}{b} - \| b \|^2 = \| a \|^2 - \|a - b \|^2,
\end{equation}
\begin{equation}
    \label{eq:scalar_product_bound}
    \dotprod{a}{b} \leq \| a \| \cdot \| b \|.
\end{equation}

\paragraph{Squared norm of the sum} For all $a_1,...,a_n \in \mathbb{R}^d$, where $n=\{2,3\}$
    \begin{equation}
        \label{eq:Squared_norm_of_the_sum}
        \|a_1 + ... + a_n \|_2^2 \leq n \| a_1 \|_2^2 + ... + n \| a_n \|_2^2.
    \end{equation}

\paragraph{Generalized-Lipschitz-smoothness.}Throughout this paper, we assume that the $(L_0,L_1)$-smoothness condition (Assumption \ref{ass:L0_L1_smooth}) is satisfied. This inequality can be represented in the equivalent form for any $x,y \in \mathbb{R}^d$:
\begin{equation}
    \label{eq:ass_smooth}
    f(y) - f(x) \leq \dotprod{\nabla f(x)}{ y - x} + \frac{L_0 + L_1 \norms{\nabla f(x)}}{2} \norms{y-x}^2,
\end{equation}
where $L_0, L_1 \geq 0$ for any $x \in \mathbb{R}^d$ and $\norms{y-x} \leq \frac{1}{L_1}$.

\paragraph{Variance decomposition.} If $\xi$ is random vector in $\mathbb{R}^d$ with bounded second moment, then 
\begin{equation}\label{eq:Variance_decomposition}
    \expect{\norms{\xi + a}^2} = \expect{\norms{\xi - \expect{\xi}}^2} + \expect{\norms{\expect{\xi} - a}^2},
\end{equation}
for any deterministic vector $a\in \mathbb{R}^d$.

\subsection{Auxiliary Lemma about Generalized Smoothness}
If Assumption~\ref{ass:L0_L1_smooth} holds, then it also holds that $\forall x \in \mathbb{R}^d$:
\begin{equation}
    \norms{\nabla f(x)}^2 \leq 2(L_0 + L_1 \norms{\nabla f(x)}) (f(x) - f^*),\label{eq:lemma_smoothness}
\end{equation}
where $f^*=\inf_{x} f(x)$.
\begin{proof}
    We start the proof by applying \eqref{eq:ass_smooth} for $y = x - \frac{1}{L_0 + L_1 \norms{\nabla f(x)}} \nabla f(x)$, where ${\norms{ y-x} = \frac{\norms{\nabla f(x)}}{L_0 + L_1 \norms{\nabla f(x)}} \leq \frac{1}{L_1}}$. Then we can obtain:
    \begin{equation*}
        f^* \leq f \left( x - \frac{1}{L_0 + L_1 \norms{\nabla f(x)}} \nabla f(x)  \right) \overset{\eqref{eq:ass_smooth}}{\leq} f(x) - \frac{1}{2 (L_0 + L_1 \norms{\nabla f(x)})} \norms{\nabla f(x)}^2.
    \end{equation*}
\end{proof}
\subsection{Wirtinger-Poincare inequality}
        Let $f$ is differentiable, then for all $x \in \mathbb{R}^d$, $\gamma e \in S^d(\gamma)$:
        \begin{equation}\label{eq:Wirtinger_Poincare}
            \expect{f(x+ \gamma e)^2} \leq \frac{\gamma^2}{d} \expect{\norms{\nabla f(x + \gamma e)}^2}.
        \end{equation}

\newpage

\section{Clipped Stochastic Gradient Descent (Proof of the Theorem~\ref{th:clipSGD})}\label{app:ClipSGD}
We start by using $(L_0,L_1)$-smoothness (see Assumption~\ref{ass:L0_L1_smooth}):
\begin{align}
    f(x^{k+1}) - f(x^k) &\overset{\eqref{eq:ass_smooth}}{\leq} \dotprod{\nabla f(x^k)}{x^{k+1} - x^k} + \frac{L_0 + L_1 \norms{\nabla f(x^k)}}{2} \norms{x^{k+1} - x^{k}}^2 \nonumber \\
    &= - \eta \dotprod{\nabla f(x^k)}{\clip{\nabla f(x^k, \bxi^k)}} \nonumber\\
    &\quad \quad \quad+ \frac{\eta^2 (L_0 + L_1 \norms{\nabla f(x^k)})}{2} \norms{\clip{\nabla f(x^k, \bxi^k)}}^2. \label{eq:clipSGD_smooth}
\end{align}

Next, we consider three cases depending on the gradient norm: $\norms{\nabla f(x^k)} \geq c$ -- the full gradient is clipped and ${\frac{c}{2} \leq \norms{\nabla f(x^k)} \leq c}$ and $\norms{\nabla f(x^k)} \leq \frac{c}{2}$ -- the full gradient is not clipped.

\subsection{First case: \texorpdfstring{$\norms{\nabla f(x^k)} \geq c$}{TEXT}}

In this case $\alpha \nabla f(x^k) = \clip{\nabla f(x^k)}$ with $\alpha = \min \left\{1, \frac{c}{\norms{\nabla f(x^k)}}\right\} = \frac{c}{\norms{\nabla f(x^k)}}$, therefore we have the following
\begin{align}
    - \eta \dotprod{\nabla f(x^k)}{\clip{f(x^k, \bxi^k)}} &\overset{\eqref{eq:qudrat_raznosti}}{=} - \frac{\alpha \eta}{2} \norms{\nabla f(x^k)}^2 - \frac{\eta}{2 \alpha} \norms{\clip{\nabla f(x^k, \bxi^k)}}^2 \nonumber\\
    &\quad\quad\quad+ \frac{\eta}{2 \alpha} \norms{\clip{\nabla f(x^k, \bxi^k)} - \alpha \nabla f(x^k)}^2 \nonumber\\
    &= - \frac{\alpha \eta}{2} \norms{\nabla f(x^k)}^2 - \frac{\eta}{2 \alpha} \norms{\clip{\nabla f(x^k, \bxi^k)}}^2 \nonumber\\
    &\quad\quad\quad+ \frac{\eta}{2 \alpha} \norms{\clip{\nabla f(x^k, \bxi^k)} - \clip{\nabla f(x^k)}}^2 \nonumber\\
    &= - \frac{c \eta}{2} \norms{\nabla f(x^k)} - \frac{\eta}{2 \alpha} \norms{\clip{\nabla f(x^k, \bxi^k)}}^2 \nonumber\\
    &\quad\quad\quad+ \frac{\eta}{2 \alpha} \norms{\clip{\nabla f(x^k, \bxi^k)} - \clip{\nabla f(x^k)}}^2. \nonumber
\end{align}
Using that clipping is a projection on onto a convex set, namely ball with radius $c$, and thus is Lipshitz operator with Lipshitz constant $1$, we can obtain:
\begin{align}
    - \eta \dotprod{\nabla f(x^k)}{\expect{\clip{f(x^k, \bxi^k)}}} &\leq - \frac{c \eta}{2} \norms{\nabla f(x^k)} - \frac{\eta}{2 \alpha} \expect{\norms{\clip{\nabla f(x^k, \bxi^k)}}^2} \nonumber\\
    &\quad\quad\quad+ \frac{\eta}{2 \alpha} \expect{\norms{\nabla f(x^k, \bxi^k) - \nabla f(x^k)}^2} \nonumber\\
    &\leq  - \frac{c \eta}{2} \norms{\nabla f(x^k)} - \frac{\eta}{2 \alpha} \expect{\norms{\clip{\nabla f(x^k, \bxi^k)}}^2} \nonumber\\
    &\quad\quad\quad+ \frac{\eta \sigma^2}{2 \alpha B} \nonumber\\
    &= - \frac{c \eta}{2} \norms{\nabla f(x^k)} - \frac{\eta}{2 \alpha} \expect{\norms{\clip{\nabla f(x^k, \bxi^k)}}^2} \nonumber\\
    &\quad\quad\quad+ \frac{\eta \norms{\nabla f(x^k)} \sigma^2}{2 c B}. \label{eq:clipSGD_case1_noise}
\end{align}
We now consider the cases depending on the relation between $c$ and $\sigma$:

\paragraph{\fbox{In the case $c \geq \sqrt{2} \sigma$}}
    
    We have in \eqref{eq:clipSGD_case1_noise}:
    \begin{align*}
        - \eta \dotprod{\nabla f(x^k)}{\expect{\clip{\nabla f(x^k, \bxi^k)}}} &\overset{\eqref{eq:clipSGD_case1_noise}}{\leq} - \frac{c \eta}{2} \norms{\nabla f(x^k)} - \frac{\eta}{2 \alpha} \expect{\norms{\clip{\nabla f(x^k, \bxi^k)}}^2} \nonumber\\
    &\quad\quad\quad+ \frac{\eta \norms{\nabla f(x^k)} \sigma^2}{2 c B} \\
        &= - \frac{\eta}{2 \alpha} \expect{\norms{\clip{\nabla f(x^k, \bxi^k)}}^2} \nonumber\\
    &\quad\quad\quad- \frac{c \eta}{2} \norms{\nabla f(x^k)} \left( 1 - \frac{\sigma^2}{c^2 B} \right)\\
        &\leq - \frac{\eta}{2 \alpha} \expect{\norms{\clip{\nabla f(x^k, \bxi^k)}}^2} \nonumber\\
    &\quad\quad\quad- \frac{c \eta}{4} \norms{\nabla f(x^k)} \\ 
        &= - \frac{\eta \norms{\nabla f(x^k)}}{2 c} \expect{\norms{\clip{\nabla f(x^k, \bxi^k)}}^2} \nonumber\\
    &\quad\quad\quad- \frac{c \eta}{4} \norms{\nabla f(x^k)}.
    \end{align*}
    Plugging this into \eqref{eq:clipSGD_smooth} and choosing $\eta \leq \frac{1}{4(L_0 + L_1 c)}$ we have:
    \begin{align}
        \expect{\nabla f(x^{k+1})} - f(x^k) &\overset{\eqref{eq:clipSGD_smooth}}{\leq} - \frac{\eta \norms{\nabla f(x^k)}}{2 c} \expect{\norms{\clip{\nabla f(x^k, \bxi^k)}}^2} - \frac{c \eta}{4} \norms{\nabla f(x^k)} \nonumber\\
        & \quad \quad\quad + \frac{\eta^2 (L_0 + L_1 \norms{\nabla f(x^k)})}{2} \expect{\norms{\clip{\nabla f(x^k, \bxi^k)}}^2}\nonumber\\
        &= - \frac{\eta \norms{\nabla f(x^k)}}{2 c} \expect{\norms{\clip{\nabla f(x^k, \bxi^k)}}^2} (1 - \eta L_1 c) \nonumber\\
    &\quad\quad\quad- \frac{c \eta}{4} \norms{\nabla f(x^k)}  + \frac{\eta^2 L_0}{2}\expect{\norms{\clip{\nabla f(x^k, \bxi^k)}}^2}\nonumber\\
        &\leq  - \frac{c \eta}{4} \norms{\nabla f(x^k)} - \frac{\eta}{2} \expect{\norms{\clip{\nabla f(x^k, \bxi^k)}}^2} \left(1 - \eta( L_0 + L_1 c)\right)\nonumber\\
        &\leq  - \frac{c \eta}{4} \norms{\nabla f(x^k)}. \label{eq:clipSGD_case1_withoutnoise_1}
    \end{align}
    Using the convexity assumption of the function, we have the following:
    \begin{align*}
        f(x^k) - f^* &\leq \dotprod{\nabla f(x^k)}{x^k - x^*}\\
        &\overset{\eqref{eq:scalar_product_bound}}{\leq} \norms{\nabla f(x^k)} \norms{x^k - x^*}\\
        & \leq \norms{\nabla f(x^k)} \underbrace{\norms{x^0 - x^*}}_{R}.
    \end{align*}
    Hence we have:
    \begin{equation}
        \norms{\nabla f(x^k)} \geq \frac{f(x^k) - f^*}{R}. \label{eq:clipSGD_case1_withoutnoise_2}
    \end{equation}

    Then substituting \eqref{eq:clipSGD_case1_withoutnoise_2} into \eqref{eq:clipSGD_case1_withoutnoise_1} we obtain:
    \begin{equation*}
        \expect{f(x^{k+1})} - f(x^k) \leq - \frac{\eta c}{4} \norms{\nabla f(x^k)} \leq - \frac{\eta c}{4 R} (f(x^k) - f^*).
    \end{equation*}
    
    This inequality is equivalent to the trailing inequality:
    \begin{equation*}
        \expect{f(x^{k+1})} - f^* \leq \left( 1 - \frac{\eta c}{4 R} \right) \left( f(x^k) - f^* \right).
    \end{equation*}
     Then for $k = 0,1,2,..., N-1$ iterations that satisfy the conditions $\norms{\nabla f(x^k)} \geq c \geq \sqrt{2} \sigma$, then ClipSGD has linear convergence
    \begin{equation*}
        \expect{f(x^{N})} - f^* \leq \left( 1 - \frac{\eta}{2 R} \right)^{N} \left( f(x^0) - f^* \right).
    \end{equation*}

\paragraph{\fbox{In the case $c \leq \sqrt{2} \sigma$}}

We have in \eqref{eq:clipSGD_case1_noise}:
    \begin{align*}
        - \eta \dotprod{\nabla f(x^k)}{\expect{\clip{f(x^k, \bxi^k)}}} &\overset{\eqref{eq:clipSGD_case1_noise}}{\leq} - \frac{c \eta}{2} \norms{\nabla f(x^k)} - \frac{\eta}{2 \alpha} \expect{\norms{\clip{\nabla f(x^k, \bxi^k)}}^2} \nonumber\\
    &\quad\quad\quad+ \frac{\eta \norms{\nabla f(x^k)} \sigma^2}{2 c B} \\
        &= - \frac{c \eta}{2} \norms{\nabla f(x^k)} - \frac{\eta}{2 \alpha} \expect{\norms{\clip{\nabla f(x^k, \bxi^k)}}^2} \nonumber\\
    &\quad\quad\quad+ \frac{\eta M \sigma^2}{2 c B}.
    \end{align*}
    Plugging this into \eqref{eq:clipSGD_smooth} and choosing $\eta \leq \frac{1}{4(L_0 + L_1 c)}$ we have:
    \begin{align}
        \expect{f(x^{k+1})} - f(x^k) &\overset{\eqref{eq:clipSGD_smooth}}{\leq} - \frac{c \eta}{2} \norms{\nabla f(x^k)} - \frac{\eta \norms{\nabla f(x^k)}}{2 c} \expect{\norms{\clip{\nabla f(x^k, \bxi^k)}}^2}  \nonumber\\
        & \quad \quad + \frac{\eta^2 (L_0 + L_1 \norms{\nabla f(x^k)})}{2} \expect{\norms{\clip{\nabla f(x^k, \bxi^k)}}^2} + \frac{\eta M \sigma^2}{2 c B} \nonumber\\
        &= - \frac{c \eta}{2} \norms{\nabla f(x^k)} - \frac{\eta \norms{\nabla f(x^k)}}{2 c} \expect{\norms{\clip{\nabla f(x^k, \bxi^k)}}^2} (1 - \eta L_1 c) \nonumber\\
        & \quad \quad + \frac{\eta^2 L_0}{2}\expect{\norms{\clip{\nabla f(x^k, \bxi^k)}}^2} + \frac{\eta M \sigma^2}{2 c B}\nonumber\\
        &\leq  - \frac{c \eta}{2} \norms{\nabla f(x^k)} - \frac{\eta}{2} \expect{\norms{\clip{\nabla f(x^k, \bxi^k)}}^2} \left(1 - \eta( L_0 + L_1 c)\right) \nonumber\\
    &\quad\quad\quad+ \frac{\eta M \sigma^2}{2 c B}\nonumber\\
        &\leq  - \frac{c \eta}{2} \norms{\nabla f(x^k)} + \frac{\eta M \sigma^2}{2 c B}. \label{eq:clipSGD_case1_with_noise_1}
    \end{align}
    Using the convexity assumption of the function, we have the following:
    \begin{align*}
        f(x^k) - f^* &\leq \dotprod{\nabla f(x^k)}{x^k - x^*}\\
        &\overset{\eqref{eq:scalar_product_bound}}{\leq} \norms{\nabla f(x^k)} \norms{x^k - x^*}\\
        & \leq \norms{\nabla f(x^k)} \underbrace{\norms{x^0 - x^*}}_{R}.
    \end{align*}
    Hence we have:
    \begin{equation}
        \norms{\nabla f(x^k)} \geq \frac{f(x^k) - f^*}{R}. \label{eq:clipSGD_case1_with_noise_2}
    \end{equation}

    Then substituting \eqref{eq:clipSGD_case1_with_noise_2} into \eqref{eq:clipSGD_case1_with_noise_1} we obtain:
    \begin{equation*}
        \expect{f(x^{k+1})} - f(x^k) \leq - \frac{\eta c}{2} \norms{\nabla f(x^k)} + \frac{\eta M \sigma^2}{2 c B} \leq - \frac{\eta c}{2 R} (f(x^k) - f^*) + \frac{\eta M \sigma^2}{2 c B}.
    \end{equation*}
    
    This inequality is equivalent to the trailing inequality:
    \begin{equation*}
        \expect{f(x^{k+1})} - f^* \leq \left( 1 - \frac{\eta c}{2 R} \right) \left( f(x^k) - f^* \right) + \frac{\eta M \sigma^2}{2 c B}.
    \end{equation*}
     Then for $k = 0,1,2,..., N-1$ iterations that satisfy the conditions $\norms{\nabla f(x^k)} \geq c$ and $c \leq \sqrt{2} \sigma$, then ClipSGD has linear convergence
    \begin{equation*}
        \expect{f(x^{N})} - f^* \leq \left( 1 - \frac{\eta c}{2 R} \right)^{N} \left( f(x^0) - f^* \right) + \frac{MR \sigma^2}{c^2 B}.
    \end{equation*}

\subsection{Second case: \texorpdfstring{$\frac{c}{2} \leq \norms{\nabla f(x^k)} \leq c$}{TEXT}}

In this case $\nabla f(x^k) = \clip{\nabla f(x^k)}$ with $\alpha = \min \left\{1, \frac{c}{\norms{\nabla f(x^k)}}\right\} = 1$, therefore we have the following
\begin{align}
    - \eta \dotprod{\nabla f(x^k)}{\clip{\nabla f(x^k, \bxi^k)}} &\overset{\eqref{eq:qudrat_raznosti}}{=} - \frac{\alpha \eta}{2} \norms{\nabla f(x^k)}^2 - \frac{\eta}{2 \alpha} \norms{\clip{\nabla f(x^k, \bxi^k)}}^2 \nonumber\\
    &\quad\quad\quad+ \frac{\eta}{2 \alpha} \norms{\clip{\nabla f(x^k, \bxi^k)} - \alpha \nabla f(x^k)}^2 \nonumber\\
    &= - \frac{\eta}{2} \norms{\nabla f(x^k)}^2 - \frac{\eta}{2} \norms{\clip{\nabla f(x^k, \bxi^k)}}^2 \nonumber\\
    &\quad\quad\quad+ \frac{\eta}{2} \norms{\clip{\nabla f(x^k, \bxi^k)} - \clip{\nabla f(x^k)}}^2 \nonumber\\
    &\leq - \frac{c \eta}{4} \norms{\nabla f(x^k)} - \frac{\eta}{2} \norms{\clip{\nabla f(x^k, \bxi^k)}}^2 \nonumber\\
    &\quad\quad\quad+ \frac{\eta}{2} \norms{\clip{\nabla f(x^k, \bxi^k)} - \clip{\nabla f(x^k)}}^2. \nonumber
\end{align}
Using that clipping is a projection on onto a convex set, namely ball with radius $c$, and thus is Lipshitz operator with Lipshitz constant $1$, we can obtain:
\begin{align*}
    - \eta \dotprod{\nabla f(x^k)}{\expect{\clip{\nabla f(x^k, \bxi^k)}}} &\leq - \frac{c \eta}{4} \norms{\nabla f(x^k)} - \frac{\eta}{2} \expect{\norms{\clip{\nabla f(x^k, \bxi^k)}}^2} \nonumber\\
    &\quad\quad\quad+ \frac{\eta}{2} \expect{\norms{\nabla f(x^k, \bxi^k) - \nabla f(x^k)}^2} \nonumber\\
    &\leq  - \frac{c \eta}{4} \norms{\nabla f(x^k)} - \frac{\eta}{2 } \expect{\norms{\clip{\nabla f(x^k, \bxi^k)}}^2} + \frac{\eta \sigma^2}{2 B} \nonumber\\
    &= - \frac{c \eta}{4} \norms{\nabla f(x^k)} - \frac{\eta}{2 } \expect{\norms{\clip{\nabla f(x^k, \bxi^k)}}^2} + \frac{\eta \sigma^2}{2B}. 
\end{align*}

Plugging this into \eqref{eq:clipSGD_smooth} and choosing $\eta \leq \frac{1}{4(L_0 + L_1 c)}$ we have:
    \begin{align}
        \expect{f(x^{k+1})} - f(x^k) &\overset{\eqref{eq:clipSGD_smooth}}{\leq} - \frac{c \eta}{4} \norms{\nabla f(x^k)} - \frac{\eta}{2} \expect{\norms{\clip{\nabla f(x^k, \bxi^k)}}^2}  \nonumber\\
        & \quad \quad + \frac{\eta^2 (L_0 + L_1 \norms{\nabla f(x^k)})}{2} \expect{\norms{\clip{\nabla f(x^k, \bxi^k)}}^2} + \frac{\eta \sigma^2}{2 B} \nonumber\\
        &= - \frac{c \eta}{4} \norms{\nabla f(x^k)} + \frac{\eta \sigma^2}{2 B} \nonumber\\
    &\quad\quad\quad- \frac{\eta}{2} \expect{\norms{\clip{\nabla f(x^k, \bxi^k)}}^2} \left(1 -  \eta (L_0 + L_1 \norms{\nabla f(x^k)}) \right)  \nonumber\\
        &\leq - \frac{c \eta}{4} \norms{\nabla f(x^k)} + \frac{\eta \sigma^2}{2 B}.\label{eq:clipSGD2_case1_1}
    \end{align}

    Using the convexity assumption of the function, we have the following:
    \begin{align*}
        f(x^k) - f^* &\leq \dotprod{\nabla f(x^k)}{x^k - x^*}\\
        &\overset{\eqref{eq:scalar_product_bound}}{\leq} \norms{\nabla f(x^k)} \norms{x^k - x^*}\\
        & \leq \norms{\nabla f(x^k)} \underbrace{\norms{x^0 - x^*}}_{R}.
    \end{align*}
    Hence we have:
    \begin{equation}
        \norms{\nabla f(x^k)} \geq \frac{f(x^k) - f^*}{R}. \label{eq:clipSGD2_case1_2}
    \end{equation}

    Then substituting \eqref{eq:clipSGD2_case1_2} into \eqref{eq:clipSGD2_case1_1} we obtain:
    \begin{equation*}
        \expect{f(x^{k+1})} - f(x^k) \leq - \frac{\eta c}{4} \norms{\nabla f(x^k)} + \frac{\eta \sigma^2}{2 B} \leq - \frac{\eta c}{4 R} (f(x^k) - f^*) + \frac{\eta \sigma^2}{2 B}.
    \end{equation*}
    
    This inequality is equivalent to the trailing inequality:
    \begin{equation*}
        \expect{f(x^{k+1})} - f^* \leq \left( 1 - \frac{\eta c}{4 R} \right) \left( f(x^k) - f^* \right) + \frac{\eta \sigma^2}{2 B}.
    \end{equation*}
     Then for $k = 0,1,2,..., N-1$ iterations that satisfy the conditions $\frac{c}{2} \leq \norms{\nabla f(x^k)} \leq c $, then ClipSGD has linear convergence
    \begin{equation*}
        \expect{f(x^{N})} - f^* \leq \left( 1 - \frac{\eta c}{4 R} \right)^{N} \left( f(x^0) - f^* \right) + \frac{2 \sigma^2 R}{c B}.
    \end{equation*}

    Let $\mathcal{T}_1 = \left\{ m^{\mathcal{T}_1}_0, m^{\mathcal{T}_1}_1, m^{\mathcal{T}_1}_2, ..., m^{\mathcal{T}_1}_{K-1}  \right\} = \left\{k \in \{0,1,2,...,N-1\} |  \norms{\nabla f(x^k,\xi^k)} \geq \frac{c}{2} \right\}$, where $K = |\mathcal{T}_1|$. Then for $k \in \mathcal{T}_1$ ClipSGD shows linear convergence:
        \begin{align*}
            F_N \cdot \mathds{1}\left[ \mathcal{T}_1 \right] &\leq \left(1 - \frac{1}{4L_1R} \right) F_N \leq \left(1 - \frac{1}{4L_1R} \right) F_{m^{\mathcal{T}_1}_{K-1}} + \frac{\eta M \sigma^2}{2 c B} + \frac{\eta \sigma^2}{2 B} \nonumber\\
    &\leq  ... \leq \left(1 - \frac{1}{4L_1R} \right)^{K} F_{m_0^{\mathcal{T}_1}} + \frac{M R \sigma^2}{ c^2 B} + \frac{2 \sigma^2 R}{c B} \\
            &\leq \left(1 - \frac{1}{4L_1R} \right)^{K} F_0 + \frac{M R \sigma^2}{ c^2 B} + \frac{2 \sigma^2 R}{c B},
        \end{align*}
        where $F_k =  \expect{f(x^k)} - f^*$, and we used that $F_{k} \leq F_{k-1}$.

\subsection{Third case: \texorpdfstring{$\norms{\nabla f(x^k)} \leq \frac{c}{2}$}{TEXT}}

    We introduce an indicative function:
    \begin{equation}
        \aleph_k = \mathds{1}\left\{ \norms{\nabla f(x^k, \bxi^k)} > c \right\}.\label{eq:indicator_function}
    \end{equation}
    Then the following is true:
    \begin{equation}
        \expect{\aleph_k} = \expect{\aleph_k^2} =\prob{\norms{\nabla f(x^k, \bxi^k)} > c} \overset{\circledOne}{\leq} \prob{\norms{\nabla f(x^k, \bxi^k) - \nabla f(x^k) } > \frac{c}{2}} \overset{\circledTwo}{\leq} \frac{4 \sigma^2}{c^2 B}, \label{eq:prob_indicator}
    \end{equation}
    where in $\circledOne$ we used $\norms{\nabla f(x^k, \bxi^k)} \leq \norms{\nabla f(x^k, \bxi^k) - \nabla f(x^k) } + \norms{\nabla f(x^k) } \leq \norms{\nabla f(x^k, \bxi^k) - \nabla f(x^k) } + \frac{c}{2}$, and in $\circledTwo$ we used Markov's inequality.
    
    Let $r_{k+1} = \expect{\norms{x^{k+1} - x^*}}$ and $F_{k+1} = \expect{f(x^{k+1}) - f^*}$, then given that
    \begin{align*}
        \clip{\nabla f(x^k, \bxi^k)} &= \nabla f(x^k, \bxi^k) (1 - \aleph_k) + \frac{c}{\norms{\nabla f(x^k, \bxi^k)}} \nabla f(x^k, \bxi^k) \aleph_k\\
        & = \nabla f(x^k, \bxi^k) + \left( \frac{c}{\norms{\nabla f(x^k, \bxi^k)}} - 1 \right) \nabla f(x^k, \bxi^k) \aleph_k
    \end{align*}
    we get with $\eta \leq \frac{1}{4(L_0 + L_1 c)}$:
    \begin{align}
        r_{k+1}^2 
        &= r_k^2 - 2 \eta \dotprod{\expect{\clip{\nabla f(x^k, \bxi^k)}}}{x^k - x^*} + \eta^2 \expect{\norms{\clip{\nabla f(x^k, \bxi^k)}}^2}\nonumber\\
        &= r_k^2 - 2 \eta \dotprod{\nabla f(x^k)}{x^k - x^*} - 2 \eta \dotprod{\expect{\left( \frac{c}{\norms{\nabla f(x^k, \bxi^k)}} - 1 \right) \nabla f(x^k, \bxi^k) \aleph_k}}{x^k - x^*} \nonumber\\
        &\quad\quad\quad + \eta^2 \expect{\norms{\clip{\nabla f(x^k, \bxi^k)}}^2}\nonumber\\
        &\overset{\eqref{eq:scalar_product_bound}}{\leq} r_k^2 - 2 \eta \dotprod{\nabla f(x^k)}{x^k - x^*} + 2 \eta \norms{\expect{\left( \frac{c}{\norms{\nabla f(x^k, \bxi^k)}} - 1 \right) \nabla f(x^k, \bxi^k) \aleph_k}}\norms{x^k - x^*} \nonumber\\
        &\quad\quad\quad + \eta^2 \expect{\norms{\clip{\nabla f(x^k, \bxi^k)}}^2}\nonumber\\
        &\overset{\circledOne}{\leq} r_k^2 - 2 \eta F_k + 2 \eta \norms{\expect{\left( \frac{c}{\norms{\nabla f(x^k, \bxi^k)}} - 1 \right) \nabla f(x^k, \bxi^k) \aleph_k}}\norms{x^0 - x^*} \nonumber\\
        &\quad\quad\quad + \eta^2 \expect{\norms{\clip{\nabla f(x^k, \bxi^k)}}^2}\nonumber\\
        &\overset{\eqref{eq:Squared_norm_of_the_sum}}{\leq}r_k^2 - 2 \eta F_k + 2 \eta \norms{\expect{\left( \frac{c}{\norms{\nabla f(x^k, \bxi^k)}} - 1 \right) \nabla f(x^k, \bxi^k) \aleph_k}}\norms{x^0 - x^*} \nonumber\\
        &\quad\quad\quad + 2 \eta^2 \expect{\norms{\clip{\nabla f(x^k, \bxi^k)} - \nabla f(x^k)}^2} + 2 \eta^2 \norms{\nabla f(x^k)}^2\nonumber\\
        &=r_k^2 - 2 \eta F_k + 2 \eta \norms{\expect{\left( \frac{c}{\norms{\nabla f(x^k, \bxi^k)}} - 1 \right) \nabla f(x^k, \bxi^k) \aleph_k}} R \nonumber\\
        &\quad\quad\quad + 2 \eta^2 \expect{\norms{\clip{\nabla f(x^k, \bxi^k)} - \clip{\nabla f(x^k)}}^2} + 2 \eta^2 \norms{\nabla f(x^k)}^2\nonumber\\
        &\overset{\circledTwo}{\leq}
        r_k^2 - 2 \eta F_k + 2 \eta \norms{\expect{\left( \frac{c}{\norms{\nabla f(x^k, \bxi^k)}} - 1 \right) \nabla f(x^k, \bxi^k) \aleph_k}}R \nonumber\\
        &\quad\quad\quad + 2 \eta^2 \expect{\norms{\nabla f(x^k, \bxi^k) - \nabla f(x^k)}^2} + 2 \eta^2 \norms{\nabla f(x^k)}^2\nonumber\\
        &\leq r_k^2 - 2 \eta F_k + 2 \eta \norms{\expect{\left( \frac{c}{\norms{\nabla f(x^k, \bxi^k)}} - 1 \right) \nabla f(x^k, \bxi^k) \aleph_k}}R \nonumber\\
        &\quad\quad\quad +  \frac{ 2 \eta^2 \sigma^2}{B} + 2 \eta^2 \norms{\nabla f(x^k)}^2\nonumber\\
        &\overset{\eqref{eq:lemma_smoothness}}{\leq} r_k^2 - 2 \eta F_k + 2 \eta \norms{\expect{\left( \frac{c}{\norms{\nabla f(x^k, \bxi^k)}} - 1 \right) \nabla f(x^k, \bxi^k) \aleph_k}}R \nonumber\\
        &\quad\quad\quad +  \frac{ 2 \eta^2 \sigma^2}{B} + 4 \eta^2 \left( L_0 + L_1 \norms{\nabla f(x^k)} \right) F_k\nonumber\\
        &\leq 
        r_k^2 - 2 \eta F_k + 2 \eta \norms{\expect{\left( \frac{c}{\norms{\nabla f(x^k, \bxi^k)}} - 1 \right) \nabla f(x^k, \bxi^k) \aleph_k}}R \nonumber\\
        &\quad\quad\quad +  \frac{ 2 \eta^2 \sigma^2}{B} + 4 \eta^2 \left( L_0 + L_1 c \right) F_k\nonumber\\
        &= r_k^2 - 2 \eta F_k \left(1 - 2 \eta \left( L_0 + L_1 c \right)  \right) +  \frac{ 2 \eta^2 \sigma^2}{B}  \nonumber\\
        &\quad\quad\quad + 2 \eta \norms{\expect{\left( \frac{c}{\norms{\nabla f(x^k, \bxi^k)}} - 1 \right) \nabla f(x^k, \bxi^k) \aleph_k}}R\nonumber\\
        &\leq r_k^2 -  \eta F_k  +  \frac{ 2 \eta^2 \sigma^2}{B} + 2 \eta \norms{\expect{\left( \frac{c}{\norms{\nabla f(x^k, \bxi^k)}} - 1 \right) \nabla f(x^k, \bxi^k) \aleph_k}}R.\label{eq:clipSGD3_case1_1}
    \end{align}

    Let's find the upper bound of the last summand:
    \begin{align}
        2 \eta R & \norms{\expect{\left( \frac{c}{\norms{\nabla f(x^k, \bxi^k)}} - 1 \right) \nabla f(x^k, \bxi^k) \aleph_k}}\nonumber\\
    &
         \overset{\eqref{eq:indicator_function}}{\leq}  2 \eta R \expect{\norms{\nabla f(x^k, \bxi^k)} \cdot \left( 1 - \frac{c}{ \norms{\nabla f(x^k, \bxi^k)}} \right) \aleph_k} \nonumber\\
        &\leq 2 \eta R  \expect{\norms{\nabla f(x^k, \bxi^k)} \cdot  \aleph_k}\nonumber\\
        & \leq 2 \eta R \left(\expect{\norms{\nabla f(x^k, \bxi^k) - \nabla f(x^k)} \cdot  \aleph_k} + \norms{\nabla f(x^k)} \expect{\aleph_k} \right)\nonumber\\
        &\leq 2 \eta R \left(\sqrt{\expect{\norms{\nabla f(x^k, \bxi^k) - \nabla f(x^k)}^2} \cdot  \expect{\aleph_k^2}} + \norms{\nabla f(x^k)} \expect{\aleph_k}\right)\nonumber\\
        &\overset{\eqref{eq:prob_indicator}}{\leq} 2 \eta R \left(\frac{2\sigma^2}{c B} + \frac{c}{2} \cdot\frac{4\sigma^2}{c^2 B}\right) \nonumber\\
        &= \frac{8 \eta \sigma^2 R}{c B}.\label{eq:clipSGD3_case1_2}
    \end{align}

    Substituting into the initial formula and rearrange the summands, we obtain
    \begin{align*}
        \eta F_k &\overset{\eqref{eq:clipSGD3_case1_1}}{\leq} r_{k}^2 - r_{k+1}^2 + \frac{2 \eta^2 \sigma^2}{B} + 2 \eta \norms{\expect{\left( \frac{c}{\norms{\nabla f(x^k, \bxi^k)}} - 1 \right) \nabla f(x^k, \bxi^k) \aleph_k}}R\\
        &\overset{\eqref{eq:clipSGD3_case1_2}}{\leq} r_{k}^2 - r_{k+1}^2 + \frac{2 \eta^2 \sigma^2}{B} + \frac{8 \eta \sigma^2 R}{c B}
    \end{align*}

    Let $\mathcal{T}_2 = \left\{m^{\mathcal{T}_2}_0, m^{\mathcal{T}_2}_1, m^{\mathcal{T}_2}_2, ..., m^{\mathcal{T}_2}_{N-K}  \right\} = \left\{k \in \{0,1,2,...,N-1\} |  \norms{\nabla f(x^k)} < \frac{c}{2} \right\}$, where $|\mathcal{T}_2| = N-K$. Then rearranging and summing over all $k \in \mathcal{T}_2$ we obtain
        \begin{align*}
            F_N \cdot \mathds{1}\left[ \mathcal{T}_2 \right] &\leq \frac{1}{N - K} \sum_{k\in \mathcal{T}_2} F_{k} \leq \frac{1}{\eta (N - K)} \sum_{k\in \mathcal{T}_2} \left( r_{k}^2 - r_{k+1}^2 \right) + \frac{1}{N - K} \sum_{k\in \mathcal{T}_2} \left(\frac{2 \eta \sigma^2}{B} + \frac{8 \sigma^2 R}{c B}\right)\\
            &=\frac{r_0^2 - r_{N}^2}{\eta (N - K)} + \frac{2 \eta \sigma^2}{B} + \frac{8 \sigma^2 R}{c B} \leq \frac{r_0^2}{\eta(N - K)}  + \frac{2 \eta \sigma^2}{B} + \frac{8 \sigma^2 R}{c B}. 
        \end{align*}
    Hence we obtain:
        \begin{equation*}
            F_N \cdot \mathds{1}\left[ \mathcal{T}_2 \right] \leq \frac{R^2}{\eta(N - K)}  + \frac{2 \eta \sigma^2}{B} + \frac{8 \sigma^2 R}{c B}.
        \end{equation*}
    
    Combining all cases we have:
    \begin{align*}
        \expect{f(x^{N})} - f^* &\leq F_N \cdot \mathds{1}\left[ \mathcal{T}_1 \right] + F_N \cdot \mathds{1}\left[ \mathcal{T}_2 \right] \\
        & \leq \left( 1 - \frac{\eta c}{2 R} \right)^{K} F_0 +  \frac{R^2}{\eta (N- K)} + \frac{\sigma^2 M R}{c^2 B} + \frac{2 \eta \sigma^2}{B} +  \frac{8 \sigma^2 R}{c B}.
    \end{align*}

\section{Normalized Stochastic Gradient Descent (Proof of the Theorem~\ref{th:NSGD})}\label{app:NSGD}

Let's introduce the notation $G(x^k, \bxi^k) = \frac{\nabla f(x^k, \bxi^k)}{\norms{\nabla f(x^k, \bxi^k)}}$, then using $(L_0,L_1)$-smoothness (see Assumption~\ref{ass:L0_L1_smooth}):
\begin{align}
    f(x^{k+1}) - f(x^k) &\overset{\eqref{eq:ass_smooth}}{\leq} \dotprod{\nabla f(x^k)}{x^{k+1} - x^k} + \frac{L_0 + L_1 \norms{\nabla f(x^k)}}{2} \norms{x^{k+1} - x^{k}}^2 \nonumber \\
    &= - \eta \dotprod{\nabla f(x^k)}{G(x^k, \bxi^k)} + \frac{\eta^2 (L_0 + L_1 \norms{\nabla f(x^k)})}{2} \norms{G(x^k, \bxi^k)}^2. \label{eq:NSGD_smooth}
\end{align}

Next, we consider 4 cases of the relation $\norms{\nabla f(x^k)}$ and $\norms{\nabla f(x^k, \bxi^k)}$ with respect to the hyperparameter $\lambda$.

\subsection{First case: \texorpdfstring{$\norms{\nabla f(x^k)} \geq \lambda$}{TEXT} and \texorpdfstring{$\norms{\nabla f(x^k, \bxi^k)} \geq \lambda$}{TEXT}}

    Let us evaluate first summand of \eqref{eq:NSGD_smooth} with $\alpha = \norms{\nabla f(x^k)}^{-1}$:
    \begin{align}
        - \eta \dotprod{\nabla f(x^k)}{G(x^k, \bxi^k)} &\overset{\eqref{eq:qudrat_raznosti}}{=} - \frac{\alpha \eta }{2} \norms{\nabla f(x^k)}^2 - \frac{\eta}{2 \alpha} \norms{G(x^k, \bxi^k)}^2 \nonumber\\
        & \quad \quad \quad+ \frac{\eta}{2 \alpha} \norms{G(x^k, \bxi^k) - \alpha \nabla f(x^k)}^2 \nonumber\\
        &= - \frac{\eta}{2 } \norms{\nabla f(x^k)} - \frac{\eta}{2 \alpha} \norms{G(x^k, \bxi^k)}^2 \nonumber\\
        & \quad \quad \quad+ \frac{\eta}{2 \lambda^2 \alpha} \norms{\lambda G(x^k, \bxi^k) - \lambda \alpha \nabla f(x^k)}^2 \nonumber\\
        &= - \frac{\eta}{2 } \norms{\nabla f(x^k)} - \frac{\eta}{2 \alpha} \norms{G(x^k, \bxi^k)}^2 \nonumber\\
        & \quad \quad \quad+ \frac{\eta}{2 \lambda^2 \alpha} \norms{\cliplam{\nabla f(x^k, \bxi^k)} - \cliplam{ \nabla f(x^k)}}^2 \nonumber
    \end{align}
    Using that clipping is a projection on onto a convex set, namely ball with radius $\lambda$, and thus is Lipshitz operator with Lipshitz constant $1$, we can obtain:
    \begin{align}
        - \eta \dotprod{\nabla f(x^k)}{\expect{G(x^k, \bxi^k)}} &\leq - \frac{\eta}{2} \norms{\nabla f(x^k)} - \frac{\eta}{2 \alpha} \expect{\norms{G(x^k, \bxi^k)}^2} \nonumber\\
        & \quad \quad \quad+ \frac{\eta }{2 \lambda^2 \alpha} \expect{\norms{\nabla f(x^k, \bxi^k) - \nabla f(x^k)}^2}.\label{eq:NSGD1_main}
    \end{align}

    \paragraph{In the case: $0 \leq \sigma \leq \frac{\lambda}{\sqrt{2}}$.} Using this in \eqref{eq:NSGD1_main}, we have the following with $\eta_k \leq \frac{\norms{\nabla f(x^k)}}{2(L_0 + L_1 \norms{\nabla f(x^k)})}$:
        \begin{align}
            \expect{f(x^{k+1})} - f(x^k) &\overset{\eqref{eq:NSGD_smooth}}{\leq}  - \eta \dotprod{\nabla f(x^k)}{\expect{G(x^k, \bxi^k)}} + \frac{\eta^2 (L_0 + L_1 \norms{\nabla f(x^k)})}{2} \expect{\norms{G(x^k, \bxi^k)}^2} \nonumber\\
            & \overset{\eqref{eq:NSGD1_main}}{\leq} - \frac{\eta}{2} \norms{\nabla f(x^k)} - \frac{\eta}{2 \alpha} \expect{\norms{G(x^k, \bxi^k)}^2} + \frac{\eta }{2 \lambda^2 \alpha} \expect{\norms{\nabla f(x^k, \bxi) - \nabla f(x^k)}^2}\nonumber\\
            &\quad\quad\quad+ \frac{\eta^2 (L_0 + L_1 \norms{\nabla f(x^k)})}{2} \expect{\norms{G(x^k, \bxi^k)}^2} \nonumber\\
            &= - \frac{\eta}{2} \norms{\nabla f(x^k)} + \frac{\eta }{2 \lambda^2 \alpha} \expect{\norms{\nabla f(x^k, \bxi^k) - \nabla f(x^k)}^2}\nonumber\\
            &\quad\quad\quad - \frac{\eta}{2} \expect{\norms{G(x^k, \bxi^k)}^2} \left( 1 - \frac{\eta (L_0 + L_1 \norms{\nabla f(x^k)})}{\norms{\nabla f(x^k)}} \right)\nonumber\\
            &\leq - \frac{\eta}{2} \norms{\nabla f(x^k)} + \frac{\eta \sigma^2}{2 \lambda^2 \alpha} \nonumber\\
            &\leq - \frac{\eta}{2} \norms{\nabla f(x^k)} + \frac{\eta}{4} \norms{\nabla f(x^k)}\nonumber\\
            &= - \frac{\eta}{4} \norms{\nabla f(x^k)}. \label{eq:NSGD1_case1_1}
        \end{align}
        
        The step size will be constant, depending on the hyperparameter $\lambda$:
        \begin{align*}
            \frac{\norms{\nabla f(x^k)}}{2\left(L_0 + L_1 \norms{\nabla f(x^k)}\right)} = \frac{1}{2\left(L_0\frac{1}{\norms{\nabla f(x^k)}} + L_1\right)} = \frac{\lambda}{2\left(L_0\frac{\lambda}{\norms{\nabla f(x^k)}} + L_1 \lambda\right)} \geq \frac{\lambda}{2\left(L_0 + L_1 \lambda\right)}.
        \end{align*}
        Thus, $\eta_k = \eta \leq \frac{\lambda}{2(L_0 + L_1 \lambda)}$.

        Using the convexity assumption of the function, we have the following:
        
            \begin{align*}
                f(x^k) - f^* &\leq \dotprod{\nabla f(x^k)}{x^k - x^*}\\
                &\overset{\eqref{eq:scalar_product_bound}}{\leq} \norms{\nabla f(x^k)} \norms{x^k - x^*}\\
                & \leq \norms{\nabla f(x^k)} \underbrace{\norms{x^0 - x^*}}_{R}.
            \end{align*}
            Hence we have:
            \begin{equation}
                \norms{\nabla f(x^k)} \geq \frac{f(x^k) - f^*}{R}. \label{eq:NSGD1_case1_2}
            \end{equation}
        
            Then substituting \eqref{eq:NSGD1_case1_2} into \eqref{eq:NSGD1_case1_1} we obtain:
            \begin{equation*}
                \expect{f(x^{k+1})} - f(x^k) \leq - \frac{\eta}{4} \norms{\nabla f(x^k)} \leq - \frac{\eta}{4R} (f(x^k) - f^*).
            \end{equation*}
            
            This inequality is equivalent to the trailing inequality:
            \begin{equation*}
                \expect{f(x^{k+1})} - f^* \leq \left( 1 - \frac{\eta}{4R} \right) \left( f(x^k) - f^* \right).
            \end{equation*}
            
                Then for $k = 0,1,2,..., N-1$ iterations that satisfy the conditions $\norms{\nabla f(x^k,\bxi^k)} \geq \sqrt{2} \sigma$ and $\norms{\nabla f(x^k)} \geq \sqrt{2} \sigma$ NSGD shows linear convergence:
                \begin{equation*}
                    \expect{f(x^{N})} - f^* \leq \left(1 - \frac{\eta}{4R} \right)^N (f(x^{0}) - f^*).
                \end{equation*}

    \paragraph{In the case: $\frac{\lambda}{\sqrt{2}} \leq \sigma $.} Using this in \eqref{eq:NSGD1_main}, we have the following with $\eta_k \leq \frac{\norms{\nabla f(x^k)}}{2(L_0 + L_1 \norms{\nabla f(x^k)})}$:
        \begin{align}
            \expect{f(x^{k+1})} - f(x^k) &\overset{\eqref{eq:NSGD_smooth}}{\leq}  - \eta \dotprod{\nabla f(x^k)}{\expect{G(x^k, \bxi^k)}} + \frac{\eta^2 (L_0 + L_1 \norms{\nabla f(x^k)})}{2} \expect{\norms{G(x^k, \bxi^k)}^2} \nonumber\\
            & \overset{\eqref{eq:NSGD1_main}}{\leq} - \frac{\eta}{2} \norms{\nabla f(x^k)} - \frac{\eta}{2 \alpha} \expect{\norms{G(x^k, \bxi^k)}^2} + \frac{\eta }{2 \lambda^2 \alpha} \expect{\norms{\nabla f(x^k, \bxi^k) - \nabla f(x^k)}^2}\nonumber\\
            &\quad\quad\quad+ \frac{\eta^2 (L_0 + L_1 \norms{\nabla f(x^k)})}{2} \expect{\norms{G(x^k, \bxi^k)}^2} \nonumber\\
            &= - \frac{\eta}{2} \norms{\nabla f(x^k)} + \frac{\eta }{2 \lambda^2 \alpha} \expect{\norms{\nabla f(x^k, \bxi^k) - \nabla f(x^k)}^2}\nonumber\\
            &\quad\quad\quad - \frac{\eta}{2} \expect{\norms{G(x^k, \bxi^k)}^2} \left( 1 - \frac{\eta (L_0 + L_1 \norms{\nabla f(x^k)})}{\norms{\nabla f(x^k)}} \right)\nonumber\\
            &\leq - \frac{\eta}{2} \norms{\nabla f(x^k)} + \frac{\eta \sigma^2}{2 \lambda^2 \alpha B} \nonumber\\
            &\leq - \frac{\eta}{2} \norms{\nabla f(x^k)} + \frac{\eta \sigma^2 M}{2 \lambda^2 B}. \label{eq:NSGD1_case2_1}
        \end{align}
        
        The step size will be constant, depending on the hyperparameter $\lambda$:
        \begin{align*}
            \frac{\norms{\nabla f(x^k)}}{2\left(L_0 + L_1 \norms{\nabla f(x^k)}\right)} = \frac{1}{2\left(L_0\frac{1}{\norms{\nabla f(x^k)}} + L_1\right)} = \frac{\lambda}{2\left(L_0\frac{\lambda}{\norms{\nabla f(x^k)}} + L_1 \lambda\right)} \geq \frac{\lambda}{2\left(L_0 + L_1 \lambda\right)}.
        \end{align*}
        Thus, $\eta_k = \eta \leq \frac{\lambda}{2(L_0 + L_1 \lambda)}$.
        
        Using the convexity assumption of the function, we have the following:

        \begin{align*}
            f(x^k) - f^* &\leq \dotprod{\nabla f(x^k)}{x^k - x^*}\\
            &\overset{\eqref{eq:scalar_product_bound}}{\leq} \norms{\nabla f(x^k)} \norms{x^k - x^*}\\
            & \leq \norms{\nabla f(x^k)} \underbrace{\norms{x^0 - x^*}}_{R}.
        \end{align*}
        Hence we have:
        \begin{equation}
            \norms{\nabla f(x^k)} \geq \frac{f(x^k) - f^*}{R}. \label{eq:NSGD1_case2_2}
        \end{equation}
    
        Then substituting \eqref{eq:NSGD1_case2_2} into \eqref{eq:NSGD1_case2_1} we obtain:
        \begin{equation*}
            \expect{f(x^{k+1})} - f(x^k) \leq - \frac{\eta}{2} \norms{\nabla f(x^k)} + \frac{\eta \sigma^2 M}{2 \lambda^2 B} \leq - \frac{\eta}{2 R} (f(x^k) - f^*) + \frac{\eta \sigma^2 M}{2 \lambda^2B}.
        \end{equation*}
        
        This inequality is equivalent to the trailing inequality:
        \begin{equation*}
            \expect{f(x^{k+1})} - f^* \leq \left( 1 - \frac{\eta}{2R} \right) \left( f(x^k) - f^* \right) + \frac{\eta \sigma^2 M}{2 \lambda^2 B}.
        \end{equation*}
    
        Then for $k = 0,1,2,..., N-1$ iterations that satisfy the conditions $\norms{\nabla f(x^k,\bxi^k)} \geq \lambda$ and $\norms{\nabla f(x^k)} \geq \lambda $ and $\sigma \geq \sqrt{2} \lambda$  NSGD shows linear convergence:
        \begin{equation*}
            \expect{f(x^{N})} - f^* \leq \left(1 - \frac{\eta}{2R} \right)^N (f(x^{0}) - f^*) + \frac{\sigma^2 M R}{\lambda^2 B}.
        \end{equation*}

\subsection{Second case: \texorpdfstring{$\norms{\nabla f(x^k)} \leq \lambda$}{TEXT} and \texorpdfstring{$\norms{\nabla f(x^k, \bxi^k)} \geq \lambda$}{TEXT}}

    Let us evaluate first summand of \eqref{eq:NSGD_smooth} with $\alpha = \lambda^{-1}$:
    \begin{align}
        - \eta \dotprod{\nabla f(x^k)}{G(x^k, \bxi^k)} &\overset{\eqref{eq:qudrat_raznosti}}{=} - \frac{\alpha \eta }{2} \norms{\nabla f(x^k)}^2 - \frac{\eta}{2 \alpha} \norms{G(x^k, \bxi^k)}^2 \nonumber\\
        & \quad \quad \quad+ \frac{\eta}{2 \alpha} \norms{G(x^k, \bxi^k) - \alpha \nabla f(x^k)}^2 \nonumber\\
        &\leq - \frac{\eta}{2} \norms{\nabla f(x^k)} - \frac{\eta}{2 \alpha} \norms{G(x^k, \bxi^k)}^2 \nonumber\\
        & \quad \quad \quad+ \frac{\eta}{2 \lambda} \norms{\lambda G(x^k, \bxi^k) - \nabla f(x^k)}^2 \nonumber\\
        &= - \frac{\eta}{2} \norms{\nabla f(x^k)} - \frac{\eta}{2 \alpha} \norms{G(x^k, \bxi^k)}^2 \nonumber\\
        & \quad \quad \quad+ \frac{\eta}{2 \lambda} \norms{\cliplam{\nabla f(x^k, \bxi^k)} - \cliplam{ \nabla f(x^k)}}^2 \nonumber
    \end{align}
    Using that clipping is a projection on onto a convex set, namely ball with radius $\lambda$, and thus is Lipshitz operator with Lipshitz constant $1$, we can obtain:
    \begin{align}
        - \eta \dotprod{\nabla f(x^k)}{\expect{G(x^k, \bxi^k)}} &\leq - \frac{\eta}{2} \norms{\nabla f(x^k)} - \frac{\eta}{2 \alpha} \expect{\norms{G(x^k, \bxi^k)}^2} \nonumber\\
        & \quad \quad \quad+ \frac{\eta }{2 \lambda} \expect{\norms{\nabla f(x^k, \bxi^k) - \nabla f(x^k)}^2}. \label{eq:NSGD2_main}
    \end{align}
     Using this, we have the following with $\eta_k \leq \frac{\norms{\nabla f(x^k)}}{2(L_0 + L_1 \norms{\nabla f(x^k)})}$:
    \begin{align}
        \expect{f(x^{k+1})} - f(x^k) &\overset{\eqref{eq:NSGD_smooth}}{\leq}  - \eta \dotprod{\nabla f(x^k)}{\expect{G(x^k, \bxi^k)}} + \frac{\eta^2 (L_0 + L_1 \norms{\nabla f(x^k)})}{2} \expect{\norms{G(x^k, \bxi^k)}^2} \nonumber\\
        & \overset{\eqref{eq:NSGD2_main}}{\leq} - \frac{\eta}{2} \norms{\nabla f(x^k)} - \frac{\eta}{2 \alpha} \expect{\norms{G(x^k, \bxi^k)}^2} + \frac{\eta }{2 \lambda} \expect{\norms{\nabla f(x^k, \bxi^k) - \nabla f(x^k)}^2}\nonumber\\
        &\quad\quad\quad+ \frac{\eta^2 (L_0 + L_1 \norms{\nabla f(x^k)})}{2} \expect{\norms{G(x^k, \bxi^k)}^2} \nonumber\\
        &= - \frac{\eta}{2} \norms{\nabla f(x^k)} + \frac{\eta }{2 \lambda} \expect{\norms{\nabla f(x^k, \bxi^k) - \nabla f(x^k)}^2}\nonumber\\
        &\quad\quad\quad - \frac{\eta}{2} \expect{\norms{G(x^k, \bxi^k)}^2} \left( 1 - \frac{\eta (L_0 + L_1 \norms{\nabla f(x^k)})}{\norms{\nabla f(x^k)}} \right)\nonumber\\
        &\leq - \frac{\eta}{2} \norms{\nabla f(x^k)} + \frac{\eta \sigma^2}{2 \lambda B} \nonumber\\
        &\leq - \frac{\eta}{2} \norms{\nabla f(x^k)} + \frac{\eta \sigma^2}{2 \lambda B}. \label{eq:NSGD2_case1_1}
    \end{align}
        
    The step size will be constant, depending on the hyperparameter $\lambda$:
    \begin{align*}
        \frac{\norms{\nabla f(x^k)}}{2\left(L_0 + L_1 \norms{\nabla f(x^k)}\right)} = \frac{1}{2\left(L_0\frac{1}{\norms{\nabla f(x^k)}} + L_1\right)} = \frac{\lambda}{2\left(L_0\frac{\lambda}{\norms{\nabla f(x^k)}} + L_1 \lambda\right)} \geq \frac{\lambda}{2\left(L_0 + L_1 \lambda\right)}.
    \end{align*}
    Thus, $\eta_k = \eta \leq \frac{\lambda}{2(L_0 + L_1 \lambda)}$.
        
    Using the convexity assumption of the function, we have the following:

    \begin{align*}
        f(x^k) - f^* &\leq \dotprod{\nabla f(x^k)}{x^k - x^*}\\
        &\overset{\eqref{eq:scalar_product_bound}}{\leq} \norms{\nabla f(x^k)} \norms{x^k - x^*}\\
        & \leq \norms{\nabla f(x^k)} \underbrace{\norms{x^0 - x^*}}_{R}.
    \end{align*}
    Hence we have:
    \begin{equation}
        \norms{\nabla f(x^k)} \geq \frac{f(x^k) - f^*}{R}. \label{eq:NSGD2_case1_2}
    \end{equation}
    
    Then substituting \eqref{eq:NSGD2_case1_2} into \eqref{eq:NSGD2_case1_1} we obtain:
    \begin{equation*}
        \expect{f(x^{k+1})} - f(x^k) \leq - \frac{\eta}{2} \norms{\nabla f(x^k)} + \frac{\eta \sigma^2}{2 \lambda B} \leq - \frac{\eta}{2 R} (f(x^k) - f^*) + \frac{\eta \sigma^2}{2 \lambda B}.
    \end{equation*}
        
    This inequality is equivalent to the trailing inequality:
    \begin{equation*}
        \expect{f(x^{k+1})} - f^* \leq \left( 1 - \frac{\eta}{2R} \right) \left( f(x^k) - f^* \right) + \frac{\eta \sigma^2}{2 \lambda B}.
    \end{equation*}
    
    Then for $k = 0,1,2,..., N-1$ iterations that satisfy the conditions $\norms{\nabla f(x^k)} \leq \lambda$ and $\norms{\nabla f(x^k, \bxi^k)} \geq \lambda $  NSGD shows linear convergence:
    \begin{equation*}
        \expect{f(x^{N})} - f^* \leq \left(1 - \frac{\eta}{2R} \right)^N (f(x^{0}) - f^*) + \frac{\sigma^2 R}{\lambda B}.
    \end{equation*}

\subsection{Third case: \texorpdfstring{$\norms{\nabla f(x^k)} \leq \lambda$}{TEXT} and \texorpdfstring{$\norms{\nabla f(x^k, \bxi^k)} \leq \lambda$}{TEXT}} 

    Using this in \eqref{eq:NSGD_smooth}, we have the following with $\eta_k \leq \frac{\norms{\nabla f(x^k)}}{2(L_0 + L_1 \norms{\nabla f(x^k)})}$ and $\alpha = \norms{\nabla f(x^k)}^{-1}$:
    \begin{align}
        \expect{f(x^{k+1})} - f(x^k) &\overset{\eqref{eq:NSGD_smooth}}{\leq}  - \eta \dotprod{\nabla f(x^k)}{\expect{G(x^k, \bxi^k)}} \nonumber\\
        & \quad \quad \quad+ \frac{\eta^2 (L_0 + L_1 \norms{\nabla f(x^k)})}{2} \expect{\norms{G(x^k, \bxi^k)}^2} \nonumber\\
        & = - \frac{\eta \alpha}{2} \norms{\nabla f(x^k)}^2 - \frac{\eta}{2 \alpha} \expect{\norms{G(x^k, \bxi^k)}^2} \nonumber\\
        & \quad \quad \quad+ \frac{\eta}{2\alpha} \expect{\norms{G(x^k, \bxi^k) - \alpha \nabla f(x^k)}^2}\nonumber\\
        &\quad\quad\quad+ \frac{\eta^2 (L_0 + L_1 \norms{\nabla f(x^k)})}{2} \expect{\norms{G(x^k, \bxi^k)}^2} \nonumber\\
        &= - \frac{\eta}{2} \norms{\nabla f(x^k)} + \frac{\eta}{2\alpha} \expect{\norms{G(x^k, \bxi^k) - \alpha \nabla f(x^k)}^2}\nonumber\\
        &\quad\quad\quad - \frac{\eta}{2} \expect{\norms{G(x^k, \bxi^k)}^2} \left( 1 - \frac{\eta (L_0 + L_1 \norms{\nabla f(x^k)})}{\norms{\nabla f(x^k)}} \right)\nonumber\\
        &\leq - \frac{\eta}{2} \norms{\nabla f(x^k)} + \frac{\eta}{2\alpha} \expect{\norms{G(x^k, \bxi^k) - \alpha \nabla f(x^k)}^2} \nonumber\\
        &\leq - \frac{\eta}{2} \norms{\nabla f(x^k)} + \frac{\eta}{\alpha} \expect{\norms{G(x^k, \bxi^k)}^2 +\norms{\alpha \nabla f(x^k)}^2}\nonumber\\
        &= - \frac{\eta}{2} \norms{\nabla f(x^k)} + \frac{\eta}{\alpha} \expect{\norms{\frac{\nabla f(x^k, \bxi^k)}{\norms{\nabla f(x^k, \bxi^k)}}}^2 +\norms{\frac{\nabla f(x^k)}{\norms{\nabla f(x^k)}}}^2}\nonumber\\
        &= - \frac{\eta}{2} \norms{\nabla f(x^k)} + \frac{2 \eta \lambda \norms{\nabla f(x^k)}}{\lambda}\nonumber\\
        &\leq - \frac{\eta}{2} \norms{\nabla f(x^k)} + 2 \eta \lambda. \label{eq:NSGD3_case1_1}
    \end{align}
    
    The step size will be constant, depending on the hyperparameter $\lambda$:
    \begin{align*}
        \frac{\norms{\nabla f(x^k)}}{2\left(L_0 + L_1 \norms{\nabla f(x^k)}\right)} = \frac{1}{2\left(L_0\frac{1}{\norms{\nabla f(x^k)}} + L_1\right)} = \frac{\lambda}{2\left(L_0\frac{\lambda}{\norms{\nabla f(x^k)}} + L_1 \lambda\right)} \geq \frac{\lambda}{2\left(L_0 + L_1 \lambda\right)}.
    \end{align*}
    Thus, $\eta_k = \eta \leq \frac{\lambda}{2(L_0 + L_1 \lambda)}$.
    
    Using the convexity assumption of the function, we have the following:

    \begin{align*}
        f(x^k) - f^* &\leq \dotprod{\nabla f(x^k)}{x^k - x^*}\\
        &\overset{\eqref{eq:scalar_product_bound}}{\leq} \norms{\nabla f(x^k)} \norms{x^k - x^*}\\
        & \leq \norms{\nabla f(x^k)} \underbrace{\norms{x^0 - x^*}}_{R}.
    \end{align*}
    Hence we have:
    \begin{equation}
        \norms{\nabla f(x^k)} \geq \frac{f(x^k) - f^*}{R}. \label{eq:NSGD3_case1_2}
    \end{equation}

    Then substituting \eqref{eq:NSGD3_case1_2} into \eqref{eq:NSGD3_case1_1} we obtain:
    \begin{equation*}
        \expect{f(x^{k+1})} - f(x^k) \leq - \frac{\eta}{2} \norms{\nabla f(x^k)} + 2 \eta \lambda \leq - \frac{\eta}{2 R} (f(x^k) - f^*) + 2 \eta \lambda.
    \end{equation*}
    
    This inequality is equivalent to the trailing inequality:
    \begin{equation*}
        \expect{f(x^{k+1})} - f^* \leq \left( 1 - \frac{\eta}{2R} \right) \left( f(x^k) - f^* \right) + 2 \eta \lambda.
    \end{equation*}
    
        Then for $k = 0,1,2,..., N-1$ iterations that satisfy the conditions $\norms{\nabla f(x^k)} \leq \lambda$  NSGD shows linear convergence:
        \begin{equation*}
            \expect{f(x^{N})} - f^* \leq \left(1 - \frac{\eta}{2R} \right)^N (f(x^{0}) - f^*) + \lambda R.
        \end{equation*}

\subsection{Fourth case: \texorpdfstring{$\norms{\nabla f(x^k)} \geq \lambda$}{TEXT} and \texorpdfstring{$\norms{\nabla f(x^k, \bxi^k)} \leq \lambda$}{TEXT}}

    Using this in \eqref{eq:NSGD_smooth}, we have the following with $\eta_k \leq \frac{\norms{\nabla f(x^k)}}{2(L_0 + L_1 \norms{\nabla f(x^k)})}$ and $\alpha = \lambda^{-1}$:
    \begin{align}
        \expect{f(x^{k+1})} - f(x^k) &\overset{\eqref{eq:NSGD_smooth}}{\leq}  - \eta \dotprod{\nabla f(x^k)}{\expect{G(x^k, \bxi^k)}} \nonumber\\
        & \quad \quad \quad+ \frac{\eta^2 (L_0 + L_1 \norms{\nabla f(x^k)})}{2} \expect{\norms{G(x^k, \bxi^k)}^2} \nonumber\\
        & = - \frac{\eta \alpha}{2} \norms{\nabla f(x^k)}^2 - \frac{\eta}{2 \alpha} \norms{\expect{G(x^k, \bxi^k)}}^2 \nonumber\\
        & \quad \quad \quad+ \frac{\eta}{2\alpha} \norms{\expect{G(x^k, \bxi^k)} - \alpha \nabla f(x^k)}^2\nonumber\\
        &\quad\quad\quad+ \frac{\eta^2 (L_0 + L_1 \norms{\nabla f(x^k)})}{2} \expect{\norms{G(x^k, \bxi^k)}^2} \nonumber\\
        &= - \frac{\eta}{2 \lambda} \norms{\nabla f(x^k)}^2 + \frac{\eta}{2 \lambda} \norms{ \expect{\lambda G(x^k, \bxi^k)} - \nabla f(x^k)}^2 \nonumber\\
        & \quad \quad \quad+ \frac{\eta^2 (L_0 + L_1 \norms{\nabla f(x^k)})}{2}\nonumber\\
        &= - \frac{\eta}{2 \lambda} \norms{\nabla f(x^k)}^2 + \frac{\eta}{2 \lambda} \norms{ \expect{\frac{\lambda \nabla f(x^k, \bxi^k)}{\norms{\nabla f(x^k, \bxi^k)}} - \nabla f(x^k,\bxi^k)}}^2 \nonumber\\
        & \quad \quad \quad+ \frac{\eta^2 (L_0 + L_1 \norms{\nabla f(x^k)})}{2}\nonumber\\
        &= - \frac{\eta}{2 \lambda} \norms{\nabla f(x^k)}^2 + \frac{\eta}{2 \lambda} \norms{ \expect{\left(\frac{\lambda}{\norms{\nabla f(x^k, \bxi^k)}} - 1\right) \nabla f(x^k,\bxi^k)}}^2 \nonumber\\
        & \quad \quad \quad+ \frac{\eta^2 (L_0 + L_1 \norms{\nabla f(x^k)})}{2}\nonumber\\
        &\leq - \frac{\eta}{2 \lambda} \norms{\nabla f(x^k)}^2 + \frac{\eta}{2 \lambda} \expect{\left(\frac{\lambda}{\norms{\nabla f(x^k, \bxi^k)}} - 1\right)^2 \norms{ \nabla f(x^k,\bxi^k)}^2} \nonumber\\
        & \quad \quad \quad+ \frac{\eta^2 (L_0 + L_1 \norms{\nabla f(x^k)})}{2}\nonumber\\
        &\leq - \frac{\eta}{2 \lambda} \norms{\nabla f(x^k)}^2 + \frac{\eta}{2 \lambda} \expect{\frac{\lambda^2}{\norms{\nabla f(x^k, \bxi^k)}^2} \norms{ \nabla f(x^k,\bxi^k)}^2} \nonumber\\
        & \quad \quad \quad+ \frac{\eta^2 (L_0 + L_1 \norms{\nabla f(x^k)})}{2}\nonumber\\
        &= - \frac{\eta}{2 \lambda} \norms{\nabla f(x^k)}^2 + \frac{\eta^2 (L_0 + L_1 \norms{\nabla f(x^k)})}{2} + \frac{\eta \lambda}{2} \nonumber\\
        &\leq - \frac{\eta}{2} \norms{\nabla f(x^k)} + \frac{\eta^2 (L_0 + L_1 \norms{\nabla f(x^k)})}{2} + \frac{\eta \lambda}{2} \nonumber\\
        &= - \frac{\eta}{2} \norms{\nabla f(x^k)} \left( 1 - \frac{\eta (L_0 + L_1 \norms{\nabla f(x^k)})}{\norms{\nabla f(x^k)}} \right)  + \frac{\eta \lambda}{2} \nonumber\\
        &\leq - \frac{\eta}{4} \norms{\nabla f(x^k)} + \frac{\eta \lambda}{2}.
        \label{eq:NSGD4_case1_1}
    \end{align}
    
    The step size will be constant, depending on the hyperparameter $\lambda$:
    \begin{align*}
        \frac{\norms{\nabla f(x^k)}}{2\left(L_0 + L_1 \norms{\nabla f(x^k)}\right)} = \frac{1}{2\left(L_0\frac{1}{\norms{\nabla f(x^k)}} + L_1\right)} = \frac{\lambda}{2\left(L_0\frac{\lambda}{\norms{\nabla f(x^k)}} + L_1 \lambda\right)} \geq \frac{\lambda}{2\left(L_0 + L_1 \lambda\right)}.
    \end{align*}
    Thus, $\eta_k = \eta \leq \frac{\lambda}{2(L_0 + L_1 \lambda)}$.
    
    Using the convexity assumption of the function, we have the following:

    \begin{align*}
        f(x^k) - f^* &\leq \dotprod{\nabla f(x^k)}{x^k - x^*}\\
        &\overset{\eqref{eq:scalar_product_bound}}{\leq} \norms{\nabla f(x^k)} \norms{x^k - x^*}\\
        & \leq \norms{\nabla f(x^k)} \underbrace{\norms{x^0 - x^*}}_{R}.
    \end{align*}
    Hence we have:
    \begin{equation}
        \norms{\nabla f(x^k)} \geq \frac{f(x^k) - f^*}{R}. \label{eq:NSGD4_case1_2}
    \end{equation}

    Then substituting \eqref{eq:NSGD4_case1_2} into \eqref{eq:NSGD4_case1_1} we obtain:
    \begin{equation*}
        \expect{f(x^{k+1})} - f(x^k) \leq - \frac{\eta}{4} \norms{\nabla f(x^k)} +\frac{\eta \lambda}{2}  \leq - \frac{\eta}{4 R} (f(x^k) - f^*) +\frac{\eta \lambda}{2}.
    \end{equation*}
    
    This inequality is equivalent to the trailing inequality:
    \begin{equation*}
        \expect{f(x^{k+1})} - f^* \leq \left( 1 - \frac{\eta}{4R} \right) \left( f(x^k) - f^* \right) + \frac{\eta \lambda}{2}.
    \end{equation*}
    
        Then for $k = 0,1,2,..., N-1$ iterations that satisfy the conditions $\norms{\nabla f(x^k)} \geq \lambda$ and $\norms{\nabla f(x^k, \bxi^k)} \leq \lambda$  NSGD shows linear convergence:
        \begin{equation*}
            \expect{f(x^{N})} - f^* \leq \left(1 - \frac{\eta}{4R} \right)^N (f(x^{0}) - f^*) + 2 \lambda R.
        \end{equation*}

Combining all the cases considered, we obtain the convergence rate of NSGD: 
\begin{equation*}
    \expect{f(x^{N})} - f^* \lesssim  \left( 1 - \frac{\eta}{R} \right)^{N} (f(x^0) - f^*)) + \frac{\sigma^2 M R}{B \lambda^2} + \lambda R.
\end{equation*}

\section{Zero-Order Clipped Stochastic Gradient Descent Method}
\label{app:ZO_ClipSGD}

This section consists of two parts: 1) a generalization of the convergence result of ClipSGD (Algorithm~\ref{algo:clipSGD}) to the biased gradient oracle $\gg(x^k, \bxi^k) = \nabla f(x^k, \bxi^k) + \bb(x^k)$, where $\bb(x^k)$ is biased bounded by $\zeta \geq 0: \norms{\bb(x^k)}\leq \zeta$; 2) deriving convergence estimates of ZO-ClipSGD directly. 

\subsection{Biased Clipped Stochastic Gradient Descent Method (Proof of the Lemma~\ref{lem:biased_ZO_ClipSGD})}
We start by using $(L_0,L_1)$-smoothness (see Assumption~\ref{ass:L0_L1_smooth}):
\begin{align}
    f(x^{k+1}) - f(x^k) &\overset{\eqref{eq:ass_smooth}}{\leq} \dotprod{\nabla f(x^k)}{x^{k+1} - x^k} + \frac{L_0 + L_1 \norms{\nabla f(x^k)}}{2} \norms{x^{k+1} - x^{k}}^2 \nonumber \\
    &= - \eta \dotprod{\nabla f(x^k)}{\clip{\gg(x^k, \bxi^k)}} \nonumber\\
    &\quad \quad \quad+ \frac{\eta^2 (L_0 + L_1 \norms{\nabla f(x^k)})}{2} \norms{\clip{\gg(x^k, \bxi^k)}}^2. \label{eq:Biased_clipSGD_smooth}
\end{align}

Next, we consider three cases depending on the gradient norm: $\norms{\nabla f(x^k)} \geq c$ -- the full gradient is clipped and ${\frac{c}{3} \leq \norms{\nabla f(x^k)} \leq c}$ and $\norms{\nabla f(x^k)} \leq \frac{c}{3}$ -- the full gradient is not clipped.

\subsubsection{First case: \texorpdfstring{$\norms{\nabla f(x^k)} \geq c$}{TEXT}}

    In this case $\alpha \nabla f(x^k) = \clip{\nabla f(x^k)}$ with $\alpha = \min \left\{1, \frac{c}{\norms{\nabla f(x^k)}}\right\} = \frac{c}{\norms{\nabla f(x^k)}}$, therefore we have the following
    \begin{align}
        - \eta \dotprod{\nabla f(x^k)}{\clip{\gg(x^k, \bxi^k)}} &\overset{\eqref{eq:qudrat_raznosti}}{=} - \frac{\alpha \eta}{2} \norms{\nabla f(x^k)}^2 - \frac{\eta}{2 \alpha} \norms{\clip{\gg(x^k, \bxi^k)}}^2 \nonumber\\
    &\quad \quad \quad+ \frac{\eta}{2 \alpha} \norms{\clip{\gg(x^k, \bxi^k)} - \alpha \nabla f(x^k)}^2 \nonumber\\
        &= - \frac{\alpha \eta}{2} \norms{\nabla f(x^k)}^2 - \frac{\eta}{2 \alpha} \norms{\clip{\gg(x^k, \bxi^k)}}^2 \nonumber\\
    &\quad \quad \quad+ \frac{\eta}{2 \alpha} \norms{\clip{\gg(x^k, \bxi^k)} - \clip{\nabla f(x^k)}}^2 \nonumber\\
        &= - \frac{c \eta}{2} \norms{\nabla f(x^k)} - \frac{\eta}{2 \alpha} \norms{\clip{\gg(x^k, \bxi^k)}}^2 \nonumber\\
    &\quad \quad \quad+ \frac{\eta}{2 \alpha} \norms{\clip{\gg(x^k, \bxi^k)} - \clip{\nabla f(x^k)}}^2. \nonumber
    \end{align}
    Using that clipping is a projection on onto a convex set, namely ball with radius $c$, and thus is Lipshitz operator with Lipshitz constant $1$, we can obtain:
    \begin{align}
        - \eta \dotprod{\nabla f(x^k)}{\expect{\clip{\gg(x^k, \bxi^k)}}} &\leq - \frac{c \eta}{2} \norms{\nabla f(x^k)} - \frac{\eta}{2 \alpha} \expect{\norms{\clip{\gg(x^k, \bxi^k)}}^2} \nonumber\\
    &\quad \quad \quad+ \frac{\eta}{2 \alpha} \expect{\norms{\gg(x^k, \bxi^k) - \nabla f(x^k)}^2} \nonumber\\
        &\overset{\eqref{eq:Variance_decomposition}}{=} - \frac{c \eta}{2} \norms{\nabla f(x^k)} - \frac{\eta}{2 \alpha} \expect{\norms{\clip{\gg(x^k, \bxi^k)}}^2} \nonumber\\
        &\quad\quad\quad + \frac{\eta}{2 \alpha} \expect{\norms{\gg(x^k, \bxi^k) - \expect{\gg(x^k, \bxi^k)}}^2} \nonumber\\
    &\quad \quad \quad+ \frac{\eta}{2 \alpha} \norms{\bb(x^k)} \nonumber\\
        &\leq  - \frac{c \eta}{2} \norms{\nabla f(x^k)} - \frac{\eta}{2 \alpha} \expect{\norms{\clip{\nabla f(x^k, \bxi^k)}}^2} \nonumber\\
    &\quad \quad \quad+ \frac{\eta \sigma^2 M}{2 c B} + \frac{\eta \norms{\nabla f(x^k)} \zeta^2}{2 c }. \label{eq:Biased_clipSGD_case1_noise}
    \end{align}
    We now consider the cases depending on the relation between $c$ and $\zeta$:

    \paragraph{\fbox{In the case $c \geq \sqrt{2} \zeta$}}
    
    We have in \eqref{eq:Biased_clipSGD_case1_noise}:
    \begin{align*}
        - \eta \dotprod{\nabla f(x^k)}{\expect{\clip{\gg(x^k, \bxi^k)}}} &\overset{\eqref{eq:Biased_clipSGD_case1_noise}}{\leq} - \frac{c \eta}{2} \norms{\nabla f(x^k)} - \frac{\eta}{2 \alpha} \expect{\norms{\clip{\gg(x^k, \bxi^k)}}^2} \nonumber\\
    &\quad \quad \quad+ \frac{\eta \sigma^2 M}{2 c B} + \frac{\eta \norms{\nabla f(x^k)} \zeta^2}{2 c} \\
        &= - \frac{\eta}{2 \alpha} \expect{\norms{\clip{\gg(x^k, \bxi^k)}}^2} \nonumber\\
    &\quad \quad \quad- \frac{c \eta}{2} \norms{\nabla f(x^k)} \left( 1 - \frac{\zeta^2}{c^2} \right) + \frac{\eta \sigma^2 M}{2 c B}\\
        &\leq - \frac{\eta}{2 \alpha} \expect{\norms{\clip{\gg(x^k, \bxi^k)}}^2} - \frac{c \eta}{4} \norms{\nabla f(x^k)} + \frac{\eta \sigma^2 M}{2 c B} \\ 
        &= - \frac{\eta \norms{\nabla f(x^k)}}{2 c} \expect{\norms{\clip{\nabla f(x^k, \bxi^k)}}^2} \nonumber\\
    &\quad \quad \quad- \frac{c \eta}{4} \norms{\nabla f(x^k)} + \frac{\eta \sigma^2 M}{2 c B}.
    \end{align*}
    Plugging this into \eqref{eq:Biased_clipSGD_smooth} and choosing $\eta \leq \frac{1}{4(L_0 + L_1 c)}$ we have:
    \begin{align}
        \expect{\nabla f(x^{k+1})} - f(x^k) &\overset{\eqref{eq:clipSGD_smooth}}{\leq} - \frac{\eta \norms{\nabla f(x^k)}}{2 c} \expect{\norms{\clip{\gg(x^k, \bxi^k)}}^2} - \frac{c \eta}{4} \norms{\nabla f(x^k)}  + \frac{\eta \sigma^2 M}{2 c B}\nonumber\\
        & \quad \quad + \frac{\eta^2 (L_0 + L_1 \norms{\nabla f(x^k)})}{2} \expect{\norms{\clip{\gg(x^k, \bxi^k)}}^2}\nonumber\\
        &= - \frac{\eta \norms{\nabla f(x^k)}}{2 c} \expect{\norms{\clip{\gg(x^k, \bxi^k)}}^2} (1 - \eta L_1 c) - \frac{c \eta}{4} \norms{\nabla f(x^k)} \nonumber\\
        & \quad \quad + \frac{\eta^2 L_0}{2}\expect{\norms{\clip{\gg(x^k, \bxi^k)}}^2} + \frac{\eta \sigma^2 M}{2 c B}\nonumber\\
        &\leq  - \frac{c \eta}{4} \norms{\nabla f(x^k)} - \frac{\eta}{2} \expect{\norms{\clip{\gg(x^k, \bxi^k)}}^2} \left(1 - \eta( L_0 + L_1 c)\right) \nonumber\\
    &\quad \quad \quad+ \frac{\eta \sigma^2 M}{2 c B}\nonumber\\
        &\leq  - \frac{c \eta}{4} \norms{\nabla f(x^k)} + \frac{\eta \sigma^2 M}{2 c B}. \label{eq:Biased_clipSGD_case1_withoutnoise_1}
    \end{align}
    Using the convexity assumption of the function, we have the following:
    \begin{align*}
        f(x^k) - f^* &\leq \dotprod{\nabla f(x^k)}{x^k - x^*} \overset{\eqref{eq:scalar_product_bound}}{\leq} \norms{\nabla f(x^k)} \norms{x^k - x^*}  \leq \norms{\nabla f(x^k)} \underbrace{\norms{x^0 - x^*}}_{R}.
    \end{align*}
    Hence we have:
    \begin{equation}
        \norms{\nabla f(x^k)} \geq \frac{f(x^k) - f^*}{R}. \label{eq:Biased_clipSGD_case1_withoutnoise_2}
    \end{equation}

    Then substituting \eqref{eq:Biased_clipSGD_case1_withoutnoise_2} into \eqref{eq:Biased_clipSGD_case1_withoutnoise_1} we obtain:
    \begin{equation*}
        \expect{f(x^{k+1})} - f(x^k) \leq - \frac{\eta c}{4} \norms{\nabla f(x^k)} + \frac{\eta \sigma^2 M}{2 c B} \leq - \frac{\eta c}{4 R} (f(x^k) - f^*) + \frac{\eta \sigma^2 M}{2 c B}.
    \end{equation*}
    
    This inequality is equivalent to the trailing inequality:
    \begin{equation*}
        \expect{f(x^{k+1})} - f^* \leq \left( 1 - \frac{\eta c}{4 R} \right) \left( f(x^k) - f^* \right) + \frac{\eta \sigma^2 M}{2 c B}.
    \end{equation*}
     Then for $k = 0,1,2,..., N-1$ iterations that satisfy the conditions $\norms{\nabla f(x^k)} \geq c \geq \sqrt{2} \zeta$, then ClipSGD with biased gradient oracle has linear convergence
    \begin{equation*}
        \expect{f(x^{N})} - f^* \leq \left( 1 - \frac{\eta}{2 R} \right)^{N} \left( f(x^0) - f^* \right) + \frac{\sigma^2 M R}{c B}.
    \end{equation*}

    \paragraph{\fbox{In the case $c \leq \sqrt{2} \zeta$}}

    We have in \eqref{eq:Biased_clipSGD_case1_noise}:
    \begin{align*}
        - \eta \dotprod{\nabla f(x^k)}{\expect{\clip{\gg(x^k, \bxi^k)}}} &\overset{\eqref{eq:Biased_clipSGD_case1_noise}}{\leq} - \frac{c \eta}{2} \norms{\nabla f(x^k)} - \frac{\eta}{2 \alpha} \expect{\norms{\clip{\gg(x^k, \bxi^k)}}^2} \nonumber\\
    &\quad \quad \quad+ \frac{\eta \sigma^2 M}{2 c B} + \frac{\eta \norms{\nabla f(x^k)} \zeta^2}{2 c} \\
        &= - \frac{c \eta}{2} \norms{\nabla f(x^k)} - \frac{\eta}{2 \alpha} \expect{\norms{\clip{\gg(x^k, \bxi^k)}}^2} \nonumber\\
    &\quad \quad \quad+ \frac{\eta M}{2 c} \left( \frac{\sigma^2}{B} + \zeta^2 \right).
    \end{align*}
    Plugging this into \eqref{eq:Biased_clipSGD_smooth} and choosing $\eta \leq \frac{1}{4(L_0 + L_1 c)}$ we have:
    \begin{align}
        \expect{f(x^{k+1})} - f(x^k) &\overset{\eqref{eq:Biased_clipSGD_smooth}}{\leq} - \frac{c \eta}{2} \norms{\nabla f(x^k)} - \frac{\eta \norms{\nabla f(x^k)}}{2 c} \expect{\norms{\clip{\gg(x^k, \bxi^k)}}^2}  \nonumber\\
        & \quad \quad + \frac{\eta^2 (L_0 + L_1 \norms{\nabla f(x^k)})}{2} \expect{\norms{\clip{\gg(x^k, \bxi^k)}}^2} + \frac{\eta M}{2 c} \left( \frac{\sigma^2}{B} + \zeta^2 \right) \nonumber\\
        &= - \frac{c \eta}{2} \norms{\nabla f(x^k)} - \frac{\eta \norms{\nabla f(x^k)}}{2 c} \expect{\norms{\clip{\gg(x^k, \bxi^k)}}^2} (1 - \eta L_1 c) \nonumber\\
        & \quad \quad + \frac{\eta^2 L_0}{2}\expect{\norms{\clip{\gg(x^k, \bxi^k)}}^2} + \frac{\eta M}{2 c} \left( \frac{\sigma^2}{B} + \zeta^2 \right) \nonumber\\
        &\leq  - \frac{c \eta}{2} \norms{\nabla f(x^k)} - \frac{\eta}{2} \expect{\norms{\clip{\gg(x^k, \bxi^k)}}^2} \left(1 - \eta( L_0 + L_1 c)\right) \nonumber\\
    &\quad \quad \quad+ \frac{\eta M}{2 c} \left( \frac{\sigma^2}{B} + \zeta^2 \right)\nonumber\\
        &\leq  - \frac{c \eta}{2} \norms{\nabla f(x^k)} + \frac{\eta M}{2 c} \left( \frac{\sigma^2}{B} + \zeta^2 \right). \label{eq:Biased_clipSGD_case1_with_noise_1}
    \end{align}
    Using the convexity assumption of the function, we have the following:
    \begin{align*}
        f(x^k) - f^* &\leq \dotprod{\nabla f(x^k)}{x^k - x^*} \overset{\eqref{eq:scalar_product_bound}}{\leq} \norms{\nabla f(x^k)} \norms{x^k - x^*} \leq \norms{\nabla f(x^k)} \underbrace{\norms{x^0 - x^*}}_{R}.
    \end{align*}
    Hence we have:
    \begin{equation}
        \norms{\nabla f(x^k)} \geq \frac{f(x^k) - f^*}{R}. \label{eq:Biased_clipSGD_case1_with_noise_2}
    \end{equation}

    Then substituting \eqref{eq:Biased_clipSGD_case1_with_noise_2} into \eqref{eq:Biased_clipSGD_case1_with_noise_1} we obtain:
    \begin{equation*}
        \expect{f(x^{k+1})} - f(x^k) \leq - \frac{\eta c}{2} \norms{\nabla f(x^k)} + \frac{\eta M}{2 c} \left( \frac{\sigma^2}{B} + \zeta^2 \right) \leq - \frac{\eta c}{2 R} (f(x^k) - f^*) + \frac{\eta M}{2 c} \left( \frac{\sigma^2}{B} + \zeta^2 \right).
    \end{equation*}
    
    This inequality is equivalent to the trailing inequality:
    \begin{equation*}
        \expect{f(x^{k+1})} - f^* \leq \left( 1 - \frac{\eta c}{2 R} \right) \left( f(x^k) - f^* \right) + \frac{\eta M}{2 c} \left( \frac{\sigma^2}{B} + \zeta^2 \right).
    \end{equation*}
     Then for $k = 0,1,2,..., N-1$ iterations that satisfy the conditions $\norms{\nabla f(x^k)} \geq c$ and $c \leq \sqrt{2} \zeta$, then ClipSGD has linear convergence
    \begin{equation*}
        \expect{f(x^{N})} - f^* \leq \left( 1 - \frac{\eta c}{2 R} \right)^{N} \left( f(x^0) - f^* \right) + \frac{M R}{c^2} \left( \frac{\sigma^2}{B} + \zeta^2 \right).
    \end{equation*}

\subsubsection{Second case: \texorpdfstring{$\frac{c}{3} \leq \norms{\nabla f(x^k)} \leq c$}{TEXT}}

In this case $\nabla f(x^k) = \clip{\nabla f(x^k)}$ with $\alpha = \min \left\{1, \frac{c}{\norms{\nabla f(x^k)}}\right\} = 1$, therefore we have the following
\begin{align}
    - \eta \dotprod{\nabla f(x^k)}{\clip{\gg(x^k, \bxi^k)}} &\overset{\eqref{eq:qudrat_raznosti}}{=} - \frac{\alpha \eta}{2} \norms{\nabla f(x^k)}^2 - \frac{\eta}{2 \alpha} \norms{\clip{\gg(x^k, \bxi^k)}}^2 \nonumber\\
    &\quad \quad \quad+ \frac{\eta}{2 \alpha} \norms{\clip{\gg(x^k, \bxi^k)} - \alpha \nabla f(x^k)}^2 \nonumber\\
    &= - \frac{\eta}{2} \norms{\nabla f(x^k)}^2 - \frac{\eta}{2} \norms{\clip{\gg(x^k, \bxi^k)}}^2 \nonumber\\
    &\quad \quad \quad+ \frac{\eta}{2} \norms{\clip{\gg(x^k, \bxi^k)} - \clip{\nabla f(x^k)}}^2 \nonumber\\
    &\leq - \frac{c \eta}{6} \norms{\nabla f(x^k)} - \frac{\eta}{2} \norms{\clip{\gg(x^k, \bxi^k)}}^2 \nonumber\\
    &\quad \quad \quad+ \frac{\eta}{2} \norms{\clip{\gg(x^k, \bxi^k)} - \clip{\nabla f(x^k)}}^2. \nonumber
\end{align}
Using that clipping is a projection on onto a convex set, namely ball with radius $c$, and thus is Lipshitz operator with Lipshitz constant $1$, we can obtain:
\begin{align*}
    - \eta \dotprod{\nabla f(x^k)}{\expect{\clip{\gg(x^k, \bxi^k)}}} &\leq - \frac{c \eta}{6} \norms{\nabla f(x^k)} - \frac{\eta}{2} \expect{\norms{\clip{\gg(x^k, \bxi^k)}}^2} \nonumber\\
    &\quad \quad \quad+ \frac{\eta}{2} \expect{\norms{\gg(x^k, \bxi^k) - \nabla f(x^k)}^2} \nonumber\\
    &\overset{\eqref{eq:Variance_decomposition}}{=}- \frac{c \eta}{6} \norms{\nabla f(x^k)} - \frac{\eta}{2} \expect{\norms{\clip{\gg(x^k, \bxi^k)}}^2} \nonumber\\
    &\quad \quad \quad+ \frac{\eta}{2} \expect{\norms{\gg(x^k, \bxi^k) - \expect{\gg(x^k, \bxi^k)}}^2} \nonumber\\
    &\quad \quad \quad+ \frac{\eta}{2} \norms{\bb(x^k)}^2\nonumber\\
    &\leq  - \frac{c \eta}{6} \norms{\nabla f(x^k)} - \frac{\eta}{2 } \expect{\norms{\clip{\gg(x^k, \bxi^k)}}^2} \nonumber\\
    &\quad \quad \quad+ \frac{\eta}{2} \left( \frac{\sigma^2}{B} + \zeta^2 \right) \nonumber\\
    &= - \frac{c \eta}{6} \norms{\nabla f(x^k)} - \frac{\eta}{2 } \expect{\norms{\clip{\gg(x^k, \bxi^k)}}^2} \nonumber\\
    &\quad \quad \quad+ \frac{\eta}{2} \left( \frac{\sigma^2}{B} + \zeta^2 \right). 
\end{align*}

Plugging this into \eqref{eq:Biased_clipSGD_smooth} and choosing $\eta \leq \frac{1}{4(L_0 + L_1 c)}$ we have:
    \begin{align}
        \expect{f(x^{k+1})} - f(x^k) &\overset{\eqref{eq:Biased_clipSGD_smooth}}{\leq} - \frac{c \eta}{6} \norms{\nabla f(x^k)} - \frac{\eta}{2} \expect{\norms{\clip{\gg(x^k, \bxi^k)}}^2}  \nonumber\\
        & \quad \quad \quad + \frac{\eta^2 (L_0 + L_1 \norms{\nabla f(x^k)})}{2} \expect{\norms{\clip{\gg(x^k, \bxi^k)}}^2} + \frac{\eta}{2} \left( \frac{\sigma^2}{B} + \zeta^2 \right) \nonumber\\
        &= - \frac{c \eta}{6} \norms{\nabla f(x^k)} - \frac{\eta}{2} \expect{\norms{\clip{\gg(x^k, \bxi^k)}}^2} \left(1 -  \eta (L_0 + L_1 \norms{\nabla f(x^k)}) \right) \nonumber\\
    &\quad \quad \quad+ \frac{\eta}{2} \left( \frac{\sigma^2}{B} + \zeta^2 \right) \nonumber\\
        &\leq - \frac{c \eta}{6} \norms{\nabla f(x^k)} + \frac{\eta}{2} \left( \frac{\sigma^2}{B} + \zeta^2 \right).\label{eq:Biased_clipSGD2_case1_1}
    \end{align}

    Using the convexity assumption of the function, we have the following:
    \begin{align*}
        f(x^k) - f^* &\leq \dotprod{\nabla f(x^k)}{x^k - x^*} \overset{\eqref{eq:scalar_product_bound}}{\leq} \norms{\nabla f(x^k)} \norms{x^k - x^*} \leq \norms{\nabla f(x^k)} \underbrace{\norms{x^0 - x^*}}_{R}.
    \end{align*}
    Hence we have:
    \begin{equation}
        \norms{\nabla f(x^k)} \geq \frac{f(x^k) - f^*}{R}. \label{eq:Biased_clipSGD2_case1_2}
    \end{equation}

    Then substituting \eqref{eq:Biased_clipSGD2_case1_2} into \eqref{eq:Biased_clipSGD2_case1_1} we obtain:
    \begin{equation*}
        \expect{f(x^{k+1})} - f(x^k) \leq - \frac{\eta c}{6} \norms{\nabla f(x^k)} + \frac{\eta}{2} \left( \frac{\sigma^2}{B} + \zeta^2 \right) \leq - \frac{\eta c}{6R} (f(x^k) - f^*) + \frac{\eta}{2} \left( \frac{\sigma^2}{B} + \zeta^2 \right).
    \end{equation*}
    
    This inequality is equivalent to the trailing inequality:
    \begin{equation*}
        \expect{f(x^{k+1})} - f^* \leq \left( 1 - \frac{\eta c}{6R} \right) \left( f(x^k) - f^* \right) + \frac{\eta}{2} \left( \frac{\sigma^2}{B} + \zeta^2 \right).
    \end{equation*}
     Then for $k = 0,1,2,..., N-1$ iterations that satisfy the conditions $\frac{c}{2} \leq \norms{\gg(x^k,\bxi^k)} \leq c $, then ClipSGD with biased gradient oracle has linear convergence
    \begin{equation*}
        \expect{f(x^{N})} - f^* \leq \left( 1 - \frac{\eta c}{6R} \right)^{N} \left( f(x^0) - f^* \right) + \frac{3R}{c} \left( \frac{\sigma^2}{B} + \zeta^2 \right).
    \end{equation*}

    Let $\mathcal{T}_1 = \left\{ m^{\mathcal{T}_1}_0, m^{\mathcal{T}_1}_1, m^{\mathcal{T}_1}_2, ..., m^{\mathcal{T}_1}_{K-1}  \right\} = \left\{k \in \{0,1,2,...,N-1\} |  \norms{\nabla f(x^k,\xi^k)} \geq \frac{c}{3} \right\}$, where $K = |\mathcal{T}_1|$. Then for $k \in \mathcal{T}_1$ ClipSGD with biased gradient oracle shows linear convergence:
        \begin{align*}
            F_N \cdot \mathds{1}\left[ \mathcal{T}_1 \right] &\lesssim \left(1 - \frac{\eta c}{R} \right) F_N \leq \left(1 - \frac{\eta c}{6R} \right) F_{m^{\mathcal{T}_1}_{K-1}} + \left(\frac{\eta M}{c} + \eta \right) \cdot\left( \frac{\sigma^2}{B} + \zeta^2 \right) \leq  ... \leq \\
            & \leq \left(1 - \frac{\eta c}{R} \right)^{K} F_{m_0^{\mathcal{T}_1}} + \left(\frac{ M R}{c^2} + \frac{R}{c} \right) \cdot\left( \frac{\sigma^2}{B} + \zeta^2 \right) \\
            &\leq \left(1 - \frac{\eta c}{R} \right)^{K} F_0 + \left(\frac{ M R}{c^2} + \frac{R}{c} \right) \cdot\left( \frac{\sigma^2}{B} + \zeta^2 \right),
        \end{align*}
        where $F_k =  \expect{f(x^k)} - f^*$, and we used that $F_{k} \leq F_{k-1}$.

\subsubsection{Third case: \texorpdfstring{$\norms{\nabla f(x^k)} \leq \frac{c}{3}$}{TEXT}}
    
    We introduce an indicative function:
    \begin{equation}
        \aleph_k = \mathds{1}\left\{ \norms{\gg(x^k, \bxi^k)} > c \right\}.\label{eq:Biased_indicator_function}
    \end{equation}
    Then the following is true:
    \begin{equation}
        \expect{\aleph_k} = \expect{\aleph_k^2} =\prob{\norms{\gg(x^k, \bxi^k)} > c} \overset{\circledOne}{\leq} \prob{\norms{\gg(x^k, \bxi^k) - \expect{\gg(x^k, \bxi^k)} } > \frac{c}{3}} \overset{\circledTwo}{\leq} \frac{9 \sigma^2}{c^2 B}, \label{eq:Biased_prob_indicator}
    \end{equation}
    where in $\circledOne$ we used $\norms{\gg(x^k, \bxi^k)} \leq \norms{\gg(x^k, \bxi^k) - \expect{\gg(x^k, \bxi^k)}} + \expect{\gg(x^k, \bxi^k)} \leq \norms{\gg(x^k, \bxi^k) - \expect{\gg(x^k, \bxi^k)}} + \frac{c}{2}$, where assume that $\zeta \leq \frac{c}{3}$: and in $\circledTwo$ we used Markov's inequality.
    
    Let $r_{k+1} = \expect{\norms{x^{k+1} - x^*}}$ and $F_{k+1} = \expect{f(x^{k+1}) - f^*}$, then given that
    \begin{align*}
        \clip{\gg(x^k, \bxi^k)} &= \gg(x^k, \bxi^k) (1 - \aleph_k) + \frac{c}{\norms{\gg(x^k, \bxi^k)}} \gg(x^k, \bxi^k) \aleph_k\\
        & = \gg(x^k, \bxi^k) + \left( \frac{c}{\norms{\gg(x^k, \bxi^k)}} - 1 \right) \gg(x^k, \bxi^k) \aleph_k
    \end{align*}
    we get with $\eta \leq \frac{1}{4(L_0 + L_1 c)}$:
    \begin{align}
        r_{k+1}^2 
        &= r_k^2 - 2 \eta \dotprod{\expect{\clip{\gg(x^k, \bxi^k)}}}{x^k - x^*} + \eta^2 \expect{\norms{\clip{\gg(x^k, \bxi^k)}}^2}\nonumber\\
        &= r_k^2 - 2 \eta \dotprod{\nabla f(x^k)}{x^k - x^*} - 2 \eta \dotprod{\expect{\left( \frac{c}{\norms{\gg(x^k, \bxi^k)}} - 1 \right) \gg(x^k, \bxi^k) \aleph_k}}{x^k - x^*} \nonumber\\
        &\quad\quad\quad + \eta^2 \expect{\norms{\clip{\gg(x^k, \bxi^k)}}^2} + 2 \eta \dotprod{\bb(x^k)}{x^k - x^*}\nonumber\\
        &\overset{\eqref{eq:scalar_product_bound}}{\leq} r_k^2 - 2 \eta \dotprod{\nabla f(x^k)}{x^k - x^*} + 2 \eta \norms{\expect{\left( \frac{c}{\norms{\gg(x^k, \bxi^k)}} - 1 \right) \gg(x^k, \bxi^k) \aleph_k}}\norms{x^k - x^*} \nonumber\\
        &\quad\quad\quad + \eta^2 \expect{\norms{\clip{\gg(x^k, \bxi^k)}}^2} + 2 \eta \norms{\bb(x^k)}\norms{x^k - x^*}\nonumber\\
        &\overset{\circledOne}{\leq} r_k^2 - 2 \eta F_k + 2 \eta \norms{\expect{\left( \frac{c}{\norms{\gg(x^k, \bxi^k)}} - 1 \right) \gg(x^k, \bxi^k) \aleph_k}}\norms{x^0 - x^*} \nonumber\\
        &\quad\quad\quad + \eta^2 \expect{\norms{\clip{\gg(x^k, \bxi^k)}}^2} + 2 \eta \norms{\bb(x^k)}\norms{x^0 - x^*}\nonumber\\
        &\overset{\eqref{eq:Squared_norm_of_the_sum}}{\leq}r_k^2 - 2 \eta F_k + 2 \eta \norms{\expect{\left( \frac{c}{\norms{\gg(x^k, \bxi^k)}} - 1 \right) \gg(x^k, \bxi^k) \aleph_k}}\norms{x^0 - x^*} \nonumber\\
        &\quad\quad\quad + 2 \eta^2 \expect{\norms{\clip{\gg(x^k, \bxi^k)} - \nabla f(x^k)}^2} + 2 \eta^2 \norms{\nabla f(x^k)}^2 + 2 \eta \zeta \norms{x^0 - x^*}\nonumber\\
        &=r_k^2 - 2 \eta F_k + 2 \eta \norms{\expect{\left( \frac{c}{\norms{\gg(x^k, \bxi^k)}} - 1 \right) \gg(x^k, \bxi^k) \aleph_k}} R \nonumber\\
        &\quad\quad\quad + 2 \eta^2 \expect{\norms{\clip{\gg(x^k, \bxi^k)} - \clip{\nabla f(x^k)}}^2} + 2 \eta^2 \norms{\nabla f(x^k)}^2 + 2 \eta \zeta R\nonumber\\
        &\overset{\circledTwo}{\leq}
        r_k^2 - 2 \eta F_k + 2 \eta \norms{\expect{\left( \frac{c}{\norms{\gg(x^k, \bxi^k)}} - 1 \right) \gg(x^k, \bxi^k) \aleph_k}}R \nonumber\\
        &\quad\quad\quad + 2 \eta^2 \expect{\norms{\gg(x^k, \bxi^k) - \nabla f(x^k)}^2} + 2 \eta^2 \norms{\nabla f(x^k)}^2 + 2 \eta \zeta R\nonumber\\
        &\overset{\eqref{eq:Variance_decomposition}}{=}r_k^2 - 2 \eta F_k + 2 \eta \norms{\expect{\left( \frac{c}{\norms{\gg(x^k, \bxi^k)}} - 1 \right) \gg(x^k, \bxi^k) \aleph_k}}R \nonumber\\
        &\quad\quad\quad + 2 \eta^2 \expect{\norms{\gg(x^k, \bxi^k) - \expect{\gg(x^k, \bxi^k)}}^2} + 2 \eta^2 \norms{\nabla f(x^k)}^2 + 2 \eta \zeta R + 2 \eta^2 \norms{\bb(x^k)}\nonumber\\
        &\leq r_k^2 - 2 \eta F_k + 2 \eta \norms{\expect{\left( \frac{c}{\norms{\gg(x^k, \bxi^k)}} - 1 \right) \gg(x^k, \bxi^k) \aleph_k}}R \nonumber\\
        &\quad\quad\quad +  \frac{ 2 \eta^2 \sigma^2}{B} + 2 \eta^2 \norms{\nabla f(x^k)}^2 + 2 \eta \zeta R + 2 \eta^2 \zeta^2\nonumber\\
        &\overset{\eqref{eq:lemma_smoothness}}{\leq} r_k^2 - 2 \eta F_k + 2 \eta \norms{\expect{\left( \frac{c}{\norms{\gg(x^k, \bxi^k)}} - 1 \right) \gg(x^k, \bxi^k) \aleph_k}}R \nonumber\\
        &\quad\quad\quad +  \frac{ 2 \eta^2 \sigma^2}{B} + 4 \eta^2 \left( L_0 + L_1 \norms{\nabla f(x^k)} \right) F_k + 2 \eta \zeta R + 2 \eta^2 \zeta^2\nonumber\\
        &\leq 
        r_k^2 - 2 \eta F_k + 2 \eta \norms{\expect{\left( \frac{c}{\norms{\gg(x^k, \bxi^k)}} - 1 \right) \gg(x^k, \bxi^k) \aleph_k}}R \nonumber\\
        &\quad\quad\quad + \frac{ 2 \eta^2 \sigma^2}{B} + 4 \eta^2 \left( L_0 + L_1 c \right) F_k + 2 \eta \zeta R + 2 \eta^2 \zeta^2\nonumber\\
        &= r_k^2 - 2 \eta F_k \left(1 - 2 \eta \left( L_0 + L_1 c \right)  \right) +  \frac{ 2 \eta^2 \sigma^2}{B} + 2 \eta \zeta R + 2 \eta^2 \zeta^2 \nonumber\\
        &\quad\quad\quad + 2 \eta \norms{\expect{\left( \frac{c}{\norms{\gg(x^k, \bxi^k)}} - 1 \right) \gg(x^k, \bxi^k) \aleph_k}}R\nonumber\\
        &\leq r_k^2 -  \eta F_k  +  \frac{ 2 \eta^2 \sigma^2}{B} + 2 \eta \zeta R + 2 \eta^2 \zeta^2 \nonumber\\
    &\quad \quad \quad+ 2 \eta \norms{\expect{\left( \frac{c}{\norms{\gg(x^k, \bxi^k)}} - 1 \right) \gg(x^k, \bxi^k) \aleph_k}}R.\label{eq:Biased_clipSGD3_case1_1}
    \end{align}

    Let's find the upper bound of the last summand:
    \begin{align}
        2 \eta R &\norms{\expect{\left( \frac{c}{\norms{\gg(x^k, \bxi^k)}} - 1 \right) \gg(x^k, \bxi^k) \aleph_k}}
        \nonumber\\
    & \overset{\eqref{eq:Biased_indicator_function}}{\leq}  2 \eta R \expect{\norms{\gg(x^k, \bxi^k)} \cdot \left( 1 - \frac{c}{ \norms{\gg(x^k, \bxi^k)}} \right) \aleph_k} \nonumber\\
        &\leq 2 \eta R  \expect{\norms{\gg(x^k, \bxi^k)} \cdot  \aleph_k}\nonumber\\
        & \leq 2 \eta R \left(\expect{\norms{\gg(x^k, \bxi^k) - \expect{\gg(x^k, \bxi^k)}} \cdot  \aleph_k} + \norms{\nabla f(x^k)} \expect{\aleph_k} + \norms{\bb(x^k)} \expect{\aleph_k} \right)\nonumber\\
        &\leq 2 \eta R \left(\sqrt{\expect{\norms{\gg(x^k, \bxi^k) - \expect{\gg(x^k, \bxi^k)}}^2} \cdot  \expect{\aleph_k^2}} + \frac{2c}{3} \expect{\aleph_k}\right)\nonumber\\
        &\overset{\eqref{eq:Biased_prob_indicator}}{\leq} 2 \eta R \left(\frac{3\sigma^2}{c B} + \frac{2c}{3} \cdot\frac{9\sigma^2}{c^2 B}\right) \nonumber\\
        &= \frac{18 \eta \sigma^2 R}{c B}.\label{eq:Biased_clipSGD3_case1_2}
    \end{align}

    Substituting into the initial formula and rearrange the summands, we obtain
    \begin{align*}
        \eta F_k &\overset{\eqref{eq:Biased_clipSGD3_case1_1}}{\leq} r_{k}^2 - r_{k+1}^2 +  \frac{ 2 \eta^2 \sigma^2}{B} + 2 \eta \zeta R + 2 \eta^2 \zeta^2 + 2 \eta \norms{\expect{\left( \frac{c}{\norms{\nabla f(x^k, \bxi^k)}} - 1 \right) \nabla f(x^k, \bxi^k) \aleph_k}}R\\
        &\overset{\eqref{eq:Biased_clipSGD3_case1_2}}{\leq} r_{k}^2 - r_{k+1}^2 + \frac{2 \eta^2 \sigma^2}{B} + \frac{18 \eta \sigma^2 R}{c B} + 2 \eta \zeta R + 2 \eta^2 \zeta^2
    \end{align*}

    Let $\mathcal{T}_2 = \left\{m^{\mathcal{T}_2}_0, m^{\mathcal{T}_2}_1, m^{\mathcal{T}_2}_2, ..., m^{\mathcal{T}_2}_{N-K}  \right\} = \left\{k \in \{0,1,2,...,N-1\} |  \norms{\nabla f(x^k)} < \frac{c}{3} \right\}$, where $|\mathcal{T}_2| = N-K$. Then rearranging and summing over all $k \in \mathcal{T}_2$ we obtain
        \begin{align*}
            F_N \cdot \mathds{1}\left[ \mathcal{T}_2 \right] &\leq \frac{1}{N - K} \sum_{k\in \mathcal{T}_2} F_{k} \leq \frac{1}{\eta (N - K)} \sum_{k\in \mathcal{T}_2} \left( r_{k}^2 - r_{k+1}^2 \right) \nonumber\\
    &\quad \quad \quad+ \frac{1}{N - K} \sum_{k\in \mathcal{T}_2}\left[ \eta \left( \frac{\sigma^2}{B} + \zeta^2 \right) + 2 R \left( \frac{9 \sigma^2}{c B} + \zeta \right)\right] \\
            &=\frac{r_0^2 - r_{N}^2}{\eta (N - K)} \eta \left( \frac{\sigma^2}{B} + \zeta^2 \right) + 2 R \left( \frac{9 \sigma^2}{c B} + \zeta \right) \leq \frac{r_0^2}{\eta(N - K)}  \eta \left( \frac{\sigma^2}{B} + \zeta^2 \right) + 2 R \left( \frac{9 \sigma^2}{c B} + \zeta \right). 
        \end{align*}
    Hence we obtain:
        \begin{equation*}
            F_N \cdot \mathds{1}\left[ \mathcal{T}_2 \right] \leq \frac{R^2}{\eta(N - K)}  \eta \left( \frac{\sigma^2}{B} + \zeta^2 \right) + 2 R \left( \frac{9 \sigma^2}{c B} + \zeta \right).
        \end{equation*}
    
    Combining all the cases considered, we obtain the convergence rate of ClipSGD with biased gradient oracle:
    \begin{align}
        \expect{f(x^{N})} - f^* &\leq F_N \cdot \mathds{1}\left[ \mathcal{T}_1 \right] + F_N \cdot \mathds{1}\left[ \mathcal{T}_2 \right] \nonumber\\
    &\lesssim  \left( 1 - \frac{\eta c}{R} \right)^{K} F_0 +  \frac{R^2}{\eta (N-K)} + \left(\frac{M R}{c^2} + \frac{R}{c} + \eta \right) \cdot\left( \frac{\sigma^2}{B} + \zeta^2 \right) + R \zeta. \label{eq:Biased_ClipSGD}
    \end{align}




\subsection{Convergence Results for ZO-ClipSGD}
In order to obtain convergence results for ZO-ClipSGD it is necessary to estimate the bias and variance of the gradient approximation~\eqref{eq:approximation_gradient}.

\paragraph{Bias of gradient approximation}
    Using the variational representation of the Euclidean norm, and definition of gradient approximation \eqref{eq:approximation_gradient} we can write:
    \begin{align}
        \norms{\expect{\gg(x, \{\xi,e\})} - \nabla f(x)} &= \norms{\expect{\frac{d}{2 \gamma}\left( \tilde{f}(x + \gamma e, \xi) - \tilde{f}(x - \gamma e, \xi) \right) e} - \nabla f(x)}
        \nonumber \\
        & \overset{\circledOne}{=} \norms{\expect{\frac{d}{\gamma}\left( f(x + \gamma e, \xi) + \delta(x + \gamma e) \right) e} - \nabla f(x)}
        \nonumber \\
        & \overset{\circledTwo}{\leq}  \norms{\expect{\frac{d}{\gamma} f(x + \gamma e, \xi) e} - \nabla f(x)} + \frac{d \Delta}{\gamma}
        \nonumber \\
        & \overset{\circledThree}{=}   \norms{\expect{\nabla f(x + \gamma u, \xi)} - \nabla f(x)} + \frac{d \Delta}{\gamma}
        \nonumber \\
        & = \sup_{z \in S^d(1)} \expect{\norms{ \nabla_z f(x + \gamma u, \xi) - \nabla_z f(x)}} + \frac{d \Delta}{\gamma}
        \nonumber \\
        & \overset{\eqref{eq:ass_smooth}}{\leq} \left(L_0 + L_1 \norms{\nabla f(x^k)}\right) \gamma \expect{\norms{u}} + \frac{d \Delta}{\gamma} 
        \nonumber \\
        & \leq \left(L_0 + L_1 M\right) \gamma + \frac{d \Delta}{\gamma}, \label{eq:proof_bias}
    \end{align}
    where $u \in B^d(1)$, $\circledOne =$ the equality is obtained from the fact, namely, distribution of $e$ is symmetric, $\circledTwo =$ the inequality is obtain from bounded noise $|\delta(x)| \leq \Delta$, $\circledThree =$ the equality is obtained from a version of Stokes’ theorem \citep[see Section 13.3.5, Exercise 14a, ][]{Zorich_2016}.

\paragraph{Bounding second moment (variance) of gradient approximation} By definition gradient approximation \eqref{eq:approximation_gradient} and Wirtinger-Poincare inequality \eqref{eq:Wirtinger_Poincare} we have
    \begin{align}
        &\expect{\norms{\gg(x,\{\xi,e\}) - \expect{\gg(x,\{\xi,e\})}}^2} \nonumber\\
        & \quad\quad \quad \leq \expect{\norms{\gg(x,\{\xi,e\})}^2}\nonumber\\
        & \quad\quad \quad= \frac{d^2}{4 \gamma^2} \expect{\norms{\left(\tilde{f}(x + \gamma e, \xi) - \tilde{f}(x - \gamma e, \xi)\right) e}^2}
        \nonumber \\
        &\quad\quad \quad = \frac{d^2}{4 \gamma^2} \expect{\left(f(x + \gamma e, \xi) - f(x - \gamma e, \xi) + \delta (x + \gamma e) - \delta (x -\gamma e)\right)^2}
        \nonumber \\
        &\quad\quad \quad \overset{\eqref{eq:Squared_norm_of_the_sum}}{\leq} \frac{d^2}{2 \gamma^2} \left( \expect{\left(f(x + \gamma e, \xi) - f(x - \gamma e, \xi)\right)^2} + 2 \Delta^2 \right)
        \nonumber \\
        &\quad\quad \quad \overset{\eqref{eq:Wirtinger_Poincare}}{\leq} \frac{d^2}{2 \gamma^2} \left( \frac{\gamma^2}{d} \expect{\norms{ \nabla f(x + \gamma e, \xi) + \nabla f(x - \gamma e, \xi)}^2} + 2 \Delta^2 \right) 
        \nonumber \\
        &\quad \quad \quad = \frac{d^2}{2 \gamma^2} \left( \frac{\gamma^2}{d} \expect{\norms{ \nabla f(x + \gamma e, \xi) + \nabla f(x - \gamma e, \xi) \pm 2 \nabla f(x, \xi)}^2} + 2 \Delta^2 \right)
        \nonumber \\
        &\quad\quad \quad \overset{\eqref{eq:ass_smooth}}{\leq} 4d \expect{\norms{\nabla f(x, \xi)}^2} + 4 d  L^2 \gamma^2 \expect{\norms{e}^2}  + \frac{d^2 \Delta^2}{\gamma^2}  
        \nonumber \\
        &\quad\quad \quad \overset{\circledOne}{\leq} 4d \tilde{\sigma}^2 + 4 d  \left(L_0 + L_1 \norms{\nabla f(x^k)}\right)^2 \gamma^2 \expect{\norms{e}^2}  + \frac{d^2 \Delta^2}{\gamma^2}\nonumber \\
        &\quad\quad \quad \leq 4d \tilde{\sigma}^2 + 4 d  \left(L_0 + L_1 M\right)^2 \gamma^2  + \frac{d^2 \Delta^2}{\gamma^2}, \label{eq:proof_variance}
    \end{align}
    where $\circledOne =$ the inequality is obtain from $\expect{\norms{\nabla f(x,\xi)}^2} \leq \tilde{\sigma}^2$.

\subsubsection{Proof of Theorem~\ref{th:ZO_ClipSGD}}
In order to obtain the convergence rate of ZO-ClipSGD in the convex setting, we need to substitute the obtained estimates~\eqref{eq:proof_bias} and \eqref{eq:proof_variance} into the convergence rate of ClipSGD~\eqref{eq:Biased_ClipSGD} instead of $\zeta$ and $\sigma^2$, respectively. Given that $\frac{M R}{c^2} + \frac{R}{c} + \eta \lesssim \frac{M R}{c^2}$ at small $c$, then the convergence of ZO-ClipSGD in the convex setup is as follows:
    \begin{align*}
        \expect{f(x^{N})} - f^* &\lesssim  \underbrace{\left( 1 - \frac{\eta}{R} \right)^{K} (f(x^0) - f^*)}_{\circledOne} + \underbrace{\frac{R^2}{\eta(N-K)}}_{\circledTwo} + \underbrace{\frac{d M R \tilde{\sigma}^2}{c^2 B}}_{\circledThree} + \underbrace{\frac{d M R \left(L_0 + L_1 M\right)^2 \gamma^2}{c^2 B}}_{\circledFour} \nonumber \\
        &\quad \quad \quad+ \underbrace{\frac{d^2 M R \Delta^2}{c^2 B \gamma^2}}_{\circledFive} + \underbrace{\frac{M R \left(L_0 + L_1 M\right)^2 \gamma^2}{c^2}}_{\circledSix} + \underbrace{\frac{d^2 M R \Delta^2}{c^2 \gamma^2}}_{\circledSeven} \nonumber \\
        &\quad \quad \quad+ \underbrace{ \left(L_0 + L_1 M\right) \gamma R}_{\circledEight} + \underbrace{ \frac{d \Delta R}{\gamma}}_{\circledNine}.
    \end{align*}

    \textbf{From term $\circledOne$}, we find the $K$:
    \begin{align}
        \circledOne: \quad \left( 1 - \frac{\eta c}{R} \right)^{K} (f(x^0) - f^*) \leq \varepsilon \quad & \Rightarrow \quad K \geq \frac{R}{\eta c} \log \frac{f(x^0) - f^*}{\varepsilon}. \label{eq:ZO_ClipSGD_KKKKK}
    \end{align}

    \textbf{From term $\circledTwo$}, we find the number of iterations $N$ required for Algorithm \ref{algo:ZO_ClipSGD} in convex setup to achieve $\varepsilon$-accuracy:
    \begin{align}
        \circledTwo: \quad \frac{R^2}{\eta(N-K)} \leq \varepsilon \quad & \Rightarrow \quad N \overset{\eqref{eq:ZO_ClipSGD_KKKKK}}{\geq} \frac{R^2}{\eta \varepsilon} +  \frac{R}{\eta c} \log \frac{f(x^0) - f^*}{\varepsilon};
        \nonumber \\
         N &= \Obound{\frac{R^2}{\eta \varepsilon} +  \frac{R}{\eta c} \log \frac{1}{\varepsilon}}. \label{eq:ZO_ClipSGD_iterations}
    \end{align}

    \textbf{From terms $\circledThree$}, we find the batch size $B$: 
    \begin{align}
        \circledThree: \quad \frac{ d MR \tilde{\sigma}^2}{c^2 B} \leq \varepsilon \quad & \Rightarrow \quad B \geq \frac{d M R \tilde{\sigma}^2}{\varepsilon c^2};
        \nonumber \\
         B &= \Obound{ \frac{d M R \tilde{\sigma}^2}{\varepsilon c^2}}. \label{eq:ZO_ClipSGD_batch_size}
    \end{align}

    \textbf{From terms $\circledFour$, $\circledSix$ and $\circledEight$} we find the smoothing parameter $\gamma$:
    \begin{align}
        &\circledFour: \quad \frac{d MR \left(L_0 + L_1 M\right)^2 \gamma^2}{c^2 B} \leq 
        \varepsilon \quad \Rightarrow \quad \gamma \leq \sqrt{\frac{\varepsilon c^2 B}{ d MR \left(L_0 + L_1 M\right)^2}} \overset{\eqref{eq:ZO_ClipSGD_batch_size}}{=} \frac{\tilde{\sigma}}{\left(L_0 + L_1 M\right)};
        \nonumber \\
        & \circledSix: \quad  \frac{MR \left(L_0 + L_1 M\right)^2 \gamma^2}{c^2} \leq \varepsilon \quad \Rightarrow \quad \gamma \leq \frac{\sqrt{\varepsilon} c}{\sqrt{MR} \left(L_0 + L_1 M\right)};
        \nonumber \\
        & \circledEight: \quad  \left(L_0 + L_1 M\right) R \gamma \leq \varepsilon \quad \Rightarrow \quad \gamma \leq \frac{\varepsilon}{R \left(L_0 + L_1 M\right)};
        \nonumber \\
        & \quad \quad \quad \gamma \leq \frac{1}{\left(L_0 + L_1 M\right)}\min \left\{\tilde{\sigma}, \frac{\sqrt{\varepsilon} c}{\sqrt{MR}}, \frac{\varepsilon}{R} \right\} = \frac{\varepsilon}{R \left(L_0 + L_1 M\right)}. \label{eq:ZO_ClipSGD_smoothing_parameter}
    \end{align}

    \textbf{From the remaining terms $\circledFive$, $\circledSeven$ and $\circledNine$}, we find the maximum allowable level of adversarial noise $\Delta$ that still guarantees the convergence of the ZO-ClipSGD to desired accuracy~$\varepsilon$ in convex setup:
    \begin{align}
        &\circledFive: \quad \frac{d^2 MR \Delta^2}{c^2 B \gamma^2} \leq \varepsilon \quad \Rightarrow \quad \Delta \leq \frac{\sqrt{\varepsilon} c \gamma \sqrt{B}}{d \sqrt{MR}} \overset{ \eqref{eq:ZO_ClipSGD_batch_size}, \eqref{eq:ZO_ClipSGD_smoothing_parameter}}{=}  \frac{\varepsilon \tilde{\sigma}}{\sqrt{d} \left(L_0 + L_1 M\right) R};
        \nonumber \\
        & \circledSeven: \quad \frac{d^2 MR \Delta^2}{\gamma^2 c^2} \leq \varepsilon \quad \Rightarrow \quad \Delta \leq \sqrt{\frac{ \gamma^2 c^2 \varepsilon }{d^2 MR}} \overset{\eqref{eq:ZO_ClipSGD_smoothing_parameter}}{=} \frac{\varepsilon^{3/2} c}{d \left(L_0 + L_1 M\right) \sqrt{M}R^{3/2}};
        \nonumber \\
        & \circledNine: \quad \frac{d  \Delta R}{\gamma} \leq \varepsilon \quad \Rightarrow \quad \Delta \leq \sqrt{\frac{ \gamma \varepsilon }{d R}} \overset{\eqref{eq:ZO_ClipSGD_smoothing_parameter}}{=} \frac{\varepsilon^2}{d \left(L_0 + L_1 M\right) R^2};
        \nonumber \\
        & \quad \quad \quad \Delta \leq  \frac{\varepsilon}{\sqrt{d} \left(L_0 + L_1 M\right) R} \min \left\{ \tilde{\sigma},\frac{\sqrt{\varepsilon} c}{\sqrt{d} \sqrt{MR}}, \frac{\varepsilon}{\sqrt{d}R}\right\} \nonumber \\
        &\quad \quad \quad\quad= \frac{\varepsilon}{\sqrt{d} \left(L_0 + L_1 M\right) R} \min \left\{ \tilde{\sigma}, \frac{\varepsilon}{\sqrt{d}R}\right\}. 
        \label{eq:ZO_ClipSGD_noise_level}
    \end{align}

    In this way, the ZO-ClipSGD achieves $\varepsilon$-accuracy: $\expect{f(x^{N}) - f^*} \leq \varepsilon$ in convex setup after 
        \begin{equation*}
            N \overset{\eqref{eq:ZO_ClipSGD_iterations}}{=} \Obound{\frac{R^2}{\eta \varepsilon} +  \frac{R}{\eta c} \log \frac{1}{\varepsilon}}, \quad T = N \cdot B \overset{\eqref{eq:ZO_ClipSGD_iterations}, \eqref{eq:ZO_ClipSGD_batch_size}}{=} \Obound{\frac{d \tilde{\sigma}^2 M R^2}{\varepsilon c^2 \eta} \left( \frac{1}{c}\log \frac{1}{\varepsilon} + \frac{R}{\varepsilon} \right)}
        \end{equation*}
        number of iterations, total number of zero-order oracle calls and at
        \begin{equation*}
            \Delta \overset{\eqref{eq:ZO_ClipSGD_noise_level}}{\lesssim} \frac{\varepsilon}{\sqrt{d} \left(L_0 + L_1 M\right) R} \min \left\{ \tilde{\sigma}, \frac{\varepsilon}{\sqrt{d}R}\right\}
        \end{equation*}
        the maximum level of noise with smoothing parameter $\frac{\varepsilon}{\left(L_0 + L_1 M\right) R}$ \eqref{eq:ZO_ClipSGD_smoothing_parameter}.

\section{Zero-Order Normalized Stochastic Gradient Descent Method}
\label{app:ZO_NSGD}

This section consists of two parts: 1) a generalization of the convergence result of NSGD (Algorithm~\ref{algo:NSGD}) to the biased gradient oracle $\gg(x^k, \bxi^k) = \nabla f(x^k, \bxi^k) + \bb(x^k)$, where $\bb(x^k)$ is biased bounded by $\zeta \geq 0: \norms{\bb(x^k)}\leq \zeta$; 2) deriving convergence estimates of ZO-NSGD directly. 

\subsection{Biased Normalized Stochastic Gradient Descent Method (Proof of the Lemma~\ref{lem:biased_ZO_NSGD})}

Let's introduce the notation $G(x^k, \bxi^k) = \frac{\gg(x^k, \bxi^k)}{\norms{\gg(x^k, \bxi^k)}}$, then using $(L_0,L_1)$-smoothness (see Assumption~\ref{ass:L0_L1_smooth}):
\begin{align}
    f(x^{k+1}) - f(x^k) &\overset{\eqref{eq:ass_smooth}}{\leq} \dotprod{\nabla f(x^k)}{x^{k+1} - x^k} + \frac{L_0 + L_1 \norms{\nabla f(x^k)}}{2} \norms{x^{k+1} - x^{k}}^2 \nonumber \\
    &= - \eta \dotprod{\nabla f(x^k)}{G(x^k, \bxi^k)} + \frac{\eta^2 (L_0 + L_1 \norms{\nabla f(x^k)})}{2} \norms{G(x^k, \bxi^k)}^2. \label{eq:Biased_NSGD_smooth}
\end{align}

Next, we consider 4 cases of the relation $\norms{\nabla f(x^k)}$ and $\norms{\gg(x^k, \bxi^k)}$ with respect to the hyperparameter $\lambda$.

\subsubsection{First case: \texorpdfstring{$\norms{\nabla f(x^k)} \geq \lambda$}{TEXT} and \texorpdfstring{$\norms{\gg(x^k, \bxi^k)} \geq \lambda$}{TEXT}}

    Let us evaluate first summand of \eqref{eq:Biased_NSGD_smooth} with $\alpha = \norms{\nabla f(x^k)}^{-1}$:
    \begin{align}
        - \eta \dotprod{\nabla f(x^k)}{G(x^k, \bxi^k)} &\overset{\eqref{eq:qudrat_raznosti}}{=} - \frac{\alpha \eta }{2} \norms{\nabla f(x^k)}^2 - \frac{\eta}{2 \alpha} \norms{G(x^k, \bxi^k)}^2 \nonumber\\
        & \quad \quad \quad + \frac{\eta}{2 \alpha} \norms{G(x^k, \bxi^k) - \alpha \nabla f(x^k)}^2 \nonumber\\
        &= - \frac{\eta}{2 } \norms{\nabla f(x^k)} - \frac{\eta}{2 \alpha} \norms{G(x^k, \bxi^k)}^2 \nonumber\\
        & \quad \quad \quad+ \frac{\eta}{2 \lambda^2 \alpha} \norms{\lambda G(x^k, \bxi^k) - \lambda \alpha \nabla f(x^k)}^2 \nonumber\\
        &= - \frac{\eta}{2 } \norms{\nabla f(x^k)} - \frac{\eta}{2 \alpha} \norms{G(x^k, \bxi^k)}^2 \nonumber\\
        & \quad \quad \quad+ \frac{\eta}{2 \lambda^2 \alpha} \norms{\cliplam{\gg(x^k, \bxi^k)} - \cliplam{ \nabla f(x^k)}}^2 \nonumber
    \end{align}
    Using that clipping is a projection on onto a convex set, namely ball with radius $\lambda$, and thus is Lipshitz operator with Lipshitz constant $1$, we can obtain:
    \begin{align}
        - \eta \dotprod{\nabla f(x^k)}{\expect{G(x^k, \bxi^k)}} &\leq - \frac{\eta}{2} \norms{\nabla f(x^k)} - \frac{\eta}{2 \alpha} \expect{\norms{G(x^k, \bxi^k)}^2} \nonumber\\
        & \quad \quad \quad+ \frac{\eta }{2 \lambda^2 \alpha} \expect{\norms{\gg(x^k, \bxi^k) - \nabla f(x^k)}^2}.\label{eq:Biased_NSGD1_main}
    \end{align}

    \paragraph{In the case: $0 \leq \zeta \leq \frac{\lambda}{\sqrt{2}}$.} Using this in \eqref{eq:Biased_NSGD1_main}, we have the following with $\eta_k \leq \frac{\norms{\nabla f(x^k)}}{2(L_0 + L_1 \norms{\nabla f(x^k)})}$:
        \begin{align}
            \expect{f(x^{k+1})} - f(x^k) &\overset{\eqref{eq:Biased_NSGD_smooth}}{\leq}  - \eta \dotprod{\nabla f(x^k)}{\expect{G(x^k, \bxi^k)}} + \frac{\eta^2 (L_0 + L_1 \norms{\nabla f(x^k)})}{2} \expect{\norms{G(x^k, \bxi^k)}^2} \nonumber\\
            & \overset{\eqref{eq:Biased_NSGD1_main}}{\leq} - \frac{\eta}{2} \norms{\nabla f(x^k)} - \frac{\eta}{2 \alpha} \expect{\norms{G(x^k, \bxi^k)}^2} + \frac{\eta }{2 \lambda^2 \alpha} \expect{\norms{\gg(x^k, \bxi^k) - \nabla f(x^k)}^2}\nonumber\\
            &\quad\quad\quad+ \frac{\eta^2 (L_0 + L_1 \norms{\nabla f(x^k)})}{2} \expect{\norms{G(x^k, \bxi^k)}^2} \nonumber\\
            &= - \frac{\eta}{2} \norms{\nabla f(x^k)} + \frac{\eta }{2 \lambda^2 \alpha} \expect{\norms{\gg(x^k, \bxi^k) - \nabla f(x^k)}^2}\nonumber\\
            &\quad\quad\quad - \frac{\eta}{2} \expect{\norms{G(x^k, \bxi^k)}^2} \left( 1 - \frac{\eta (L_0 + L_1 \norms{\nabla f(x^k)})}{\norms{\nabla f(x^k)}} \right)\nonumber\\
            &\leq - \frac{\eta}{2} \norms{\nabla f(x^k)} + \frac{\eta }{2 \lambda^2 \alpha} \expect{\norms{\gg(x^k, \bxi^k) - \nabla f(x^k)}^2}\nonumber\\
            &\overset{\eqref{eq:Variance_decomposition}}{=}- \frac{\eta}{2} \norms{\nabla f(x^k)} + \frac{\eta }{2 \lambda^2 \alpha} \expect{\norms{\gg(x^k, \bxi^k) - \expect{\gg(x^k, \bxi^k)}}}  + \frac{\eta }{2 \lambda^2 \alpha} \norms{\bb(x^k)}^2\nonumber\\
            &\leq - \frac{\eta}{2} \norms{\nabla f(x^k)} + \frac{\eta \sigma^2}{2 \lambda^2 \alpha B} + \frac{\eta \zeta^2}{2 \lambda^2 \alpha} \nonumber\\
            &\leq - \frac{\eta}{2} \norms{\nabla f(x^k)} + \frac{\eta}{4} \norms{\nabla f(x^k)} + \frac{\eta \sigma^2 M}{2 \lambda^2 B}\nonumber\\
            &= - \frac{\eta}{4} \norms{\nabla f(x^k)} + \frac{\eta \sigma^2 M}{2 \lambda^2 B}. \label{eq:Biased_NSGD1_case1_1}
        \end{align}
        
        The step size will be constant, depending on the hyperparameter $\lambda$:
        \begin{align*}
            \frac{\norms{\nabla f(x^k)}}{2\left(L_0 + L_1 \norms{\nabla f(x^k)}\right)} = \frac{1}{2\left(L_0\frac{1}{\norms{\nabla f(x^k)}} + L_1\right)} = \frac{\lambda}{2\left(L_0\frac{\lambda}{\norms{\nabla f(x^k)}} + L_1 \lambda\right)} \geq \frac{\lambda}{2\left(L_0 + L_1 \lambda\right)}.
        \end{align*}
        Thus, $\eta_k = \eta \leq \frac{\lambda}{2(L_0 + L_1 \lambda)}$.

        Using the convexity assumption of the function, we have the following:
        
            \begin{align*}
                f(x^k) - f^* &\leq \dotprod{\nabla f(x^k)}{x^k - x^*} \overset{\eqref{eq:scalar_product_bound}}{\leq} \norms{\nabla f(x^k)} \norms{x^k - x^*} \leq \norms{\nabla f(x^k)} \underbrace{\norms{x^0 - x^*}}_{R}.
            \end{align*}
            Hence we have:
            \begin{equation}
                \norms{\nabla f(x^k)} \geq \frac{f(x^k) - f^*}{R}. \label{eq:Biased_NSGD1_case1_2}
            \end{equation}
        
            Then substituting \eqref{eq:Biased_NSGD1_case1_2} into \eqref{eq:Biased_NSGD1_case1_1} we obtain:
            \begin{equation*}
                \expect{f(x^{k+1})} - f(x^k) \leq - \frac{\eta}{4} \norms{\nabla f(x^k)} + \frac{\eta \sigma^2 M}{2 \lambda^2 B} \leq - \frac{\eta}{4R} (f(x^k) - f^*) + \frac{\eta \sigma^2 M}{2 \lambda^2 B}.
            \end{equation*}
            
            This inequality is equivalent to the trailing inequality:
            \begin{equation*}
                \expect{f(x^{k+1})} - f^* \leq \left( 1 - \frac{\eta}{4R} \right) \left( f(x^k) - f^* \right) + \frac{\eta \sigma^2 M}{2 \lambda^2 B}.
            \end{equation*}
            
                Then for $k = 0,1,2,..., N-1$ iterations that satisfy the conditions $\norms{\gg(x^k,\bxi^k)} \geq \sqrt{2} \zeta$ and $\norms{\nabla f(x^k)} \geq \sqrt{2} \zeta$ NSGD with biased gradient oracle shows linear convergence:
                \begin{equation*}
                    \expect{f(x^{N})} - f^* \leq \left(1 - \frac{\eta}{4R} \right)^N (f(x^{0}) - f^*) + \frac{2 \sigma^2 M R}{ \lambda^2 B}.
                \end{equation*}

    \paragraph{In the case: $\frac{\lambda}{\sqrt{2}} \leq \zeta $.} Using this in \eqref{eq:Biased_NSGD1_main}, we have the following with $\eta_k \leq \frac{\norms{\nabla f(x^k)}}{2(L_0 + L_1 \norms{\nabla f(x^k)})}$:
        \begin{align}
            \expect{f(x^{k+1})} - f(x^k) &\overset{\eqref{eq:Biased_NSGD_smooth}}{\leq}  - \eta \dotprod{\nabla f(x^k)}{\expect{G(x^k, \bxi^k)}} + \frac{\eta^2 (L_0 + L_1 \norms{\nabla f(x^k)})}{2} \expect{\norms{G(x^k, \bxi^k)}^2} \nonumber\\
            & \overset{\eqref{eq:Biased_NSGD1_main}}{\leq} - \frac{\eta}{2} \norms{\nabla f(x^k)} - \frac{\eta}{2 \alpha} \expect{\norms{G(x^k, \bxi^k)}^2} + \frac{\eta }{2 \lambda^2 \alpha} \expect{\norms{\gg(x^k, \bxi^k) - \nabla f(x^k)}^2}\nonumber\\
            &\quad\quad\quad+ \frac{\eta^2 (L_0 + L_1 \norms{\nabla f(x^k)})}{2} \expect{\norms{G(x^k, \bxi^k)}^2} \nonumber\\
            &= - \frac{\eta}{2} \norms{\nabla f(x^k)} + \frac{\eta }{2 \lambda^2 \alpha} \expect{\norms{\gg(x^k, \bxi^k) - \nabla f(x^k)}^2}\nonumber\\
            &\quad\quad\quad - \frac{\eta}{2} \expect{\norms{G(x^k, \bxi^k)}^2} \left( 1 - \frac{\eta (L_0 + L_1 \norms{\nabla f(x^k)})}{\norms{\nabla f(x^k)}} \right)\nonumber\\
            &\leq - \frac{\eta}{2} \norms{\nabla f(x^k)} + \frac{\eta }{2 \lambda^2 \alpha} \expect{\norms{\gg(x^k, \bxi^k) - \nabla f(x^k)}^2}\nonumber\\
            &\overset{\eqref{eq:Variance_decomposition}}{=}- \frac{\eta}{2} \norms{\nabla f(x^k)} + \frac{\eta }{2 \lambda^2 \alpha} \expect{\norms{\gg(x^k, \bxi^k) - \expect{\gg(x^k, \bxi^k)}}}  + \frac{\eta }{2 \lambda^2 \alpha} \norms{\bb(x^k)}^2\nonumber\\
            &\leq - \frac{\eta}{2} \norms{\nabla f(x^k)} + \frac{\eta \sigma^2}{2 \lambda^2 \alpha B} + \frac{\eta \zeta^2}{2 \lambda^2 \alpha} \nonumber\\
            &\leq - \frac{\eta}{2} \norms{\nabla f(x^k)} + \frac{\eta \sigma^2 M}{2 \lambda^2 B} + \frac{\eta \zeta^2 M}{2 \lambda^2} \nonumber\\
            &= - \frac{\eta}{2} \norms{\nabla f(x^k)} + \frac{\eta \sigma^2 M}{2 \lambda^2 B} + \frac{\eta \zeta^2 M}{2 \lambda^2}.  \label{eq:Biased_NSGD1_case2_1}
        \end{align}
        
        The step size will be constant, depending on the hyperparameter $\lambda$:
        \begin{align*}
            \frac{\norms{\nabla f(x^k)}}{2\left(L_0 + L_1 \norms{\nabla f(x^k)}\right)} = \frac{1}{2\left(L_0\frac{1}{\norms{\nabla f(x^k)}} + L_1\right)} = \frac{\lambda}{2\left(L_0\frac{\lambda}{\norms{\nabla f(x^k)}} + L_1 \lambda\right)} \geq \frac{\lambda}{2\left(L_0 + L_1 \lambda\right)}.
        \end{align*}
        Thus, $\eta_k = \eta \leq \frac{\lambda}{2(L_0 + L_1 \lambda)}$.
        
        Using the convexity assumption of the function, we have the following:

        \begin{align*}
            f(x^k) - f^* &\leq \dotprod{\nabla f(x^k)}{x^k - x^*} \overset{\eqref{eq:scalar_product_bound}}{\leq} \norms{\nabla f(x^k)} \norms{x^k - x^*} \leq \norms{\nabla f(x^k)} \underbrace{\norms{x^0 - x^*}}_{R}.
        \end{align*}
        Hence we have:
        \begin{equation}
            \norms{\nabla f(x^k)} \geq \frac{f(x^k) - f^*}{R}. \label{eq:Biased_NSGD1_case2_2}
        \end{equation}
    
        Then substituting \eqref{eq:Biased_NSGD1_case2_2} into \eqref{eq:Biased_NSGD1_case2_1} we obtain:
        \begin{equation*}
            \expect{f(x^{k+1})} - f(x^k) \leq - \frac{\eta}{2} \norms{\nabla f(x^k)} + \frac{\eta \sigma^2 M}{2 \lambda^2 B} + \frac{\eta \zeta^2 M}{2 \lambda^2} \leq - \frac{\eta}{2 R} (f(x^k) - f^*) + \frac{\eta \sigma^2 M}{2 \lambda^2B} + \frac{\eta \zeta^2 M}{2 \lambda^2}.
        \end{equation*}
        
        This inequality is equivalent to the trailing inequality:
        \begin{equation*}
            \expect{f(x^{k+1})} - f^* \leq \left( 1 - \frac{\eta}{2R} \right) \left( f(x^k) - f^* \right) + \frac{\eta \sigma^2 M}{2 \lambda^2 B} + \frac{\eta \zeta^2 M}{2 \lambda^2}.
        \end{equation*}
    
        Then for $k = 0,1,2,..., N-1$ iterations that satisfy the conditions $\norms{\gg(x^k,\bxi^k)} \geq \lambda$ and $\norms{\nabla f(x^k)} \geq \lambda $ and $\zeta \geq \sqrt{2} \lambda$  NSGD with biased gradient oracle shows linear convergence:
        \begin{equation*}
            \expect{f(x^{N})} - f^* \leq \left(1 - \frac{\eta}{2R} \right)^N (f(x^{0}) - f^*) + \frac{\sigma^2 M R}{\lambda^2 B} + \frac{\zeta^2 M R}{\lambda^2}.
        \end{equation*}

\subsubsection{Second case: \texorpdfstring{$\norms{\nabla f(x^k)} \leq \lambda$}{TEXT} and \texorpdfstring{$\norms{\gg(x^k, \bxi^k)} \geq \lambda$}{TEXT}}

    Let us evaluate first summand of \eqref{eq:Biased_NSGD_smooth} with $\alpha = \lambda^{-1}$:
    \begin{align}
        - \eta \dotprod{\nabla f(x^k)}{G(x^k, \bxi^k)} &\overset{\eqref{eq:qudrat_raznosti}}{=} - \frac{\alpha \eta }{2} \norms{\nabla f(x^k)}^2 - \frac{\eta}{2 \alpha} \norms{G(x^k, \bxi^k)}^2 \nonumber\\
        & \quad \quad \quad+ \frac{\eta}{2 \alpha} \norms{G(x^k, \bxi^k) - \alpha \nabla f(x^k)}^2 \nonumber\\
        &\leq - \frac{\eta}{2} \norms{\nabla f(x^k)} - \frac{\eta}{2 \alpha} \norms{G(x^k, \bxi^k)}^2 \nonumber\\
        & \quad \quad \quad+ \frac{\eta}{2 \lambda} \norms{\lambda G(x^k, \bxi^k) - \nabla f(x^k)}^2 \nonumber\\
        &= - \frac{\eta}{2} \norms{\nabla f(x^k)} - \frac{\eta}{2 \alpha} \norms{G(x^k, \bxi^k)}^2 \nonumber\\
        & \quad \quad \quad+ \frac{\eta}{2 \lambda} \norms{\cliplam{\gg(x^k, \bxi^k)} - \cliplam{ \nabla f(x^k)}}^2 \nonumber
    \end{align}
    Using that clipping is a projection on onto a convex set, namely ball with radius $\lambda$, and thus is Lipshitz operator with Lipshitz constant $1$, we can obtain:
    \begin{align}
        - \eta \dotprod{\nabla f(x^k)}{\expect{G(x^k, \bxi^k)}} &\leq - \frac{\eta}{2} \norms{\nabla f(x^k)} - \frac{\eta}{2 \alpha} \expect{\norms{G(x^k, \bxi^k)}^2} \nonumber\\
        & \quad \quad \quad+ \frac{\eta }{2 \lambda} \expect{\norms{\gg(x^k, \bxi^k) - \nabla f(x^k)}^2}.
        \label{eq:Biased_NSGD2_main}
    \end{align}
     Using this, we have the following with $\eta_k \leq \frac{\norms{\nabla f(x^k)}}{2(L_0 + L_1 \norms{\nabla f(x^k)})}$:
    \begin{align}
        \expect{f(x^{k+1})} - f(x^k) &\overset{\eqref{eq:Biased_NSGD_smooth}}{\leq}  - \eta \dotprod{\nabla f(x^k)}{\expect{G(x^k, \bxi^k)}} + \frac{\eta^2 (L_0 + L_1 \norms{\nabla f(x^k)})}{2} \expect{\norms{G(x^k, \bxi^k)}^2} \nonumber\\
        & \overset{\eqref{eq:Biased_NSGD2_main}}{\leq} - \frac{\eta}{2} \norms{\nabla f(x^k)} - \frac{\eta}{2 \alpha} \expect{\norms{G(x^k, \bxi^k)}^2} + \frac{\eta }{2 \lambda} \expect{\norms{\gg(x^k, \bxi^k) - \nabla f(x^k)}^2}\nonumber\\
        &\quad\quad\quad+ \frac{\eta^2 (L_0 + L_1 \norms{\nabla f(x^k)})}{2} \expect{\norms{G(x^k, \bxi^k)}^2} \nonumber\\
        &\overset{\eqref{eq:Variance_decomposition}}{=} - \frac{\eta}{2} \norms{\nabla f(x^k)} + \frac{\eta }{2 \lambda} \expect{\norms{\gg(x^k, \bxi^k) - \expect{\gg(x^k, \bxi^k)}}^2} + \frac{\eta }{2 \lambda} \norms{\bb(x^k)}^2\nonumber\\
        &\quad\quad\quad - \frac{\eta}{2} \expect{\norms{G(x^k, \bxi^k)}^2} \left( 1 - \frac{\eta (L_0 + L_1 \norms{\nabla f(x^k)})}{\norms{\nabla f(x^k)}} \right)\nonumber\\
        &\leq - \frac{\eta}{2} \norms{\nabla f(x^k)} + \frac{\eta \sigma^2}{2 \lambda B} + \frac{\eta \zeta^2}{2 \lambda}. \label{eq:Biased_NSGD2_case1_1}
    \end{align}
        
    The step size will be constant, depending on the hyperparameter $\lambda$:
    \begin{align*}
        \frac{\norms{\nabla f(x^k)}}{2\left(L_0 + L_1 \norms{\nabla f(x^k)}\right)} = \frac{1}{2\left(L_0\frac{1}{\norms{\nabla f(x^k)}} + L_1\right)} = \frac{\lambda}{2\left(L_0\frac{\lambda}{\norms{\nabla f(x^k)}} + L_1 \lambda\right)} \geq \frac{\lambda}{2\left(L_0 + L_1 \lambda\right)}.
    \end{align*}
    Thus, $\eta_k = \eta \leq \frac{\lambda}{2(L_0 + L_1 \lambda)}$.
        
    Using the convexity assumption of the function, we have the following:

    \begin{align*}
        f(x^k) - f^* &\leq \dotprod{\nabla f(x^k)}{x^k - x^*} \overset{\eqref{eq:scalar_product_bound}}{\leq} \norms{\nabla f(x^k)} \norms{x^k - x^*}  \leq \norms{\nabla f(x^k)} \underbrace{\norms{x^0 - x^*}}_{R}.
    \end{align*}
    Hence we have:
    \begin{equation}
        \norms{\nabla f(x^k)} \geq \frac{f(x^k) - f^*}{R}. \label{eq:Biased_NSGD2_case1_2}
    \end{equation}
    
    Then substituting \eqref{eq:Biased_NSGD2_case1_2} into \eqref{eq:Biased_NSGD2_case1_1} we obtain:
    \begin{equation*}
        \expect{f(x^{k+1})} - f(x^k) \leq - \frac{\eta}{2} \norms{\nabla f(x^k)} + \frac{\eta \sigma^2}{2 \lambda B} + \frac{\eta \zeta^2}{2 \lambda} \leq - \frac{\eta}{2 R} (f(x^k) - f^*) + \frac{\eta \sigma^2}{2 \lambda B} + \frac{\eta \zeta^2}{2 \lambda}.
    \end{equation*}
        
    This inequality is equivalent to the trailing inequality:
    \begin{equation*}
        \expect{f(x^{k+1})} - f^* \leq \left( 1 - \frac{\eta}{2R} \right) \left( f(x^k) - f^* \right) + \frac{\eta}{2 \lambda} \left( \frac{\sigma^2}{B} + \zeta^2 \right).
    \end{equation*}
    
    Then for $k = 0,1,2,..., N-1$ iterations that satisfy the conditions $\norms{\nabla f(x^k)} \leq \lambda$ and $\norms{\gg(x^k, \bxi^k)} \geq \lambda $  NSGD with biased gradient oracle shows linear convergence:
    \begin{equation*}
        \expect{f(x^{N})} - f^* \leq \left(1 - \frac{\eta}{2R} \right)^N (f(x^{0}) - f^*) + \frac{R}{\lambda} \left( \frac{\sigma^2}{B} + \zeta^2 \right).
    \end{equation*}

\subsubsection{Third case: \texorpdfstring{$\norms{\nabla f(x^k)} \leq \lambda$}{TEXT} and \texorpdfstring{$\norms{\gg(x^k, \bxi^k)} \leq \lambda$}{TEXT}} 

    Using this in \eqref{eq:Biased_NSGD_smooth}, we have the following with $\eta_k \leq \frac{\norms{\nabla f(x^k)}}{2(L_0 + L_1 \norms{\nabla f(x^k)})}$ and $\alpha = \norms{\nabla f(x^k)}^{-1}$:
    \begin{align}
        \expect{f(x^{k+1})} - f(x^k) &\overset{\eqref{eq:Biased_NSGD_smooth}}{\leq}  - \eta \dotprod{\nabla f(x^k)}{\expect{G(x^k, \bxi^k)}} + \frac{\eta^2 (L_0 + L_1 \norms{\nabla f(x^k)})}{2} \expect{\norms{G(x^k, \bxi^k)}^2} \nonumber\\
        & \overset{\eqref{eq:qudrat_raznosti}}{=} - \frac{\eta \alpha}{2} \norms{\nabla f(x^k)}^2 - \frac{\eta}{2 \alpha} \expect{\norms{G(x^k, \bxi^k)}^2} + \frac{\eta}{2\alpha} \expect{\norms{G(x^k, \bxi^k) - \alpha \nabla f(x^k)}^2}\nonumber\\
        &\quad\quad\quad+ \frac{\eta^2 (L_0 + L_1 \norms{\nabla f(x^k)})}{2} \expect{\norms{G(x^k, \bxi^k)}^2} \nonumber\\
        &= - \frac{\eta}{2} \norms{\nabla f(x^k)} + \frac{\eta}{2\alpha} \expect{\norms{G(x^k, \bxi^k) - \alpha \nabla f(x^k)}^2}\nonumber\\
        &\quad\quad\quad - \frac{\eta}{2} \expect{\norms{G(x^k, \bxi^k)}^2} \left( 1 - \frac{\eta (L_0 + L_1 \norms{\nabla f(x^k)})}{\norms{\nabla f(x^k)}} \right)\nonumber\\
        &\leq - \frac{\eta}{2} \norms{\nabla f(x^k)} + \frac{\eta}{2\alpha} \expect{\norms{G(x^k, \bxi^k) - \alpha \nabla f(x^k)}^2} \nonumber\\
        &\leq - \frac{\eta}{2} \norms{\nabla f(x^k)} + \frac{\eta}{\alpha} \expect{\norms{G(x^k, \bxi^k)}^2 +\norms{\alpha \nabla f(x^k)}^2}\nonumber\\
        &= - \frac{\eta}{2} \norms{\nabla f(x^k)} + \frac{\eta}{\alpha} \expect{\norms{\frac{\gg(x^k, \bxi^k)}{\norms{\gg(x^k, \bxi^k)}}}^2 +\norms{\frac{\nabla f(x^k)}{\norms{\nabla f(x^k)}}}^2}\nonumber\\
        &= - \frac{\eta}{2} \norms{\nabla f(x^k)} + \frac{2 \eta \lambda \norms{\nabla f(x^k)}}{\lambda}\nonumber\\
        &\leq - \frac{\eta}{2} \norms{\nabla f(x^k)} + 2 \eta \lambda. \label{eq:Biased_NSGD3_case1_1}
    \end{align}
    
    The step size will be constant, depending on the hyperparameter $\lambda$:
    \begin{align*}
        \frac{\norms{\nabla f(x^k)}}{2\left(L_0 + L_1 \norms{\nabla f(x^k)}\right)} = \frac{1}{2\left(L_0\frac{1}{\norms{\nabla f(x^k)}} + L_1\right)} = \frac{\lambda}{2\left(L_0\frac{\lambda}{\norms{\nabla f(x^k)}} + L_1 \lambda\right)} \geq \frac{\lambda}{2\left(L_0 + L_1 \lambda\right)}.
    \end{align*}
    Thus, $\eta_k = \eta \leq \frac{\lambda}{2(L_0 + L_1 \lambda)}$.
    
    Using the convexity assumption of the function, we have the following:

    \begin{align*}
        f(x^k) - f^* &\leq \dotprod{\nabla f(x^k)}{x^k - x^*} \overset{\eqref{eq:scalar_product_bound}}{\leq} \norms{\nabla f(x^k)} \norms{x^k - x^*}  \leq \norms{\nabla f(x^k)} \underbrace{\norms{x^0 - x^*}}_{R}.
    \end{align*}
    Hence we have:
    \begin{equation}
        \norms{\nabla f(x^k)} \geq \frac{f(x^k) - f^*}{R}. \label{eq:Biased_NSGD3_case1_2}
    \end{equation}

    Then substituting \eqref{eq:Biased_NSGD3_case1_2} into \eqref{eq:Biased_NSGD3_case1_1} we obtain:
    \begin{equation*}
        \expect{f(x^{k+1})} - f(x^k) \leq - \frac{\eta}{2} \norms{\nabla f(x^k)} + 2 \eta \lambda \leq - \frac{\eta}{2 R} (f(x^k) - f^*) + 2 \eta \lambda.
    \end{equation*}
    
    This inequality is equivalent to the trailing inequality:
    \begin{equation*}
        \expect{f(x^{k+1})} - f^* \leq \left( 1 - \frac{\eta}{2R} \right) \left( f(x^k) - f^* \right) + 2 \eta \lambda.
    \end{equation*}
    
        Then for $k = 0,1,2,..., N-1$ iterations that satisfy the conditions $\norms{\nabla f(x^k)} \leq \lambda$  NSGD with biased gradient oracle shows linear convergence:
        \begin{equation*}
            \expect{f(x^{N})} - f^* \leq \left(1 - \frac{\eta}{2R} \right)^N (f(x^{0}) - f^*) + \lambda R.
        \end{equation*}

\subsubsection{Fourth case: \texorpdfstring{$\norms{\nabla f(x^k)} \geq \lambda$}{TEXT} and \texorpdfstring{$\norms{\gg(x^k, \bxi^k)} \leq \lambda$}{TEXT}}

    Using this in \eqref{eq:Biased_NSGD_smooth}, we have the following with $\eta_k \leq \frac{\norms{\nabla f(x^k)}}{2(L_0 + L_1 \norms{\nabla f(x^k)})}$ and $\alpha = \lambda^{-1}$:
    \begin{align}
        \expect{f(x^{k+1})} - f(x^k) &\overset{\eqref{eq:Biased_NSGD_smooth}}{\leq}  - \eta \dotprod{\nabla f(x^k)}{\expect{G(x^k, \bxi^k)}} \nonumber\\
        & \quad \quad \quad+ \frac{\eta^2 (L_0 + L_1 \norms{\nabla f(x^k)})}{2} \expect{\norms{G(x^k, \bxi^k)}^2} \nonumber\\
        & \overset{\eqref{eq:qudrat_raznosti}}{=} - \frac{\eta \alpha}{2} \norms{\nabla f(x^k)}^2 - \frac{\eta}{2 \alpha} \norms{\expect{G(x^k, \bxi^k)}}^2 \nonumber\\
        & \quad \quad \quad+ \frac{\eta}{2\alpha} \norms{\expect{G(x^k, \bxi^k)} - \alpha \nabla f(x^k)}^2\nonumber\\
        &\quad\quad\quad+ \frac{\eta^2 (L_0 + L_1 \norms{\nabla f(x^k)})}{2} \expect{\norms{G(x^k, \bxi^k)}^2} \nonumber\\
        &= - \frac{\eta}{2 \lambda} \norms{\nabla f(x^k)}^2 + \frac{\eta}{2 \lambda} \norms{ \expect{\lambda G(x^k, \bxi^k)} - \nabla f(x^k)}^2 \nonumber\\
        & \quad \quad \quad+ \frac{\eta^2 (L_0 + L_1 \norms{\nabla f(x^k)})}{2}\nonumber\\
        &= - \frac{\eta}{2 \lambda} \norms{\nabla f(x^k)}^2 + \frac{\eta}{\lambda} \norms{ \expect{\frac{\lambda \gg(x^k, \bxi^k)}{\norms{\gg(x^k, \bxi^k)}} - \gg(x^k,\bxi^k)}}^2 \nonumber\\
        &\quad\quad\quad + \frac{\eta}{\lambda} \norms{\bb(x^k)}^2 + \frac{\eta^2 (L_0 + L_1 \norms{\nabla f(x^k)})}{2}\nonumber \\
        &= - \frac{\eta}{2 \lambda} \norms{\nabla f(x^k)}^2 + \frac{\eta}{2 \lambda} \norms{ \expect{\left(\frac{\lambda}{\norms{\gg(x^k, \bxi^k)}} - 1\right) \gg(x^k,\bxi^k)}}^2 \nonumber\\
        &\quad\quad\quad + \frac{\eta}{\lambda} \norms{\bb(x^k)}^2 + \frac{\eta^2 (L_0 + L_1 \norms{\nabla f(x^k)})}{2}\nonumber \\
        &\leq - \frac{\eta}{2 \lambda} \norms{\nabla f(x^k)}^2 + \frac{\eta}{2 \lambda} \expect{\left(\frac{\lambda}{\norms{\gg(x^k, \bxi^k)}} - 1\right)^2 \norms{ \gg(x^k,\bxi^k)}^2} \nonumber\\
        &\quad\quad\quad + \frac{\eta}{\lambda} \norms{\bb(x^k)}^2 + \frac{\eta^2 (L_0 + L_1 \norms{\nabla f(x^k)})}{2}\nonumber \\
        &\leq - \frac{\eta}{2 \lambda} \norms{\nabla f(x^k)}^2 + \frac{\eta}{2 \lambda} \expect{\frac{\lambda^2}{\norms{\gg(x^k, \bxi^k)}^2} \norms{ \gg(x^k,\bxi^k)}^2} \nonumber\\
        &\quad\quad\quad + \frac{\eta}{\lambda} \norms{\bb(x^k)}^2 + \frac{\eta^2 (L_0 + L_1 \norms{\nabla f(x^k)})}{2}\nonumber \\
        &= - \frac{\eta}{2 \lambda} \norms{\nabla f(x^k)}^2 + \frac{\eta^2 (L_0 + L_1 \norms{\nabla f(x^k)})}{2} + \frac{\eta \lambda}{2}  + \frac{\eta}{\lambda} \norms{\bb(x^k)}^2 \nonumber\\
        &\leq - \frac{\eta}{2} \norms{\nabla f(x^k)} + \frac{\eta^2 (L_0 + L_1 \norms{\nabla f(x^k)})}{2} + \frac{\eta \lambda}{2} + \frac{\eta}{\lambda} \norms{\bb(x^k)}^2 \nonumber\\
        &= - \frac{\eta}{2} \norms{\nabla f(x^k)} \left( 1 - \frac{\eta (L_0 + L_1 \norms{\nabla f(x^k)})}{\norms{\nabla f(x^k)}} \right)  + \frac{\eta \lambda}{2} + \frac{\eta}{\lambda} \norms{\bb(x^k)}^2 \nonumber\\
        &\leq - \frac{\eta}{4} \norms{\nabla f(x^k)} + \frac{\eta \lambda}{2} + \frac{\eta \zeta^2}{\lambda} .
        \label{eq:Biased_NSGD4_case1_1}
    \end{align}
    
    The step size will be constant, depending on the hyperparameter $\lambda$:
    \begin{align*}
        \frac{\norms{\nabla f(x^k)}}{2\left(L_0 + L_1 \norms{\nabla f(x^k)}\right)} = \frac{1}{2\left(L_0\frac{1}{\norms{\nabla f(x^k)}} + L_1\right)} = \frac{\lambda}{2\left(L_0\frac{\lambda}{\norms{\nabla f(x^k)}} + L_1 \lambda\right)} \geq \frac{\lambda}{2\left(L_0 + L_1 \lambda\right)}.
    \end{align*}
    Thus, $\eta_k = \eta \leq \frac{\lambda}{2(L_0 + L_1 \lambda)}$.
    
    Using the convexity assumption of the function, we have the following:

    \begin{align*}
        f(x^k) - f^* &\leq \dotprod{\nabla f(x^k)}{x^k - x^*} \overset{\eqref{eq:scalar_product_bound}}{\leq} \norms{\nabla f(x^k)} \norms{x^k - x^*}  \leq \norms{\nabla f(x^k)} \underbrace{\norms{x^0 - x^*}}_{R}.
    \end{align*}
    Hence we have:
    \begin{equation}
        \norms{\nabla f(x^k)} \geq \frac{f(x^k) - f^*}{R}. \label{eq:Biased_NSGD4_case1_2}
    \end{equation}

    Then substituting \eqref{eq:Biased_NSGD4_case1_2} into \eqref{eq:Biased_NSGD4_case1_1} we obtain:
    \begin{equation*}
        \expect{f(x^{k+1})} - f(x^k) \leq - \frac{\eta}{4} \norms{\nabla f(x^k)} +\frac{\eta \lambda}{2} + \frac{\eta \zeta^2}{\lambda}  \leq - \frac{\eta}{4 R} (f(x^k) - f^*) +\frac{\eta \lambda}{2} + \frac{\eta \zeta^2}{\lambda}.
    \end{equation*}
    
    This inequality is equivalent to the trailing inequality:
    \begin{equation*}
        \expect{f(x^{k+1})} - f^* \leq \left( 1 - \frac{\eta}{4R} \right) \left( f(x^k) - f^* \right) + \frac{\eta \lambda}{2} + \frac{\eta \zeta^2}{\lambda}.
    \end{equation*}
    
        Then for $k = 0,1,2,..., N-1$ iterations that satisfy the conditions $\norms{\nabla f(x^k)} \geq \lambda$ and $\norms{\gg(x^k, \bxi^k)} \leq \lambda$  NSGD with biased gradient oracle shows linear convergence:
        \begin{equation*}
            \expect{f(x^{N})} - f^* \leq \left(1 - \frac{\eta}{4R} \right)^N (f(x^{0}) - f^*) + 2 \lambda R + \frac{2 \zeta^2 R}{\lambda}.
        \end{equation*}

Combining all the cases considered, we obtain the convergence rate of NSGD with biased gradient oracle: 
\begin{equation}
    \expect{f(x^{N})} - f^* \lesssim  \left( 1 - \frac{\eta}{R} \right)^{N} (f(x^0) - f^*) + \frac{M R}{\lambda^2} \left( \frac{\sigma^2}{B} + \zeta^2 \right) + \lambda R.\label{eq:Biased_NSGD}
\end{equation}

\subsection{Convergence Results for ZO-NSGD (Proof of the Theorem~\ref{th:ZO_NSGD})}
In order to obtain the convergence rate of ZO-NSGD in the convex setting, we need to substitute the obtained estimates~\eqref{eq:proof_bias} and \eqref{eq:proof_variance} into the convergence rate of NSGD~\eqref{eq:Biased_NSGD} instead of $\zeta$ and $\sigma^2$, respectively. Then the convergence of ZO-NSGD in the convex setup is as follows:
    \begin{align*}
        \expect{f(x^{N})} - f^* &\lesssim  \underbrace{\left( 1 - \frac{\eta}{R} \right)^{N} (f(x^0) - f^*)}_{\circledOne} + \underbrace{\frac{d M R \tilde{\sigma}^2}{\lambda^2 B}}_{\circledTwo} + \underbrace{\frac{d M R \left(L_0 + L_1 M\right)^2 \gamma^2}{\lambda^2 B}}_{\circledThree} + \underbrace{\frac{d^2 M R \Delta^2}{\lambda^2 B \gamma^2}}_{\circledFour} \\
        &\quad\quad\quad + \underbrace{\frac{M R \left(L_0 + L_1 M\right)^2 \gamma^2}{\lambda^2}}_{\circledFive} + \underbrace{\frac{d^2 M R \Delta^2}{\lambda^2 \gamma^2}}_{\circledSix} + \underbrace{\lambda R}_{\circledSeven}.
    \end{align*}

    \textbf{From term $\circledSeven$}, we find the hyperparameter $\lambda$:
    \begin{align}
        \circledOne: \quad \lambda R \leq \varepsilon \quad & \Rightarrow \quad \lambda \leq \frac{\varepsilon}{R}. \label{eq:ZO_NSGD_lambda}
    \end{align}

    \textbf{From term $\circledOne$}, we find the number of iterations $N$ required for Algorithm \ref{algo:ZO_NSGD} in convex setup to achieve $\varepsilon$-accuracy:
    \begin{align}
        \circledOne: \quad \left( 1 - \frac{\eta}{R} \right)^{N} (f(x^0) - f^*) \leq \varepsilon \quad & \Rightarrow \quad N \geq \frac{R}{\eta} \log \frac{(f(x^0) - f^*)}{\varepsilon};
        \nonumber \\
         N &= \OboundTilde{\frac{R}{\eta}}. \label{eq:ZO_NSGD_iterations}
    \end{align}

    \textbf{From terms $\circledTwo$}, we find the batch size $B$: 
    \begin{align}
        \circledTwo: \quad \frac{ d MR \tilde{\sigma}^2}{\lambda^2 B} \leq \varepsilon \quad & \Rightarrow \quad B \overset{\eqref{eq:ZO_NSGD_lambda}}{\geq} \frac{d M R^3 \tilde{\sigma}^2}{\varepsilon^3};
        \nonumber \\
         B &= \Obound{ \frac{d M R^3 \tilde{\sigma}^2}{\varepsilon^3}}. \label{eq:ZO_NSGD_batch_size}
    \end{align}

    \textbf{From terms $\circledThree$ and $\circledFive$} we find the smoothing parameter $\gamma$:
    \begin{align}
        &\circledThree: \quad \frac{d MR \left(L_0 + L_1 M\right)^2 \gamma^2}{\lambda^2 B} \leq 
        \varepsilon \quad \Rightarrow \quad \gamma \leq \sqrt{\frac{\varepsilon \lambda^2 B}{ d MR \left(L_0 + L_1 M\right)^2}} \overset{\eqref{eq:ZO_NSGD_batch_size},\eqref{eq:ZO_NSGD_lambda}}{=} \frac{\tilde{\sigma}}{\left(L_0 + L_1 M\right)};
        \nonumber \\
        & \circledFive: \quad  \frac{MR \left(L_0 + L_1 M\right)^2 \gamma^2}{\lambda^2} \leq \varepsilon \quad \Rightarrow \quad \gamma \leq \frac{\sqrt{\varepsilon^3}}{\sqrt{M}R^{3/2} \left(L_0 + L_1 M\right)};
        \nonumber \\
        & \quad \quad \quad \gamma \leq \frac{1}{\left(L_0 + L_1 M\right)}\min \left\{\tilde{\sigma}, \frac{\varepsilon^{3/2}}{\sqrt{M}R^{3/2}}\right\} = \frac{\varepsilon^{3/2}}{\left(L_0 + L_1 M\right) \sqrt{M}R^{3/2}}. \label{eq:ZO_NSGD_smoothing_parameter}
    \end{align}

    \textbf{From the remaining terms $\circledFour$ and $\circledSix$}, we find the maximum allowable level of adversarial noise $\Delta$ that still guarantees the convergence of the ZO-NSGD to desired accuracy~$\varepsilon$ in convex setup:
    \begin{align}
        &\circledFour: \quad \frac{d^2 MR \Delta^2}{\lambda^2 B \gamma^2} \leq \varepsilon \quad \Rightarrow \quad \Delta \leq \frac{\sqrt{\varepsilon} \lambda \gamma \sqrt{B}}{d \sqrt{MR}} \overset{ \eqref{eq:ZO_NSGD_batch_size}, \eqref{eq:ZO_NSGD_smoothing_parameter}, \eqref{eq:ZO_NSGD_lambda}}{=}  \frac{\varepsilon^{3/2} \tilde{\sigma}}{\sqrt{d} \left(L_0 + L_1 M\right) R^{3/2}};
        \nonumber \\
        & \circledSix: \quad \frac{d^2 MR \Delta^2}{\gamma^2 \lambda^2} \leq \varepsilon \quad \Rightarrow \quad \Delta \leq \sqrt{\frac{ \gamma^2 \lambda^2 \varepsilon }{d^2 MR}} \overset{ \eqref{eq:ZO_NSGD_lambda},\eqref{eq:ZO_NSGD_smoothing_parameter}}{=} \frac{\varepsilon^3}{d \left(L_0 + L_1 M\right) R^3};
        \nonumber \\
        & \quad \quad \quad \Delta \leq  \frac{\varepsilon^{3/2}}{\sqrt{d} \left(L_0 + L_1 M\right) R^{3/2}} \min \left\{ \tilde{\sigma},\frac{\varepsilon^{3/2}}{\sqrt{d} R^{3/2}}\right\}. 
        \label{eq:ZO_NSGD_noise_level}
    \end{align}

    In this way, the ZO-NSGD achieves $\varepsilon$-accuracy: $\expect{f(x^{N}) - f^*} \leq \varepsilon$ in convex setup after 
        \begin{equation*}
            N \overset{\eqref{eq:ZO_NSGD_iterations}}{=} \OboundTilde{\frac{R}{\eta}}, \quad T = N \cdot B \overset{\eqref{eq:ZO_NSGD_iterations}, \eqref{eq:ZO_NSGD_batch_size}}{=} \Obound{\frac{d\tilde{\sigma}^2 M R^4}{\varepsilon^3 \eta}}
        \end{equation*}
        number of iterations, total number of zero-order oracle calls and at
        \begin{equation*}
            \Delta \overset{\eqref{eq:ZO_NSGD_noise_level}}{\lesssim} \frac{\varepsilon^{3/2}}{\sqrt{d} \left(L_0 + L_1 M\right) R^{3/2}} \min \left\{ \tilde{\sigma},\frac{\varepsilon^{3/2}}{\sqrt{d} R^{3/2}}\right\}
        \end{equation*}
        the maximum level of noise with smoothing parameter $\frac{\varepsilon^{3/2}}{\left(L_0 + L_1 M\right) \sqrt{M}R^{3/2}}$ \eqref{eq:ZO_NSGD_smoothing_parameter}.

\section{Additional Clarification}\label{app:Additional Clarification}
In this section, we would like to clarify the convergence in the case $L_0 = 0$ (Remark~\ref{rem:Smoothness_only_L0}). In this case the problem does not reach a minimum (hence $R = \arginf f(x) = + \infty$). Therefore, we exemplify the special case of NSGD (when $\norms{\nabla f(x^k,\bxi^k)} \geq \sqrt{2} \sigma$ and $\norms{\nabla f(x^k)} \geq \sqrt{2} \sigma$), shows that it is possible to achieve the desired accuracy $\varepsilon$ in a finite number of iterations.

Let's introduce the notation $G(x^k, \bxi^k) = \frac{\nabla f(x^k, \bxi^k)}{\norms{\nabla f(x^k, \bxi^k)}}$, then using $(L_0,L_1)$-smoothness (see Assumption~\ref{ass:L0_L1_smooth}):
\begin{align}
    f(x^{k+1}) - f(x^k) &\overset{\eqref{eq:ass_smooth}}{\leq} \dotprod{\nabla f(x^k)}{x^{k+1} - x^k} + \frac{L_0 + L_1 \norms{\nabla f(x^k)}}{2} \norms{x^{k+1} - x^{k}}^2 \nonumber \\
    &= - \eta \dotprod{\nabla f(x^k)}{G(x^k, \bxi^k)} + \frac{\eta^2 (L_0 + L_1 \norms{\nabla f(x^k)})}{2} \norms{G(x^k, \bxi^k)}^2. \label{eq:NSGD_smooth_additional}
\end{align}

    Let us evaluate first summand of \eqref{eq:NSGD_smooth_additional} with $\alpha = \norms{\nabla f(x^k)}^{-1}$:
    \begin{align}
        - \eta \dotprod{\nabla f(x^k)}{G(x^k, \bxi^k)} &\overset{\eqref{eq:qudrat_raznosti}}{=} - \frac{\alpha \eta }{2} \norms{\nabla f(x^k)}^2 - \frac{\eta}{2 \alpha} \norms{G(x^k, \bxi^k)}^2 \nonumber\\
        &\quad \quad \quad + \frac{\eta}{2 \alpha} \norms{G(x^k, \bxi^k) - \alpha \nabla f(x^k)}^2 \nonumber\\
        &= - \frac{\eta}{2 } \norms{\nabla f(x^k)} - \frac{\eta}{2 \alpha} \norms{G(x^k, \bxi^k)}^2 \nonumber\\
        &\quad \quad \quad+ \frac{\eta}{2 \lambda^2 \alpha} \norms{\lambda G(x^k, \bxi^k) - \lambda \alpha \nabla f(x^k)}^2 \nonumber\\
        &= - \frac{\eta}{2 } \norms{\nabla f(x^k)} - \frac{\eta}{2 \alpha} \norms{G(x^k, \bxi^k)}^2 \nonumber\\
        &\quad \quad \quad+ \frac{\eta}{2 \lambda^2 \alpha} \norms{\cliplam{\nabla f(x^k, \bxi^k)} - \cliplam{ \nabla f(x^k)}}^2 \nonumber
    \end{align}
    Using that clipping is a projection on onto a convex set, namely ball with radius $\lambda$, and thus is Lipshitz operator with Lipshitz constant $1$, we can obtain:
    \begin{align}
        - \eta \dotprod{\nabla f(x^k)}{\expect{G(x^k, \bxi^k)}} &\leq - \frac{\eta}{2} \norms{\nabla f(x^k)} - \frac{\eta}{2 \alpha} \expect{\norms{G(x^k, \bxi^k)}^2} \nonumber\\
        &\quad \quad \quad+ \frac{\eta }{2 \lambda^2 \alpha} \expect{\norms{\nabla f(x^k, \bxi^k) - \nabla f(x^k)}^2}.\label{eq:NSGD1_main_additional}
    \end{align}
    
     Using this in \eqref{eq:NSGD1_main_additional}, we have the following with $\eta_k \leq \frac{\norms{\nabla f(x^k)}}{2(L_0 + L_1 \norms{\nabla f(x^k)})}$:
        \begin{align}
            \expect{f(x^{k+1})} - f(x^k) &\overset{\eqref{eq:NSGD_smooth_additional}}{\leq}  - \eta \dotprod{\nabla f(x^k)}{\expect{G(x^k, \bxi^k)}} + \frac{\eta^2 (L_0 + L_1 \norms{\nabla f(x^k)})}{2} \expect{\norms{G(x^k, \bxi^k)}^2} \nonumber\\
            & \overset{\eqref{eq:NSGD1_main_additional}}{\leq} - \frac{\eta}{2} \norms{\nabla f(x^k)} - \frac{\eta}{2 \alpha} \expect{\norms{G(x^k, \bxi^k)}^2} + \frac{\eta }{2 \lambda^2 \alpha} \expect{\norms{\nabla f(x^k, \bxi) - \nabla f(x^k)}^2}\nonumber\\
            &\quad\quad\quad+ \frac{\eta^2 (L_0 + L_1 \norms{\nabla f(x^k)})}{2} \expect{\norms{G(x^k, \bxi^k)}^2} \nonumber\\
            &= - \frac{\eta}{2} \norms{\nabla f(x^k)} + \frac{\eta }{2 \lambda^2 \alpha} \expect{\norms{\nabla f(x^k, \bxi^k) - \nabla f(x^k)}^2}\nonumber\\
            &\quad\quad\quad - \frac{\eta}{2} \expect{\norms{G(x^k, \bxi^k)}^2} \left( 1 - \frac{\eta (L_0 + L_1 \norms{\nabla f(x^k)})}{\norms{\nabla f(x^k)}} \right)\nonumber\\
            &\leq - \frac{\eta}{2} \norms{\nabla f(x^k)} + \frac{\eta \sigma^2}{2 \lambda^2 \alpha} \nonumber\\
            &\leq - \frac{\eta}{2} \norms{\nabla f(x^k)} + \frac{\eta}{4} \norms{\nabla f(x^k)}\nonumber\\
            &= - \frac{\eta}{4} \norms{\nabla f(x^k)}. \label{eq:NSGD1_case1_1_additional}
        \end{align}
        
        The step size will be constant, depending on the hyperparameter $\lambda$:
        \begin{align*}
            \frac{\norms{\nabla f(x^k)}}{2\left(L_0 + L_1 \norms{\nabla f(x^k)}\right)} = \frac{1}{2\left(L_0\frac{1}{\norms{\nabla f(x^k)}} + L_1\right)} = \frac{\lambda}{2\left(L_0\frac{\lambda}{\norms{\nabla f(x^k)}} + L_1 \lambda\right)} \geq \frac{\lambda}{2\left(L_0 + L_1 \lambda\right)}.
        \end{align*}
        Thus, $\eta_k = \eta \leq \frac{\lambda}{2(L_0 + L_1 \lambda)}$.

        We introduce the hyperparameter of the algorithm $R_s = \norms{x^0 - s}$. Then using the convexity assumption of the function, we have the following:
        
           \begin{align*}
                f(x^k) - f(s) &\leq \dotprod{\nabla f(x^k)}{x^k - s} \\
                &\overset{\eqref{eq:scalar_product_bound}}{\leq} \norms{\nabla f(x^k)} \norms{x^k - s}\\
                &\leq \norms{\nabla f(x^k)} \underbrace{\norms{x^0 - s}}_{R_s}.
            \end{align*}
            Hence we have:
            \begin{equation}
                \label{eq:NSGD1_case1_2_additional}
                \norms{\nabla f(x^k)} \geq \frac{f(x^k) - f(s)}{R_s}.
            \end{equation}

            Then substituting \eqref{eq:NSGD1_case1_2_additional} into \eqref{eq:NSGD1_case1_1_additional} we obtain:
            \begin{equation*}
                \expect{f(x^{k+1})} - f(x^k) \leq - \frac{\eta}{4} \norms{\nabla f(x^k)} \leq - \frac{\eta}{4R_s} (f(x^k) - f(s)).
            \end{equation*}
            
            This inequality is equivalent to the trailing inequality:
            \begin{equation*}
                \expect{f(x^{k+1})} - f^* \leq \left( 1 - \frac{\eta}{4R_s} \right) \left( f(x^k) - f^* \right) + \frac{\eta}{4R_s}(f(s) - f^*).
            \end{equation*}
            
                Then for $k = 0,1,2,..., N-1$ iterations that satisfy the conditions $\norms{\nabla f(x^k,\bxi^k)} \geq \sqrt{2} \sigma$ and $\norms{\nabla f(x^k)} \geq \sqrt{2} \sigma$ NSGD shows linear convergence:
                \begin{equation*}
    f(x^{N}) - f^*  \leq \left( 1 - \frac{\eta}{4R_s}\right)^N (f(x^0) - f^*) + f(s) - f^*.
\end{equation*}

Thus, we have shown that it is indeed possible to converge to a linear rate of convergence on logistic regression using the hyperparameter $R_s$.



\end{document}